\documentclass[11pt]{amsart}

\usepackage{amsmath,amsthm, amscd, amssymb, amsfonts}
\usepackage[all]{xy}
\usepackage[dvips, dvipsnames, usenames]{color}

\newtheorem{theorem}{Theorem}[section]

\newtheorem{lemma}[theorem]{Lemma}
\newtheorem{corollary}[theorem]{Corollary}
\newtheorem{conjecture}{Conjecture}
\newtheorem{teo intro}{Theorem}
\newtheorem{question}{Question}
\newtheorem{proposition}[theorem]{Proposition}

\theoremstyle{definition}
\newtheorem{definition}[theorem]{Definition}

\theoremstyle{remark}
\newtheorem{remark}[theorem]{Remark}

\def\yd{^{H}_{H}\mathcal{YD}}

\def\zt{\Z^{\theta}}

\def\G{\mathbb{G}}
\def\xx{\mathbb{X}}
\def\SS{\mathcal{S}}
\newcommand{\g}{\mathfrak g}
\newcommand\id{\operatorname{id}}

\newcommand\ord{\operatorname{ord}}
\newcommand\Hom{\operatorname{Hom}}

\newcommand\Aut{\operatorname{Aut}}

\newcommand\Sh{{\operatorname{Sh}}}
\newcommand\cop{\operatorname{cop}}

\newcommand\gr{\operatorname{gr}}

\newcommand\op{\operatorname{op}}
\newcommand\ad{\operatorname{ad}}
\newcommand\sop{\operatorname{supp}}

\def\gyd{{}^{G}_{G}\mathcal{YD}}
\newcommand{\ku}{ \mathbf{k}}

\def\ot{\otimes}

\def\bq{\mathfrak{q}}

\def\Z{\mathbb{Z}}
\def\N{\mathbb{N}}

\def\cB{\mathcal{B}}
\def\cO{\mathcal{O}}
\def\cU{\mathcal{U}}
\def\cC{\mathcal{C}}
\def\cI{\mathcal{I}}
\def\cJ{\mathcal{J}}
\def\cR{\mathcal{R}}

\def\cH{\mathcal{H}}
\def\cW{\mathcal{W}}
\def\cP{\mathcal{P}}
\def\u{\mathfrak{u}}
\def\Xf{\mathfrak{X}}

\newcommand{\Eb}{\underline E}
\newcommand{\Fb}{\underline F}
\newcommand{\Kb}{\underline K}
\newcommand{\Lb}{\underline L}
\newcommand{\qb}{\underline q}
\newcommand{\cb}{\underline c}
\def\xb{\mathbf{x}}
\def\yb{\mathbf{y}}
\def\qmb{\mathbf{q}}
\newcommand{\schi}{{s_p^*\chi}}

\newcommand{\eps}{\varepsilon}

\newcommand{\vi}{${\mathsf {(i)}\;}$}
\newcommand{\vii}{${\mathsf {(ii)}\;}$}
\newcommand{\viii}{${\mathsf {(iii)}\;}$}
\newcommand{\viv}{${\mathsf {(iv)}\;}$}
\newcommand{\vv}{${\mathsf {(v)}\;}$}
\newcommand{\vvi}{${\mathsf {(vi)}\;}$}
\newcommand{\vvii}{${\mathsf {(vii)}\;}$}

\def\pf{\begin{proof}}
\def\epf{\end{proof}}

\begin{document}


\title[On Nichols algebras of diagonal type]{On Nichols algebras of diagonal type}
\author[Iv\'an Angiono]{Iv\'an Angiono}

\address{FaMAF-CIEM (CONICET), Universidad Nacional de C\'ordoba,
Medina A\-llen\-de s/n, Ciudad Universitaria, 5000 C\' ordoba, Rep\'
ublica Argentina.} \email{angiono@famaf.unc.edu.ar}

\thanks{The work was partially supported by CONICET,
FONCyT-ANPCyT, Secyt (UNC), Mincyt (C\'ordoba)}

\begin{abstract}
We give an explicit and essentially minimal list of defining relations of a Nichols algebra of diagonal type with finite root system. This list contains
the well-known quantum Serre relations but also many new variations. A conjecture by Andruskiewitsch and Schneider states that any finite-dimensional
pointed Hopf algebra over an algebraically closed field of characteristic zero is generated as an algebra by its group-like and skew-primitive elements. As
an application of our main result, we prove the conjecture when the group of group-like elements is abelian.
\end{abstract}

\maketitle

\section*{Introduction}

\subsection*{1}
Let $\ku$ be an algebraically closed field of characteristic zero
and let $\theta$ be a natural number. Let $\qmb=(q_{ij})_{1\leq
i,j\leq\theta}$ be a matrix with invertible entries on $\ku$ and let
$V$ be a vector space of dimension $\theta$. The Nichols algebra
associated to $\qmb$ is a graded connected Hopf algebra
$\cB(V)=\oplus_{n\geq0}\cB^n(V)$ with many favorable properties. It
plays a fundamental role in the classification of finite-dimensional
(or finite growth) pointed Hopf algebras. Precisely, a basic
question in the classification Program \cite{AS lifting} is the
following:
\begin{question}\label{question}
\cite[Question 5.9]{An}: Given $(V,\qmb)$, determine if the associated Nichols algebra $\cB(V)$ is finite-dimensional. In such case, compute the dimension
of $\cB(V)$ and give a presentation by generators and relations.
\end{question}
The first part of this question has been answered by Heckenberger \cite{H-classif RS}, who obtained the list of all matrices $\qmb$  whose associated
Nichols
algebra has a finite root system. Roughly, this list contains three classes of matrices:
\begin{itemize}
 \item Standard matrices \cite{AA}: they are associated with finite Weyl groups. Their root systems coincide with root systems of finite Cartan matrices.
This family includes properly the so-called braidings of Cartan type, in particular the matrices related with the positive part of the small quantum groups.
 \item Matrices of super type \cite{AAY}, related with the positive part of quantized enveloping algebras of contragradient Lie superalgebras. Their root
systems become from the corresponding Lie superalgebras.
 \item A finite list of exceptional matrices, whose associated diagram has connected components with at most 7 vertices, and the scalars defining these
braidings are roots of unity of low order.
\end{itemize}
There are several answers to the second part of Question \ref{question} under particular assumptions:
\begin{itemize}
 \item[$\rhd$] \cite{L} for the positive part of quantized enveloping algebras of semisimple Lie algebras and small quantum groups, using the full representation
theory of quantum groups;
 \item[$\rhd$] \cite{AS Fin-qg} for braidings of Cartan type;
 \item[$\rhd$] \cite{A-standard} for braidings of standard type;
 \item[$\rhd$] \cite{Y} for the positive part of quantized enveloping algebras of contragradient Lie superalgebras;
 \item[$\rhd$] \cite{AAY} for braidings of super type;
 \item[$\rhd$] \cite{H-rk two}, giving a general form of relations for matrices of rank two;
 \item[$\rhd$] \cite{He} for some examples of rank two matrices, but giving explicit relations.
\end{itemize}
In \cite{A-presentation} we gave general formulae for defining relations of Nichols algebras of diagonal type, see Theorem \ref{Thm:presentacion} below.
The expression of those relations and the proof that they generate the defining ideal are independent of Heckenberger's classification; they  rely in
Kharchenko's and Rosso's PBW bases \cite{Kh,R-lyndon words} and a detailed study of convex orders in generalized root systems \cite{A-presentation}, through
the classification of coideal subalgebras \cite{HS}. In this paper we refine the main result of \cite{A-presentation} and prove:
\begin{teo intro}\label{thm: presentacion intro}
A minimal set of relations of $\cB(V)$ is obtained by considering relations of the following type:
\begin{enumerate}
  \item Quantum Serre relations, and powers of generators $x_i$ corresponding to non-Cartan vertices; they are needed to introduce Lusztig's isomorphisms
at the level of doubles of tensor algebras.
  \item Relations in the image of the previous ones by the Lusztig's isomorphisms, and correspond to relations \eqref{quantumSerregeneralizadas1} in
Theorem \ref{Thm:presentacion}.
  \item Relations that guarantee that the ideal generated by the previous relations is a braided biideal: they appear in the coproduct of relations of the
item $(2)$ in the tensor algebra $T(V)$.
  \item Powers of root vectors (generators of the PBW basis) corresponding to roots in the orbit of Cartan vertices.
\end{enumerate}
\end{teo intro}
See Theorem \ref{thm:presentacion minima} for a complete and explicit set of relations. In this set we distinguish relations appearing in
\cite{A-presentation} for standard braidings, and relations in \cite{Y} related with braidings of super type. There exists also a large list of new
relations, related with the set of exceptional braidings or with braidings of super type evaluated in roots of unity of small order.

\subsection*{2} The knowledge of the explicit relations of a Nichols algebra has several potential applications to the theory of pointed Hopf algebras,
that we discuss now:
\begin{itemize}
 \item One of the basic question in the Lifting Method \cite{AS lifting,AS Pointed HA} for the classification of Hopf algebras is the following:
\end{itemize}
\begin{conjecture}{\cite[Conjecture 1.4]{AS Fin-qg}}\label{conj:AS, generacion en grado 1}
Let $\Gamma$ be a finite group and $\ku$ an algebraically closed field of characteristic 0. If $H$ is a finite-dimensional pointed Hopf algebra over $\ku$
such that $G(H)=\Gamma$, then $H$ is generated as an algebra by $\Gamma$ and its skew-primitive elements.
\end{conjecture}
This question was answered in \cite{AS Class} for braidings of Cartan type under some mild conditions. This result was extended to the case of
standard braidings in \cite{AnGa}. In Section \ref{section:conjetura AS} we obtain as a consequence of Theorem \ref{thm:presentacion minima}:
\begin{teo intro}\label{thm: generacion grado 1 intro}
Let $H$ be a finite dimensional pointed Hopf algebra over an abelian group $\Gamma$. Then $H$ is generated as an algebra by $\Gamma$ and its skew primitive
elements.
\end{teo intro}
That is, we answer positively Conjecture \ref{conj:AS, generacion en grado 1} in a general context: when $G(H)$ is any abelian group. This Theorem is also
applied to the known cases of finite-dimensional Nichols algebras over non-abelian groups

\begin{itemize}
\item Another crucial step of the Lifting Method is to obtain all deformations of the pointed Hopf algebras $\cB(V)\#\ku\Gamma$; that is, all the pointed
Hopf algebras such that their associated coradically graded algebras are isomorphic to $\cB(V)\#\ku\Gamma$.
\end{itemize}
This problem was solved for $\Gamma$ abelian in \cite{AS Class} --
under the restriction that the order is not divisible by 2,3,5,7. We
believe that the explicit presentation in this paper would be
substantial to solve the question for any abelian group.
\begin{itemize}
\item The explicit relations would be also useful in the study of various elements in the representation theory of pointed Hopf algebras. In this direction,
the theory of Nichols algebras of diagonal type provides an uniform approach to the study of quantum groups and quantum supergroups.
\end{itemize}

\subsection*{3}
The plan of this paper is the following. We introduce the notion of Nichols algebras in Section \ref{section:preliminares}. We give a PBW basis of any
Nichols algebra and some properties of this basis following \cite{Kh,R-lyndon words}. Next we recall the notions of Weyl groupoid and its associated root
system following \cite{CH2,HS}, and make a connection with the theory of Nichols algebras of diagonal type. We present the needed material from
\cite{A-presentation}, in particular Theorem \ref{Thm:presentacion}, a key result for our purposes.

Section \ref{section:Lusztig isomorphisms} is devoted to Lusztig's isomorphisms in the general context of braidings of diagonal type \cite{H-isom},
extending analogous isomorphisms from \cite{L}.

In Section \ref{section:presentacion explicita} we give the
mentioned presentation by generators and relations, based in the
classification of braidings of diagonal type with finite root system
\cite{H-classif RS}. The strategy of proof consists first to define
Lusztig isomorphisms for the Drinfeld doubles of the braided Hopf
algebras $U^+$ obtained by quotient by the relations in Theorem
\ref{thm: presentacion intro}, except the group in $(4)$. This
quotient is analogous to the algebra $U_q^+(\g)$; the Drinfeld
double $\u_q(\g)$ of the Nichols algebra is a quotient of the
previous algebra, as it was considered by Lusztig and
Andruskiewitsch-Schneider. We denote these two algebras by $U^+$ and
$\u^+$, respectively, so $\u^+=\cB(V)$. The existence of the
Lusztig's isomorphisms prove that the PBW generators corresponding
to the algebras $U^+$ and their quotients $\u^+$ are the same, but
the heights of some generators are not the order of the associated
scalar in $U^+$. Therefore we obtain $\u^+$ after to quotient $U^+$
by some powers of root vectors as in $(4)$.

Theorem \ref{thm:presentacion minima} extends the presentation obtained in \cite{A-standard} for standard braidings, and in \cite{AAY} for braidings of
super type, and gives a new proof in the case of braidings of Cartan type, in particular quantized enveloping algebras $U_q(\mathfrak{g})$ and small
quantum groups $\mathfrak{u}_q(\mathfrak{g})$.

Finally, Section \ref{section:conjetura AS} is devoted to the proof of Theorem \ref{thm: generacion grado 1 intro}. We prove first that any finite
dimensional braided graded Hopf algebra of diagonal type
$$S=\oplus_{n\geq 0}S_n, \qquad S_0=\ku 1, \ S_1\cong V,$$
generated as an algebra by $V$ is isomorphic to the Nichols algebra $\cB(V)$; this result extends \cite[Thm. 5.5]{AS Class}, \cite[Thm. 2.5]{AnGa}, but the
proof follows the same scheme.

\subsection*{Acknowledges}
This work is part of the author's PhD Thesis. I want to thank specially to my advisor Nicol\'as Andruskiewitsch for his inspiring guidance, patience
and supervision during these years. I want to thank also to my family for all their support, and to Antonela for all her love.

\section{Preliminaries}\label{section:preliminares}

In this Section we recall results from different works needed in the sequel. First we consider the existence of PBW bases for Nichols algebras of diagonal
type \cite{Kh,R-lyndon words}, and the rich combinatoric related to them. Next we recall the definitions of Weyl groupoid, the
associated root systems and some properties thereof \cite{HS,HY}. We close this Section stating a general presentation of Nichols algebras coming from
\cite{A-presentation}.

\subsection{Lyndon words and PBW bases for Nichols algebras of diagonal type}\label{subsection:pbw}
\
To begin with, we recall the definition of a Nichols algebras and show a characterization in the case of a diagonal braiding.

\begin{definition}{\cite{AS Pointed HA}} \label{def:algebra de Nichols}
Given $V\in \yd$, the tensor algebra $T(V)$ admits a unique
structure of graded braided Hopf algebra in $\yd$ such that $V
\subseteq \cP(V)$. Consider the family $\mathfrak{S}$ of all the
homogeneous Hopf ideals $I\subseteq T(V)$ such that
\begin{itemize}
    \item $I$ is generated by homogeneous elements of degree $\geq 2$,
    \item $I$ is a Yetter-Drinfeld submodule of $T(V)$.
\end{itemize}
The \emph{Nichols algebra} $\cB(V)$ associated to $V$ is the quotient of $T(V)$ by the biggest ideal $I(V)$ of $\mathfrak{S}$.
\end{definition}

\bigbreak \emph{Let $(V,c)$ be a braided vector space of diagonal type such that $q_{ij}=q_{ji}$ for any $i,j$.} Let $\Gamma= \zt$, and $\alpha_1,\ldots,
\alpha_{\theta}$ be the canonical basis. We set the characters $\chi_1,\ldots,\chi_{\theta}$ of $\Gamma$ given by $ \chi_j(\alpha_i)=q_{ij}$, $1 \leq i,j
\leq \theta$.

Consider $V$ as a Yetter-Drinfeld module over $\ku \Gamma$ such that
$x_i \in V_{\alpha_i}^{\chi_i}$. In this context we can characterize
the Nichols algebra as a quotient that admits a certain
non-degenerate bilinear form. We use Swedler's notation for the
coproduct in $T(V)$: $\Delta(x)= x_1\otimes x_2$, where we omit the
summation symbol.

\begin{proposition}{\cite[Prop. 1.2.3]{L}, \cite[Prop. 2.10]{AS Pointed HA}}\label{prop:formabilineal}
There exists a unique bilinear form $(\cdot|\cdot): T(V) \times T(V) \rightarrow \ku$ such that
$(1|1)=1$, and:
\begin{eqnarray}
    (x_i | x_j) &=& \delta_{ij}, \quad \mbox{for any } i,j; \label{eq:bilinearprop1}
    \\ (x|yy') &=& (x_{1} | y) (x_{2} | y'), \quad \mbox{for any } x,y,y' \in
    T(V);\label{eq:bilinearprop2}
    \\ (xx'|y) &=& (x|y_{1}) (x'|y_{2}), \quad \mbox{for any } x,x',y \in T(V). \label{eq:bilinearprop3}
\end{eqnarray}
This is a symmetric form, for which we have:
\begin{equation}
    (x|y)=0, \quad \mbox{for any } x \in T(V)_g, \ y \in T(V)_h, \ g,h \in \Gamma, \ g \neq h. \label{eq:bilinearprop4}
\end{equation}
The radical of this form $\left\{ x \in T(V):(x|y)=0,\ \forall y\in T(V)\right\}$ coincides with $I(V)$, so $(\cdot|\cdot)$ induces a non-degenerate
bilinear form on $\cB(V)=T(V)/I(V)$, denoted also by $(\cdot|\cdot)$. \qed
\end{proposition}

Therefore $I(V)$ is a $\zt$-homogeneous ideal, and then $\cB(V)$ is $\zt$-graded.

\bigbreak

Let $A$ be an algebra, $P,S \subset A$ and $h: S \mapsto \N \cup\{ \infty \}$. We fix a linear order $<$ on $S$. $B(P,S,<,h)$ will denote the set
\begin{align*}
\big\{ &p\,s_1^{e_1}\dots s_t^{e_t}: t \in \N_0, \quad s_1>\dots >s_t,\quad s_i \in S, \quad 0<e_i<h(s_i), \quad p \in P \big\}.
\end{align*}
If $B(P,S,<,h)$ is a $\ku$-linear basis, we say that $(P,S,<,h)$ is a set of \emph{PBW generators}, whose \emph{height} is $h$, and $B(P,S,<,h)$ is a
\emph{PBW basis} of $A$.

\bigbreak We will describe a particular PBW basis for any graded braided Hopf algebra $\cB=\oplus_{n\in \N} \cB^n$ generated by $\cB^1 \cong V$ as an
algebra, where $V$ is a braided vector space; we will follow the results in \cite{Kh}.

Fix $\theta\in \N$, and a set $X= \{x_1,\dots, x_{\theta}\}$. Let $\xx$ be the set of words with letters in $X$ and consider the lexicographical order
on $\xx$.

\begin{definition}
An element $u \in \xx$, $u\neq 1$ is a \emph{Lyndon word} if for any decomposition $u=vw$, $v,w \in\xx - \left\{ 1 \right\}$, we have $u<w$.
We will denote the set of all Lyndon words by $L$.
\end{definition}

\begin{remark}
\begin{itemize}
 \item Each Lyndon word begin with its smallest letter.
 \item Each $u \in \xx-X$ is a Lyndon word if and only if for each decomposition $u=u_1 u_2$ with $u_1,u_2 \in \xx\setminus{1}$,
we have $u_1u_2=u < u_2u_1$.
 \item If $u_1,u_2 \in L$ and $u_1<u_2$, then $u_1u_2 \in L$.
\end{itemize}
\end{remark}

\bigbreak A basic Lyndon's result says that any word $u \in \xx$ admits a unique decomposition as non-increasing product of Lyndon words:
\begin{equation}\label{eq:descly}
u=l_1l_2\dots  l_r, \qquad l_i \in L, l_r \leq \dots \leq l_1.
\end{equation}
It is called the \emph{Lyndon decomposition} of $u \in \xx$, and the $l_i\in L$ in \eqref{eq:descly} are called the \emph{Lyndon letters} of $u$.

Another characterization of Lyndon words is the following:

\begin{lemma}{\cite[p.6]{Kh}}
Let $u \in\xx-X$. Then $u \in L$ if and only if there exist $u_1,u_2 \in L$ such that $u_1<u_2$ and $u=u_1u_2$.\qed
\end{lemma}

\begin{definition}
For each $u \in L-X$, the \emph{Shirshov decomposition} of $u$ is the decomposition $u=u_1u_2$, $u_1,u_2 \in L$, such that $u_2$ is the smallest end
of $u$ between all the possible decompositions with these conditions.
\end{definition}

Given a finite-dimensional vector space $V$, fix a basis $X =
\{x_1,\dots ,x_{\theta}\}$ of $V$; we can identify $\ku \xx$ with
$T(V)$. In what follows we consider two graduations for the algebra
$T(V)$: the usual $\N_0$-graduation $T(V) = \oplus_{n\geq 0}T^n(V)$,
and $\zt$-graduation of $T(V)$, determined by the condition $\deg
x_i = \alpha_i$, $1\le i \le \theta$, where
$\{\alpha_1,\dots,\alpha_\theta \}$ is the canonical basis of $\zt$.
\medskip

Consider a braiding $c$ for $V$. The \emph{braided bracket} of $x,y\in T(V)$ is defined by
\begin{equation}\label{eq:braidedcommutator}
[x,y]_c := \text{multiplication } \circ \left( \id - c \right) \left( x \ot y \right).
\end{equation}

Assume that $(V,c)$ is of diagonal type, and let $\chi: \zt\times
\zt \to \ku^{\times}$ be the bicharacter determined by the condition
\begin{equation}\label{eq:forma bilineal diagonal}
\chi(\alpha_i, \alpha_j) = q_{ij}, \quad \mbox{ for each pair }1\le i, j \le \theta.
\end{equation}
Then, for each pair of $\zt$-homogeneous elements $u,v \in \xx$,
\begin{equation}\label{eq:braiding tipo diagonal}
    c(u \ot v)= q_{u,v} v \ot u, \qquad q_{u,v} = \chi(\deg u, \deg v)\in \ku^{\times}.
\end{equation}
In such case, the braided bracket satisfies a ``braided Jacobi
identity'' and determines skew-derivations as follows:
\begin{align}\label{eq:identidad jacobi}
\left[\left[ u, v \right]_c, w \right]_c &= \left[u, \left[ v, w\right]_c \right]_c-\chi(\alpha,\beta) v \ \left[u,w\right]_c + \chi(\beta,\gamma)
\left[u,w\right]_c \ v,
\\ \label{eq:derivacion} \left[ u,v \ w \right]_c &= \left[ u,v \right]_c w + \chi( \alpha, \beta ) v \ \left[ u,w \right]_c,
\\ \label{eq:derivacion 2} \left[ u \ v, w \right]_c &= \chi( \beta, \gamma ) \left[ u,w \right]_c \ v + u \ \left[ v,w \right]_c,
\end{align}
where $u,v,w \in T(V)$ are homogeneous of degree $\alpha,\beta,\gamma \in \N^{\theta}$, respectively.

Using the previous decompositions, we can define the $\ku$-linear
endomorphism $\left[ - \right]_c$ of $\ku \xx$ as follows:
$$ \left[ u \right]_c := \begin{cases} u,& \text{if } u = 1 \text{ or }u \in X;\\
[\left[ v \right]_c, \left[ w \right]_c]_c,  & \text{if } u \in L, \, \ell(u)>1, \ u=vw \text{ is the Shirshov decomposition};\\
\left[u_1\right]_c \dots \left[u_t\right]_c,& \text{ if } u\in \xx-L \text{ and its Lyndon decomposition is }u=u_1\dots u_t.
\end{cases}
$$

\begin{definition} The \emph{hyperletter} corresponding to $l \in L$ is $\left[l\right]_c$. An \emph{hyperword} is a word whose letters are hyperletters,
and a \emph{monotone hyperword} is an hyperword $W=\left[u_1\right]_c^{k_1}\dots\left[u_m\right]_c^{k_m}$ such that $u_1>\dots >u_m$.
\end{definition}

\begin{remark}\label{corchete}
For any $u \in L$, $\left[ u \right]_c$ is a $\Z\left[q_{ij}\right]$-linear combination of words with the same $\zt$-graduation than $u$, such that
$\left[u\right]_c\in u+ \ku \xx^{\ell(u)}_{>u}$.
\end{remark}

\begin{theorem}{\cite[Thm. 10]{R-lyndon words}} \label{thm:corchete hiperletras}
Let $u,v \in L$, $u<v$. Then $\left[\left[u\right]_c,\left[v\right]_c \right]_c$ is a $\Z\left[q_{ij}\right]$-linear combination of monotone hyperwords
$\left[l_1\right]_c \dots  \left[l_r\right]_c$, $l_i \in L$, such that the corresponding hyperletters satisfy $v>l_i \geq uv$. Moreover, $\left[uv\right]_c$
appears in such combination with non-zero coefficient and each hyperword has the same $\zt$-graduation than $uv$. \qed
 \end{theorem}

We consider the polynomials in $\Z[t]$:
$$ (n)_t:=1+t+\cdots+t^{n-1},\quad (n)_t!=(1)_t (2)_t\cdots (n)_t, \qquad n\in \N. $$
The q-combinatorial numbers are defined as the following quotient:
$$ \binom{n}{i}_t:=\frac{(n)_t!}{(n-i)_t!(i)_t!}, \qquad 0 \leq i \leq n. $$
It follows inductively that $\binom{n}{i}_t\in\Z[t]$, for all $n$
and all $i\in\{0,1,\ldots,n\}$. Therefore $\binom{n}{i}_q$ will
denote the evaluation of $\binom{n}{i}_t$ for $t=q$, where
$q\in\ku$.

The comultiplication of hyperwords in $T(V)$ has a nice expression, as we can see in the following result.

\begin{lemma}{\cite[Thm.13]{R-lyndon words}}\label{thm:coproducto TV en hiperpalabras} Let $u_1,\dots,u_r,v \in L$, with $v<u_r \leq \dots  \leq u_1$.
Then,
\begin{eqnarray*}
\Delta \left( [u_1]_c \cdots [u_r]_c [v]_c^m \right) &=& 1 \ot [u_1]_c \cdots [u_r]_c [v]_c^m + \sum ^{m}_{i=0} \binom{ m }{ i } _{q_{v,v}}
\left[u_1\right]_c\dots\left[u_r\right]_c \left[ v \right]_c ^i \ot \left[ v \right]_c^{m-i}
\\ && + \sum_{ \substack{ l_1\geq \dots  \geq l_p >v, \  l_i \in L \\ 0\leq j \leq m } } x_{l_1,\dots ,l_p}^{(j)} \ot \left[l_1\right]_c \cdots
\left[l_p\right]_c\left[v\right]_c^j;
\end{eqnarray*}
Where $x_{l_1,\dots ,l_p}^{(j)}$ is $\zt$-homogeneous, and
$\deg(x_{l_1,\dots ,l_p}^{(j)})+\deg(l_1\dots  l_p v^j)=
\deg(u_1\dots u_rv^m)$. \qed
\end{lemma}

Another useful result from \cite{R-lyndon words} is the following one.

\begin{lemma}\label{Lemma:subalgebrasWl}
For each $l\in L$ let $W_l$ be the subspace of $T(V)$ generated by
\begin{equation}\label{elementsWl}
    [l_1]_c [l_2]_c\cdots [l_k]_c, \quad k \in \N_0, \, l_i \in L, \, l_1 \geq \ldots \geq l_k \geq l.
\end{equation}
Then $W_l$ is a left coideal subalgebra of $T(V)$.\qed
\end{lemma}
\medskip

We consider another order in $\xx$ as in \cite{U}; it was implicitly used in \cite{Kh}. Let $u,v \in \xx$. We say that $u \succ v$ if $\ell(u)<\ell(v)$, or
$\ell(u)=\ell(v)$ and $u>v$ for the lexicographical order. This order $\succ$ is total, and it is called the \emph{deg-lex order}.

The empty word $1$ is the maximal element for $\succ$, and this order is invariant by left and right multiplication.
\medskip

Let $I$ be an ideal of $T(V)$, and $R=T(V)/I$. Let $\pi: T(V) \rightarrow R$ be the canonical projection. We set:
$$G_I:= \left\{ u \in \xx: u \notin \\ \ku \xx_{\succ u}+I  \right\}.$$
Note that if $u \in G_I$ and $u=vw$, then $v,w \in G_I$. Therefore each $u \in G_I$ is a non-increasing product of Lyndon words of $G_I$.

\begin{proposition}{\cite{Kh,R-lyndon words}}\label{thm:firstPBWbasis} The set $\pi(G_I)$ is a basis of $R$. \qed
\end{proposition}

\noindent In what follows $I$ will denote a Hopf ideal. Consider the set $S_I:=G_I\cap L$. Define $h_I:S_I\to\left\{2,3,\dots \right\}\cup \{\infty\}$
according to the following condition:
\begin{equation}\label{defaltura}
    h_I(u):= \min \left\{ t \in \N : u^t  \in \ku \xx_{\succ u^t} + I \right\}.
\end{equation}
We recall the following result and its corollaries following \cite{Kh}.

\begin{theorem}\label{thm:base PBW Kharchenko} $B_I':= B\left( \left\{1+I\right\}, \left[ S_I \right]_c+I, <, h_I \right)$ is a PBW basis of $H=T(V)/I$.
\qed
\end{theorem}

\begin{corollary}\label{cor:primero}
A word $u$ belongs to $G_I$ if and only if the corresponding hyperletter $\left[u\right]_c$ is not a linear combination, modulo $I$, of greater hyperwords
$\left[ w \right]_c$, $w \succ u$, whose hyperletters are in $S_I$. \qed
\end{corollary}

\begin{corollary}\label{cor:segundo}
If $v\in S_I$ is such that $h_I(v)<\infty$, then $q_{v,v}$ is a root of unity. Moreover, if $\ord q_{v,v}=h$, then $h_I(v)=h$, and $\left[v\right]^h$ is a
linear combination of hyperwords $\left[w\right]_c$, $w\succ v^h$.
\qed
\end{corollary}
\medskip

\subsection{Weyl groupoids and root systems}\label{subsection:weylgrupoide}

We follow the notation in \cite{CH1}. Fix a non-empty set $\Xf$, and
a finite set $I$. For each $i \in I$ we fix a bijective function
$r_i:\Xf\rightarrow \Xf$, and for each $X \in \Xf$ a generalized
Cartan matrix $A^X= (a^X_{ij})_{i,j \in I}$. Let $(\alpha_i)_{i\in
I}$ be the canonical basis of $\Z^I$.

\begin{definition}{\cite{HY,CH1}}\label{def:grupoide de Weyl}
The 4-uple $\cC:= \cC(I, \Xf, (r_i)_{i \in I}, (A^X)_{X \in \Xf})$ is a \emph{Cartan scheme} if it holds:
\begin{itemize}
  \item for any $i \in I$, $r_i^2=id$, and
  \item for any $X \in \Xf$ and any pair $i,j \in I$: $a^X_{ij}=a^{r_i(X)}_{ij}$.
\end{itemize}
For each $i \in I$ and each $X \in \Xf$ we denote by $s_i^X$ the automorphism of $\Z^I$ given by
$$ s_i^X(\alpha_j)=\alpha_j-a_{ij}^X\alpha_i, \qquad j \in I. $$
The \emph{Weyl groupoid} of $\cC$ is the groupoid $\cW(\cC)$ whose set of objects is $\Xf$ and whose morphisms are generated by $s_i^X$, considered as
elements $s_i^X \in \Hom(X, r_i(X))$, $i \in I$, $X \in \Xf$.
\end{definition}

In general we denote $\cW(\cC)$ simply by $\cW$, and for each $X \in \Xf$:
\begin{align}\label{eq:def Hom grupoide}
    \Hom(\cW,X) & := \cup_{Y \in \Xf} \Hom(Y,X),
    \\ \label{eq:raiz real} \Delta^{X \ re} &:= \{ w(\alpha_i): \ i \in I, \ w \in \Hom(\cW,X) \}.
\end{align}
$\Delta^{X \ re}$ is the set of \emph{real roots} of $X$. Each
$w\in\Hom(\cW,X_1)$ is written as a product
$s_{i_1}^{X_2}s_{i_2}^{X_3}\cdots s_{i_m}^{X_{m+1}}$, where
$X_j=r_{i_{j-1}} \cdots r_{i_1}(X_1)$, $i \geq 2$. We denote it by
$w= \id_{X_1} s_{i_1} \cdots s_{i_m}$: it means that $w \in
\Hom(\cW,X_1)$, because each $X_j \in \Xf$ is univocally determined
by this condition. The \emph{length} of $w$ is defined by
$$ \ell(w)= \min \{ n \in \N_0: \ \exists i_1, \ldots, i_n \in I \mbox{ such that }w=\id_X s_{i_1} \cdots s_{i_n} \}. $$
We assume that $\cW$ is \emph{connected}: that is, $\Hom(Y,X) \neq
\emptyset$, for any pair $X,Y \in \Xf$.

\begin{definition}{\cite{HY,CH1}}\label{def:sistema de raices}
Given a Cartan scheme $\cC$, consider for each $X \in \Xf$ a set $\Delta^X \subset \Z^I$. We say that $\cR:= \cR(\cC, (\Delta^X)_{X \in \Xf} )$ is a
\emph{root system of type} $\cC$ if
\begin{enumerate}
  \item for any $X \in \Xf$, $\Delta^X= (\Delta^X \cap \N_0^I) \cup -(\Delta^X \cap \N_0^I)$,
  \item for any $i \in I$ and any $X \in \Xf$, $\Delta^X \cap \Z \alpha_i= \{\pm \alpha_i \}$,
  \item for any $i \in I$ and any $X \in \Xf$, $s_i^X(\Delta^X)=\Delta^{r_i(X)}$,
  \item if $m_{ij}^X:= |\Delta^X \cap (\N_0\alpha_i+\N_0 \alpha_j)|$, then $(r_ir_j)^{m_{ij}^X}(X)=X$ for any pair $i \neq j \in I$ and any $X\in \Xf$.
\end{enumerate}
$\Delta^X_+:= \Delta^X \cap \N_0^I$ is called the set of
\emph{positive roots}, and $\Delta^X_-:=-\Delta^X_+$ is the set of
\emph{negative roots}.
\end{definition}
\begin{remark}\label{obs:raices reales incluidas}
From (2) and (3) we deduce that $\Delta^{X \ re} \subset \Delta^X$, for any $X\in\Xf$.
\end{remark}

For each positive root $\alpha=\sum_i n_i \alpha_i$, the \emph{support} of $\alpha$ is the set
$$\sop \alpha:= \{i: \ 1\leq i\leq\theta, n_i\neq 0 \}. $$

By (3) we have that $w(\Delta^X)= \Delta^Y$ for any $w \in \Hom(Y,X)$. We say that $\cR$ is \emph{finite} if $\Delta^X$ is finite for some $X\in \Xf$.
By \cite[Lemma 2.11]{CH1}, it is equivalent to the fact that all the sets $\Delta^X$ are finite, $X \in \Xf$, and also that $\cW$ is finite. Moreover,
for any pair  $i\neq j\in I$ and any $X \in \Xf$, we have that $k\alpha_i + \alpha_j \in \Delta^X$ if and only if $0 \leq k \leq -a_{ij}^X$. Therefore,
\begin{equation}\label{eq: caracterizacion aij}
 a_{ij}^X = -\max \{ k \in \N_0: \ k\alpha_i + \alpha_j \in \Delta^X \}.
\end{equation}

A fundamental result involving root systems is the following one:

\begin{theorem}{\cite[Thm. 2.10]{CH2}}  \label{theorem:decomposition-of-roots}
For every $\alpha \in \Delta^X_+ \setminus \{ \alpha_i: i=1, \ldots \theta\}$, there exist $\beta, \gamma \in \Delta^X_+$ such that $\alpha=\beta+ \gamma$.
\qed
\end{theorem}

We give now some results about real roots and the length of elements.

\begin{lemma}{\cite[Cor. 3]{HY}}  \label{Lemma:longitudHY}
Let $m \in \N$, $X, Y \in \Xf$, $i_1, \ldots, i_m,j \in I$, $w=\id_X s_{i_1} \cdots s_{i_m} \in \Hom(Y,X)$, where $\ell (w)=m$. Then,
\begin{itemize}
  \item $\ell (w s_j)=m+1$ if and only if $w(\alpha_j) \in \Delta^X_+$,
  \item $\ell (w s_j)=m-1$ if and only if $w(\alpha_j) \in \Delta^X_-$.\qed
\end{itemize}
\end{lemma}

\begin{proposition}{\cite[Prop. 2.12]{CH1}}  \label{Prop:maxlongitudCH}
For any $w=\id_X s_{i_1} \cdots s_{i_m}\in\cW$ such that
$\ell(w)=m$, the roots $\beta_j=s_{i_1}\cdots
s_{i_{j-1}}(\alpha_{i_j})\in\Delta^X$ are positive and all
different. If $\cR$ is finite and $w$ is an element of maximal
length, say $N$, then $\Delta^X_+= \{\beta_j | 1 \leq j\leq N\}$.
The roots $\beta_1,$,...,$\beta_N$ are pairwise different, and hence
for each $\alpha \in \Delta^X_+$ there exist $i_1,\ldots,i_k,j\in I$
such that $\alpha=s_{i_k}\cdots s_{i_1}(\alpha_j)$. \qed
\end{proposition}

Call $\Delta^V_+ $ the set of degrees of a PBW basis of $\cB(V)$, counted with their multiplicities, as in \cite{H-Weyl grp}. It does not depend on the PBW
basis, see \cite{H-Weyl grp,AA}. We can attach a Cartan scheme $\cC$, a Weyl groupoid $\cW$ and a root system $\cR$, see \cite[Thms. 6.2, 6.9]{HS}. To do
this, define for each $1 \leq i \neq j \leq \theta$,
\begin{equation}
-a_{ij}:= \min \left\{ n \in \mathbb{N}_0: (n+1)_{q_{ii}}
(1-q_{ii}^n q_{ij}q_{ji} )=0 \right\}, \label{defn:mij}
\end{equation}
and set $a_{ii}=2$, $s_i\in\Aut(\Z^\theta)$ such that $s_i(\alpha_j)=\alpha_j-a_{ij}\alpha_i$.

Set $\overline{q}_{rs}=\chi(s_i(\alpha_r),s_i(\alpha_s))$. Let $V_i$
be another vector space of dimension $\theta$, and attach to it the
matrix $\overline{\qmb}=(\overline{q}_{rs})$. By \cite{H-Weyl grp},
$\Delta^{V_i}_+ = s_i\left( \Delta^V_+ \setminus\{\alpha_i\}\right)
\cup \{\alpha_i\}$. If we consider $\Delta^V=\Delta^V_+\cup
(-\Delta^V_+) $, last equation lets us to define the Weyl groupoid
of $V$, whose root system is defined by the sets $\Delta^{V'}$, $V'$
obtained after to apply some reflections to the matrix of $V$.

\subsection[Defining relations of Nichols algebras]{Defining relations of Nichols algebras of diagonal type}
\label{section:relaciones} \

\begin{proposition}{\cite[Prop. 3.1]{A-presentation}}\label{Prop:basePBWortogonal}
Assume that the braiding matrix is symmetric. Then a PBW basis of Lyndon hyperwords of $\cB(V)$ is orthogonal with respect to the bilinear form in
Proposition \ref{prop:formabilineal}.\qed
\end{proposition}

\begin{corollary}\cite[Cor. 3.2]{A-presentation}\label{Coro:normageneradores}
If $u= x_{\beta_M}^{n_M} \cdots x_{\beta_1}^{n_1}$, where $0 \leq n_j < N_{\beta_j}$, then
\begin{equation}\label{normaPBWgenerador}
    c_u:= (u|u)= \prod_{j=1}^M n_j!_{q_{\beta_j}} c_{x_{\beta_j}}^{n_j} \neq 0.
\end{equation}\qed
\end{corollary}

\begin{remark}\label{Rem:calculoformabilineal} Notice that:
$$ (x_{\beta_i}x_{\beta_j} | u)= (x_{\beta_i} | u_{(1)})(x_{\beta_j} | u_{(2)}) =  d_{i,j} c_{x_{\beta_i}}c_{x_{\beta_j}}, $$
where $d_{i,j}$ is the coefficient $x_{\beta_i} \ot x_{\beta_j}$ for the expression of $\Delta(u)$ in terms of the PBW basis (both factors of the tensor
product).
\end{remark}

For each pair $1\leq i \leq j \leq \theta$, we denote
$$B_{ij}:= \left\{ x_{\beta_j}^{n_j} \cdots x_{\beta_i}^{n_i}: \ 0\leq n_k <N_{\beta_k} \right\};$$
that is, the set of hyperwords whose hyperletters are between $x_{\beta_i}$ and $x_{\beta_j}$

Let $(W,d)$ be a braided vector space of diagonal type that admits a basis $\widehat x_1, \ldots, \widehat x_\theta$ such that, for some $\widehat q_{ij}
\in \ku^\times$, $d(\widehat x_i \otimes \widehat x_j) = \widehat q_{ij} \widehat x_j \otimes \widehat x_i$, where $\widehat q_{ij}=\widehat q_{ji}$ ,
and $(V,c)$, $(W,d)$ are twist equivalent:
$$ q_{ij}q_{ji}= \widehat q_{ij} \widehat q_{ji}, \quad q_{ii}= \widehat q_{ii}, \qquad 1 \leq i \neq j \leq \theta.$$
Call $\widehat x_\beta= [l_\beta]_d$, the hyperletter corresponding
to $l_\beta$ for the braiding $d$. If $u= x_{\beta_M}^{n_M} \cdots
x_{\beta_1}^{n_1}$, call also
$$\widehat u= \widehat x_{\beta_M}^{n_M} \cdots \widehat x_{\beta_1}^{n_1}.$$

Let $\sigma: \zt \times \zt \rightarrow \ku^\times$ be the bicharacter determined by the condition
\begin{equation}\label{cociclo}
 \sigma(\alpha_i,\alpha_j) = \left\{ \begin{array}{lr} \widehat q_{ij}q_{ij}^{-1}, & i \leq j \\ 1, & i>j \end{array} \right.
\end{equation}
Define $t_{\alpha_i}=1$ for any $1 \leq i \leq \theta$, and inductively,
$$ t_\beta= \sigma(\beta_1, \beta_2) t_{\beta_1} t_{\beta_2}, \qquad \Sh(l_\beta)= (l_{\beta_1}, l_{\beta_2}) . $$
For each $u= x_{\beta_M}^{n_M} \cdots x_{\beta_1}^{n_1}$ call
\begin{equation}\label{formulacoeficientes}
 f(u):= \prod_{1\leq i<j\leq M}\sigma(\beta_j,\beta_i)^{n_in_j} \ \prod_{1\leq i\leq M}\sigma(\beta_i,\beta_i)^{\binom{n_i}{2}}t_{\beta_i}^{n_i}.
\end{equation}

Finally, for each pair $1\leq i<j\leq\theta$ and $u= x_{\beta_M}^{n_M} \cdots x_{\beta_1}^{n_1}$, let
\begin{equation}\label{escalares}
c_{i,j}^u := \frac{f(u) \ (\widehat x_{\beta_i} \widehat x_{\beta_j} | \widehat u)}{\sigma(\beta_i, \beta_j) t_{\beta_i} t_{\beta_j} c_{\widehat u}},
\end{equation}
where $(\cdot|\cdot)$ is the bilinear form corresponding to $(W,d)$, and $c_{\widehat u}$ is the scalar in Corollary \ref{Coro:normageneradores}. Such
scalars let us to give a presentation by generators and relations as follows.

\begin{theorem}{\cite[Thm. 3.9]{A-presentation}}\label{Thm:presentacion}
Let $(V,c)$ be a finite-dimensional braided vector space of diagonal type such that $\Delta^V_+$ is finite. Let $x_1, \cdots, x_\theta$ be a basis of $V$
such that $c(x_i \ot x_j)= q_{ij} x_j \ot x_i$, where $(q_{ij}) \in (\ku^\times)^{\theta \times \theta}$ is the braiding matrix, and let
$\{x_{\beta_k} \}_{\beta_k \in \Delta^V_+}$ be the set of hyperletters corresponding to the fixed order of the basis of $V$.

Then $\cB(V)$ is presented by generators $x_1, \ldots, x_\theta$, and relations
\begin{align}
x_{\beta}^{N_{\beta}}& =0, \qquad \beta \in \Delta^V_+, \ \ord(q_\beta)=N_\beta < \infty, \label{powerrootvector1}
\\ \left[ x_{\beta_i}, x_{\beta_j} \right]_c &= \sum_{u \in B_{ij}- \{ x_{\beta_j}x_{\beta_i} \}: \ \ \deg u= \beta_i+\beta_j} c_{i,j}^u \ u,
\label{quantumSerregeneralizadas1}
\\ & 1 \leq i<j \leq M, \ \Sh(l_{\beta_i} l_{\beta_j}) =( l_{\beta_i}, l_{\beta_j} ), \ l_{\beta_i} l_{\beta_j} \neq l_{\beta_k}, \forall k, \nonumber
\end{align}
where $c_{i,j}^u$ are as in \eqref{escalares}. Moreover, $\{ x_{\beta_M}^{n_M} \cdots x_{\beta_1}^{n_1}: \ 0 \leq n_j < N_{\beta_j} \}$ is a basis of
$\cB(V)$.\qed
\end{theorem}

\medskip

\section[Lusztig Isomorphisms]{Lusztig Isomorphisms of Nichols algebras of diagonal type}\label{section:Lusztig isomorphisms}

In this Section we recall the Lusztig isomorphisms \cite{H-isom} of Nichols algebras of diagonal type, which are a generalization of the isomorphisms of
quantized enveloping algebras in \cite{L}. We shall consider different quotients of the tensor algebra of a braided vector space of diagonal type and
the Drinfeld doubles of their bosonizations by a free abelian group.
\medskip

\textbf{Notation:} Let $\chi : \Z^{\theta} \times \Z^{\theta} \to \ku^{\times}$ be a bicharacter, $q_{ij}=\chi(\alpha_i, \alpha_j)$. Then $\chi^{\op}$
and $\chi^{-1}$ will denote the bicharacters:
$$ \chi^{op}(\alpha,\beta):= \chi(\beta, \alpha), \quad \chi^{-1}(\alpha,\beta):= \chi(\alpha, \beta)^{-1}, \qquad \alpha, \beta \in \Z^\theta. $$
Also, for any automorphism $s:\Z^\theta \to \Z^\theta$, $s^*\chi$ will denote the bicharacter defined by
\begin{equation}\label{eqn:accion w sobre chi}
(s^*\chi)(\alpha,\beta):= \chi\left(s^{-1}(\alpha), s^{-1}(\beta) \right), \qquad \alpha, \beta \in \Z^\theta.
\end{equation}
\bigbreak

Let $(V,c)$ a braided vector space of diagonal type, whose braiding matrix is $(q_{ij})$. We consider $T(V)$ as an algebra in the category of
Yetter-Drinfeld modules over $\ku \Z^\theta$ as above. We follow the results in \cite[Section 4.1]{H-isom}.

\begin{definition}
The Drinfeld double $\cU(\chi)$ of the Hopf algebra
$T(V)\#\ku\Z^{\theta}$ is the algebra generated by elements $E_i$,
$F_i$, $K_i$, $K_i^{-1}$, $L_i$, $L_i^{-1}$, $1 \leq i \leq \theta$,
and relations
\begin{align*}
&XY=YX, \qquad X,Y \in  \{ K_i^{\pm}, L_i^{\pm}: 1 \leq i \leq \theta \},
\\ &K_iK_i^{-1}=L_iL_i^{-1}=1,
\\ &K_iE_jK_i^{-1}=q_{ij}E_j, \quad & L_iE_jL_i^{-1}=q_{ji}^{-1}E_j,
\\ &K_iF_jK_i^{-1}=q_{ij}^{-1}F_j, \quad & L_iF_jL_i^{-1}=q_{ji}F_j,
\\ &E_iF_j-F_jE_i=\delta_{i,j}(K_i-L_i).
\end{align*}
It admits a Hopf algebra structure, where the comultiplication satisfies
\begin{align*}
&\Delta(K_i)=K_i \ot K_i, &\Delta(E_i)=E_i \ot 1 + K_i \ot E_i,
\\ &\Delta(L_i)=L_i \ot L_i, &\Delta(F_i)=F_i \ot L_i + 1 \ot F_i,
\end{align*}
and then $\eps(K_i)=\eps(L_i)=1$, $\eps(E_i)=\eps(F_i)=0$.
\end{definition}
Notice that $\cU(\chi)$ is a $\Z^{\theta}$-graded Hopf algebra, where the graduation is characterized by the following conditions:
$$ \deg(K_i)=\deg(L_i)=0, \qquad \deg(E_i)=\alpha_i, \qquad \deg(F_i)=-\alpha_i. $$

$\cU^+(\chi)$ (respectively, $\cU^-(\chi)$) denotes the subalgebra
generated by $E_i$ (respectively, $F_i$), $1 \leq i \leq \theta$,
$\cU^{+0}(\chi)$ (respectively, $\cU^{-0}(\chi)$) is the subalgebra
generated by $K_i$, $K_i^{-1}$ (respectively, $L_i$, $L_i^{-1}$), $1
\leq i \leq \theta$, and finally $\cU^0(\chi)$ is the subalgebra
generated by $K_i$, $K_i^{-1}$, $L_i$ and $L_i^{-1}$. Note that
$\cU^0(\chi)$ is isomorphic to $\ku \Z^{2\theta}$ as Hopf algebras.
Moreover, the subalgebra generated by $\cU^+(\chi)$, $K_i$ and
$K_i^{-1}$, $1 \leq i \leq \theta$, is isomorphic to
$T(V)\#\ku\Z^{\theta}$, so $\cU^+(\chi)$ is isomorphic to $T(V)$ as
braided graded Hopf algebras in the category of Yetter-Drinfeld
modules over $\ku \Z^\theta$, where we consider the actions and
coactions:
$$ K_i \cdot E_j=q_{ij} E_j, \qquad \delta(E_i)= K_i \ot E_i. $$
\medskip

We will consider a family of useful isomorphisms as in \cite[Section 4.1]{H-isom}.

\begin{proposition}\label{prop:morfismos phi}
\begin{enumerate}
  \item[(a)] For any $\underline{a}=(a_1,\ldots,a_\theta)\in(\ku^\times)^\theta$ there exists a unique algebra automorphism $\varphi_{\underline a}$ of
$\cU(\chi)$ such that
\begin{equation}\label{eqn:varphi a}
\varphi_{\underline a}(K_i)=K_i, \quad \varphi_{\underline a}(L_i)=L_i, \quad \varphi_{\underline a}(E_i)=a_iE_i, \quad
\varphi_{\underline a}(F_i)=a_i^{-1}F_i.
\end{equation}
  \item[(b)] There exists a unique algebra automorphism $\phi_1$ of $\cU(\chi)$ such that
\begin{equation}\label{eqn:phi1}
 \phi_1(K_i)=K_i^{-1}, \quad \phi_{1}(L_i)=L_i^{-1}, \quad \phi_1(E_i)=F_iL_i^{-1}, \quad \phi_1(F_i)=K_i^{-1}E_i.
\end{equation}
  \item[(c)] There exists a unique algebra isomorphism $\phi_2: \cU(\chi) \to \cU(\chi^{-1})$ such that
\begin{equation}\label{eqn:phi2}
 \phi_2(K_i)=K_i, \quad \phi_2(L_i)=L_i, \quad \phi_2(E_i)=F_i, \quad \phi_2(F_i)=-E_i.
\end{equation}
  \item[(d)] There exists a unique Hopf algebra isomorphism $\phi_3:\cU(\chi) \to \cU(\chi^{\op})^{\cop}$ such that
\begin{equation}\label{eqn:phi3}
 \phi_3(K_i)=L_i, \quad \phi_3(L_i)=K_i, \quad \phi_3(E_i)=F_i, \quad \phi_3(F_i)=E_i.
\end{equation}
  \item[(e)] There exists a unique algebra antiautomorphism $\phi_4$ of $\cU(\chi)$ such that
\begin{equation}\label{eqn:phi4}
 \phi_4(K_i)=K_i, \quad \phi_4(L_i)=L_i, \quad \phi_4(E_i)=F_i, \quad \phi_(F_i)=E_i.
\end{equation}
  \item[(f)] Let $\underline a=(-1, \cdots, -1)$. The antipode $\SS$ of $\cU(\chi)$ is given by the composition $\SS=\phi_1\phi_4\varphi_{\underline a}$.
Also, $\phi_4^2=\id$. \qed
\end{enumerate}
\end{proposition}

As in \cite{H-isom} we will consider some skew-derivations.
$\underline \Delta$ will denote the braided comultiplication of
$\cU^+(\chi)$, which is $\N_0$-graded: if $E\in\cU^+(\chi)$ is
homogeneous of degree $n$, and $k\in\{0,1,\ldots,n\}$,
$\underline\Delta_{n-k,k}(E)$ will denote the component of
$\underline \Delta(E)$ in $\cU^+(\chi)_{n-k}\ot\cU^+(\chi)_k$.

\begin{proposition}\label{prop:derivadas torcidas}
For any $i\in\{1,\ldots,\theta\}$ there exist linear endomorphisms $\partial^K_i$, $\partial^L_i$ of $\cU^+(\chi)$ such that
$$ EF_i-F_iE= \partial^K_i(E)K_i- L_i \partial^L_i(E) \quad \mbox{ for all }E \in \cU^+(\chi). $$
Such endomorphisms are given by:
$$ \underline \Delta_{n-1,1}(E)= \sum_{i=1}^\theta \partial_i^K(E) \ot E_i,  \qquad \underline \Delta_{1,n-1}(E)= \sum_{i=1}^\theta E_i \ot
\partial_i^L(E), \qquad E \in \cU^+(\chi)_n,$$
and satisfy the following conditions:
\begin{align*}
\partial_i^K(1)&= \partial_i^L(1)=0,
\\ \partial_i^K(E_j)&= \partial_i^L(E_j)=\delta_{i,j},
\\ \partial_i^K(EE')&= \partial_i^K(E)(K_i \cdot E')+ E \partial_i^K(E'),
\\ \partial_i^L(EE')&= \partial_i^L(E)E'+ (L_i^{-1}\cdot E) \partial_i^L(E'),
\end{align*}
for all $j\in\{1,\ldots,\theta\}$, and all pair of elements $E,E' \in \cU^+(\chi)$.\qed
\end{proposition}

We recall now a characterization of quotients of the algebra $\cU(\chi)$ with a triangular decomposition \cite[Section 4.1]{H-isom}. According to
\cite[Prop. 4.14]{H-isom}, the multiplication
\begin{equation}\label{eqn:multiplicacion iso}
m: \cU^+(\chi) \ot \cU^0(\chi) \ot \cU^-(\chi) \to \cU(\chi)
\end{equation}
is an isomorphism of $\Z^\theta$-graded vector spaces.

\begin{proposition}\label{prop:desc triangular cocientes}
Let $\cI^+\subset\cU^+\cap\ker\eps$ (respectively, $\cI^-\subset\cU^-\cap\ker\eps$) be an ideal of $\cU^+(\chi)$ (respectively, $\cU^-(\chi)$). The
following conditions are equivalent:
\begin{itemize}
  \item The multiplication \eqref{eqn:multiplicacion iso} induces an isomorphism
  $$ m: \cU^+(\chi)/\cI^+ \ot \cU^0(\chi) \ot \cU^-(\chi)/\cI^- \to \cU(\chi)/ (\cI^+ +\cI^-). $$
  \item The vector spaces $\cI^+ \cU^0(\chi) \cU^-(\chi)$ and $\cU^+(\chi)\cU^0(\chi)\cI^-$ are ideals of $\cU(\chi)$.
  \item For all $X\in\cU^0(\chi)$ and all $i\in\{1,\ldots,\theta\}$ we have
  \begin{equation*}
  \begin{array}{ccc} X \cdot \cI^+ \subset \cI^+,  & \partial_i^K(\cI^+) \subset \cI^+, & \partial_i^K\left( \phi_4(\cI^-) \right) \subset \phi_4(\cI^-),
  \\ X \cdot \cI^- \subset \cI^-,  & \partial_i^L(\cI^+) \subset \cI^+, & \partial_i^L\left( \phi_4(\cI^-) \right) \subset \phi_4(\cI^-).
  \end{array}
  \end{equation*}
\end{itemize}\qed
\end{proposition}

\begin{lemma}{\cite[Cor. 4.20]{H-isom}}\label{Lemma:caracterizacion ideales triangulacion}
Let $I^+$ be an braided biideal of $\cU^+(\chi)$, which is also a Yetter-Drinfeld $\cU^0(\chi)$-submodule and satisfies
$I^+ \subset \oplus_{n \geq 2} \cU^+(\chi)_n$. Then $I^+\cU^0(\chi)\cU^-(\chi)$ is a Hopf ideal of $\cU(\chi)$.\qed
\end{lemma}

We will assume that all the integers $m_{ij}:=-a_{ij}$ of
\eqref{defn:mij} associated with the matrices $(q_{ij})$ are
defined. Then we consider the automorphisms
$s_{p,\chi}:\Z^\theta\to\Z^\theta$. We define also the scalars
\begin{equation}\label{eqn:escalares lambda}
\lambda_i(\chi):= (-a_{pi})_{q_{pp}}\prod_{s=0}^{-a_{pi}-1}(q_{pp}^sq_{pi}q_{ip}-1) ,  \qquad i \neq p.
\end{equation}
Denote by $\Eb_i$, $\Fb_i$, $\Kb_i$, $\Lb_i$ the generators corresponding to $\cU(\schi)$, and by $\qb_{ij}=\schi(\alpha_i,\alpha_j)$ the coefficients
of the braiding matrix corresponding to $\schi$.

\begin{definition}
We say that $p\in\{1,\ldots,\theta\}$ is a \emph{Cartan vertex} if it satisfies
$$ q_{pp}^{-a_{pj}}q_{pj}q_{jp}=1,\quad \mbox{ for every }j \neq p.$$
In such case, note that the existence of the integers $-a_{pj}$
implies that $\ord q_{pp}\geq -a_{pj}+1$.
\end{definition}

We denote by $\cO(\chi)$ the union of the orbits of the simple roots $\alpha_p$ by the action of the Weyl groupoid, where $p$ is a Cartan vertex.
\bigbreak

Fix $p\in\{1,\ldots,\theta\}$. For any $i\neq p$ we define as in \cite{H-isom},
$$ E_{i,0(p)}^+, E_{i,0(p)}^-:=E_i, \qquad F_{i,0(p)}^+, F_{i,0(p)}^-:=F_i, $$
and recursively,
\begin{align}
E_{i,m+1(p)}^+ &:= E_pE_{i,m(p)}^+ - (K_p \cdot E_{i,m(p)}^+)E_p = (\ad_c E_p)^{m+1}E_i,
\\ E_{i,m+1(p)}^- &:= E_pE_{i,m(p)}^- - (L_p \cdot E_{i,m(p)}^-)E_p,
\\ F_{i,m+1(p)}^+ &:= F_pF_{i,m(p)}^+ - (L_p \cdot F_{i,m(p)}^+)F_p,
\\ F_{i,m+1(p)}^- &:= F_pF_{i,m(p)}^- - (K_p \cdot F_{i,m(p)}^-)F_p.
\end{align}
When $p$ is explicit, we simply denote $E_{i,m(p)}^\pm$ by $E_{i,m}^\pm$. By \cite[Cor. 5.4]{H-isom} the following identity holds for any $m \in \N_0$:
\begin{equation}\label{eqn:corchete E+ con Fp}
E_{i,m}^+ F_i - F_i E_{i,m}^+= (m)_{q_{pp}}(q_{pp}^{m-1}q_{pi}q_{ip}-1)L_p E_{i,m-1}^+.
\end{equation}

Fix a braided graded Hopf algebra $\cB\cong T(V)/I$, where $I$ is a graded Hopf ideal generated by homogeneous elements of degree $\geq 2$. For each
$1\leq j\leq\theta$, $p\neq j$, we define
\begin{equation}\label{eqn:conjuntos Mij}
     M_{p,j}^{\pm}(\cB):=\left\{ E_{j,m}^{\pm}: m \in \N_0 \right\}.
\end{equation}
In what follows we consider $\ord(1)=1$.

\begin{remark}\label{obs:nilpotency}
If $E_i^N=0$ in $\cB$, with $N$ minimal (it is called the nilpotency
order of $x_i$), then $q_{ii}$ is a primitive root of unity of order
$N$. Moreover, $(\ad_c E_i)^NE_j=0$.

Also, the nilpotency order of $E_i$ is infinite if $E_i^n\neq 0$ for
all $n\in\N$.
\end{remark}

We recall a result from \cite{A-standard} extending \cite[Prop. 1, Eqn. (18)]{H-Weyl grp}.

\begin{lemma}\label{lema:kernel derivaciones generado por ad}
For each $p\in\{1,\ldots,\theta\}$, let $\cB_{\pm p}$ be the
subalgebra of $\cB$ generated by $\cup_{j\neq p}M_{p,j}^{\pm}(\cB)$,
and denote $N_p=\ord q_{pp}$. There exist isomorphisms of graded
vector spaces:
\begin{itemize}
    \item $\ker(\partial_p^K)\cong\cB_{+p} \ot\ku\left[E_p^{N_p}\right]$,
    $\ker(\partial_p^L)\cong\cB_{-p} \ot\ku\left[E_p^{N_p}\right]$, if
$1<\ord(q_{pp})<\infty$ but $E_p$ is not nilpotent, or
    \item $\ker(\partial_p^K)\cong\cB_{+p}$, $\ker(\partial_p^L)\cong\cB_{-p}$, if $\ord(q_{pp})$ is the nilpotency order of
    $E_p$ or $q_{pp}=1$.
\end{itemize}
Moreover, the multiplication induces an isomorphism of graded vector
spaces $\cB \cong \cB_p^{\pm} \ot \ku\left[E_p \right]$.\qed
\end{lemma}
\medskip

Set $N_p= \ord q_{pp}$. We call, following \cite{H-isom}, $\cI_p^+(\chi)$ (respectively, $\cI_p^-(\chi)$) to the ideal of $\cU^+(\chi)$ (respectively,
$\cU^-(\chi)$) generated by
\begin{itemize}
  \item[(a)] $E_p^{N_p}$ (respectively, $F_p^{N_p}$), if $p$ is not a Cartan vertex,
  \item[(b)] $E_{i,m_{pi}+1}^+$ (respectively, $F_{i,m_{pi}+1}^+$), for each $i$ such that $N_p>m_{pi}+1$.
\end{itemize}
Notice that $E_{i,m_{pi}+1}^+\in\cI_p^+(\chi)$ for any $i$ such that $N_p=m_{pi}+1$. We denote:
$$\cU_p(\chi):=\cU(\chi)/\left(\cI_p^+(\chi)+\cI_p^-(\chi) \right), \quad \cU_p^+(\chi):=\cU^+(\chi)/\cI_p^+(\chi), \quad
\cU_p^-(\chi):= \cU^-(\chi)/\cI_p^-(\chi).$$

$I^+(\chi)$ will denote the ideal of $\cU^+(\chi)$ such that the quotient $\cU^+(\chi)/I^+(\chi)$ is isomorphic to the Nichols algebra of $V$; that is,
the greatest braided Hopf ideal of $\cU^+(\chi)$ generated by elements of degree $\geq2$. Call $I^-(\chi)= \phi_4(I^+(\chi))$, where $\phi_4$
is defined by \eqref{eqn:phi4}, and
$$ \u^{\pm}(\chi):= \cU^{\pm}(\chi) / I^{\pm}(\chi), \qquad \u(\chi):= \cU(\chi) / (I^-(\chi)+I^+(\chi)).$$
In such case, $\u(\chi)$ is the Drinfeld double of the algebra $\u^+(\chi)\#\ku\Z^\theta$, where $\ku\Z^\theta=\cU^{+0}(\chi)$.

The Lusztig isomorphisms can be defined in this general context.

\begin{theorem}{\cite[Lemma 6.5, Theorem 6.12]{H-isom}}\label{thm: iso Lusztig-Heck}
There exist algebra morphisms
\begin{equation}\label{eqn:iso lusztig para cUp}
T_p, T_p^-: \cU_p(\chi) \to \cU_p(\schi)
\end{equation}
univocally determined by the following conditions:
\begin{align*}
&T_p(K_p)=T_p^-(K_p)=\Kb_p^{-1}, &T_p(K_i)=T_p^-(K_i)=\Kb_p^{m_{pi}}\Kb_i,
\\ &T_p(L_p)=T_p^-(L_p)=\Lb_p^{-1}, &T_p(L_i)=T_p^-(L_i)=\Lb_p^{m_{pi}}\Lb_i,
\\ &T_p(E_p)=\Fb_p\Lb_p^{-1}, &T_p(E_i)=\Eb^+_{i,m_{pi}},
\\ &T_p(F_p)=\Kb_p^{-1}\Eb_p, &T_p(F_i)=\lambda_p(\schi)^{-1}\Fb^+_{i,m_{pi}},
\\ &T_p^-(E_p)=\Kb_p^{-1}\Fb_p, &T_p^-(E_i)=\lambda_p(\schi^{-1})^{-1}\Eb^-_{i,m_{pi}},
\\ &T_p^-(F_p)=\Eb_p\Lb_p^{-1}, &T_p^-(F_i)=\Fb^-_{i,m_{pi}}.
\end{align*}
for every $i\neq p$. Both are isomorphisms satisfying
$$ T_pT_p^-=T_p^-T_p=\id, \qquad T_p(\cU_{+p}^+(\chi))=\cU_{-p}^+(\schi).$$
Moreover, there exists $\mu \in (\ku^\times)^\theta$ such that
\begin{equation}\label{eqn:Tp con phi4}
T_p \circ \phi_4 = \phi_4 \circ T_p^- \circ \varphi_{\mu}.
\end{equation}
Such isomorphisms induce algebra isomorphisms (denoted by the same name):
$$ T_p, T_p^-: \u(\chi) \to \u(s_p^* \chi).$$\qed
\end{theorem}

\begin{remark}\label{rem:corchete via Tp}
If the homogeneous elements $X,Y\in\cU_p^+(\chi)$ are such that $T_p(X),T_p(Y)\in \cU_p^+(\schi)$, as $\deg T_p(X)=s_p(\deg X)$, it follows that
$$ T_p \left( [X,Y]_c \right) = \left[ T_p(X), T_p(y) \right]_{\cb}. $$
\end{remark}

\section[Explicit presentation by generators and relations]{An explicit presentation by generators and relations of Nichols algebras of diagonal type}
\label{section:presentacion explicita}

We shall obtain a family of isomorphisms induced by the ones in the previous Section. In this case we shall consider a quotient of $\cU(\chi)$ by an ideal
which is smaller than $(I^-(\chi)+I^+(\chi))$. Such ideal will be generated by some of the relations in Theorem \ref{Thm:presentacion}, and will be the
smallest one such that it is possible to define all the family of isomorphisms over the Weyl groupoid. It will give us a relation between the Hilbert
series of these algebras, and new sets of roots. We shall use at the end the uniqueness of the root system, when the Weyl groupoid is finite.
\medskip

We introduce some notation. We denote $\widetilde{q_{ij}}=q_{ij}q_{ji}$. Also,
$$ x_{i_1i_2\cdots i_k}=(\ad_c x_{i_1})\cdots(\ad_c x_{i_{k-1}})x_{i_k}, \qquad i_j\in\{1,\ldots,\theta\}. $$
For each $m\in \N$, we define the elements $x_{(m+1)\alpha_i+m\alpha_j}\in\cU(\chi)$ recursively:
\begin{itemize}
  \item if $m=1$, $x_{2\alpha_i+\alpha_j}:= (\ad_c x_i)^2 x_j=x_{iij}$,
  \item $x_{(m+2)\alpha_i+(m+1)\alpha_j}:= [ x_{(m+1)\alpha_i+m\alpha_j}, x_{ij} ]_c$.
\end{itemize}
We give now the main result of this section, which is Theorem \ref{thm:presentacion minima}: it gives an explicit presentation by generators and
relations of any Nichols algebra of diagonal type with finite root system. We begin by proving several Lemmata to show the existence of Lusztig
isomorphisms for some Hopf algebras. These Hopf algebras are intermediate between the tensor algebra and the Nichols algebras of a given braided vector
space. Finally we use those Lusztig isomorphisms to prove the Theorem.

\begin{theorem}\label{thm:presentacion minima}
Let $(V,c)$ be a finite-dimensional braided vector space of diagonal type, with braiding matrix $(q_{ij})_{1\leq i,j\leq\theta}$, $\theta=\dim V$, and fix
a basis $x_1,\ldots,x_\theta$ of $V$ such that $c(x_i\ot x_j)=q_{ij} x_j\ot x_i$. Let $\chi$ be the bicharacter associated to $(q_{ij})$. Assume that the
root system $\Delta^\chi$ is finite. Then $\cB(V)$ is presented by generators $x_1,\ldots,x_\theta$ and relations:
\begin{align}
&x_\alpha^{N_\alpha}, &\alpha \in \cO(\chi);\label{eqn:potencia raices}
\\ &(\ad_cx_i)^{m_{ij}+1}x_j, & q_{ii}^{m_{ij}+1} \neq 1;\label{eqn:relacion quantum Serre}
\\ &x_i^{N_i}, & i \mbox{ is not a Cartan vertex};\label{eqn:potencia raices simples}
\end{align}

\noindent $\diamond$ if $i,j \in \{1, \ldots, \theta \}$ satisfy
$q_{ii}=\widetilde{q_{ij}}=q_{jj}=-1$, and there exists $k\neq i,j$
such that $\widetilde{q_{ik}}^2\neq1$ or
$\widetilde{q_{jk}}^2\neq1$,
\begin{equation}\label{eqn:relacion dos vertices con -1}
x_{ij}^2;
\end{equation}

\noindent $\diamond$ if $i,j,k \in \{1, \ldots, \theta \}$ satisfy
$q_{jj}=-1$,
$\widetilde{q_{ik}}=\widetilde{q_{ij}}\widetilde{q_{kj}}=1$,
$\widetilde{q_{ij}}\neq -1$,
\begin{equation}\label{eqn:relacion vertice -1}
\left[ x_{ijk} , x_j \right]_c;
\end{equation}

\noindent $\diamond$ if $i,j \in \{1, \ldots, \theta \}$ satisfy
$q_{jj}=-1$, $q_{ii}\widetilde{q_{ij}}\in \G_6$,
$\widetilde{q_{ij}}\neq -1$, and also $q_{ii}\in \G_3$ or
$m_{ij}\geq 3$,
\begin{equation}\label{eqn:relacion estandar B2}
\left[ x_{iij}, x_{ij} \right]_c;
\end{equation}

\noindent $\diamond$ if $i,j,k \in \{1, \ldots, \theta \}$ satisfy $q_{ii}=\pm \widetilde{q_{ij}}\in\G_3$, $\widetilde{q_{ik}}=1$, and also
$-q_{jj}=\widetilde{q_{ij}}\widetilde{q_{jk}}=1$ or $q_{jj}^{-1}=\widetilde{q_{ij}}=\widetilde{q_{jk}}\neq -1$,
\begin{equation}\label{eqn:relacion estandar B3}
\left[ x_{iijk} , x_{ij} \right]_c;
\end{equation}

\noindent $\diamond$ if $i,j,k \in \{1, \ldots, \theta \}$ satisfy $\widetilde{q_{ik}}, \widetilde{q_{ij}}, \widetilde{q_{jk}} \neq 1$,
\begin{equation}\label{eqn:relacion triangulo}
x_{ijk}-\frac{1-\widetilde{q_{jk}}}{q_{kj}(1-\widetilde{q_{ik}})}\left[x_{ik},x_j\right]_c-q_{ij}(1-\widetilde{q_{jk}}) \ x_jx_{ik};
\end{equation}

\noindent $\diamond$ if $i,j,k\in\{1,\ldots,\theta\}$ satisfy one of the following situations
\begin{itemize}
 \item[$\circ$] $q_{ii}=q_{jj}=-1$, $\widetilde{q_{ij}}^2= \widetilde{q_{jk}}^{-1}$, $\widetilde{q_{ik}}=1$, or
 \item[$\circ$] $\widetilde{q_{ij}}=q_{jj}=-1$, $q_{ii}= -\widetilde{q_{jk}}^2\in\G_3$, $\widetilde{q_{ik}}=1$, or
 \item[$\circ$] $q_{kk}=\widetilde{q_{jk}}=q_{jj}=-1$, $q_{ii}= -\widetilde{q_{ij}}\in\G_3$, $\widetilde{q_{ik}}=1$, or
 \item[$\circ$] $q_{jj}=-1$, $\widetilde{q_{ij}}=q_{ii}^{-2}$, $\widetilde{q_{jk}}=-q_{ii}^{-3}$, $\widetilde{q_{ik}}=1$, or
 \item[$\circ$] $q_{ii}=q_{jj}=q_{kk}=-1$, $\pm\widetilde{q_{ij}}=\widetilde{q_{jk}}\in\G_3$,
 $\widetilde{q_{ik}}=1$,
\end{itemize}
\begin{equation}\label{eqn:relacion super C3}
\left[ \left[x_{ij}, x_{ijk} \right]_c, x_j \right]_c;
\end{equation}

\noindent $\diamond$ if $i,j,k\in\{1,\ldots,\theta\}$ satisfy $q_{ii}=q_{jj}=-1$, $\widetilde{q_{ij}}^3=\widetilde{q_{jk}}^{-1}$, $\widetilde{q_{ik}}=1$,
\begin{equation}\label{eqn:relacion super G3}
\left[ \left[x_{ij}, \left[x_{ij}, x_{ijk} \right]_c \right]_c, x_j \right]_c;
\end{equation}

\noindent $\diamond$ if $i,j,k\in\{1,\ldots,\theta\}$ satisfy
$q_{jj}=\widetilde{q_{ij}}^2=\widetilde{q_{jk}}\in \G_3$,
$\widetilde{q_{ik}}=1$,
\begin{equation}\label{eqn:relacion super C3 raiz de orden 3}
\left[ \left[ x_{ijk} , x_j \right]_c x_j \right]_c;
\end{equation}

\noindent $\diamond$ if $i,j,k\in\{1,\ldots,\theta\}$ satisfy
$q_{kk}=q_{jj}=\widetilde{q_{ij}}^{-1}=\widetilde{q_{jk}}^{-1}\in
\G_9$, $\widetilde{q_{ik}}=1$, $q_{ii}=q_{kk}^6$
\begin{equation}\label{eqn:relacion fila 18 rango 3}
\left[ \left[ x_{iij} , x_{iijk} \right]_c, x_{ij} \right]_c;
\end{equation}

\noindent $\diamond$ if $i,j,k\in\{1,\ldots,\theta\}$ satisfy
$q_{ii}=\widetilde{q_{ij}}^{-1}\in \G_9$,
$q_{jj}=\widetilde{q_{jk}}^{-1}=q_{ii}^5$, $\widetilde{q_{ik}}=1$,
$q_{kk}=q_{ii}^6$
\begin{equation}\label{eqn:relacion fila 18 rango 3, caso 2}
[\left[x_{ijk}, x_{j} \right]_c, x_k]_c -(1 +
\widetilde{q_{jk}})^{-1}q_{jk} \left[ \left[x_{ijk}, x_{k} \right]_c
, x_{j} \right]_c;
\end{equation}

\noindent $\diamond$ if $i,j,k\in\{1,\ldots,\theta\}$ satisfy
$q_{jj}=\widetilde{q_{ij}}^3=\widetilde{q_{jk}}\in \G_4$,
$\widetilde{q_{ik}}=1$,
\begin{equation}\label{eqn:relacion super G3 raiz de orden 4}
\left[ \left[ \left[ x_{ijk} , x_j \right]_c, x_j \right]_c, x_j
\right]_c;
\end{equation}

\noindent $\diamond$ if $i,j,k\in\{1,\ldots,\theta\}$ satisfy
$q_{ii} = \widetilde{q_{ij}} =-1$, $q_{jj}= \widetilde{q_{jk}}^{-1}
\neq-1$, $\widetilde{q_{ik}}=1$,
\begin{equation}\label{eqn:relacion parecida a super C3}
\left[x_{ij}, x_{ijk} \right]_c;
\end{equation}

\noindent $\diamond$ if $i,j,k \in\{1,\ldots,\theta\}$ satisfy
$q_{ii}= q_{kk} =-1$, $\widetilde{q_{ik}}=1$, $\widetilde{q_{ij}}
\in \G_3$, $q_{jj}= -\widetilde{q_{jk}} = \pm \widetilde{q_{ij}}$,
\begin{equation}\label{eqn:relacion parecida a super C3-bis}
[x_i, x_{jjk}]_c -(1 + q_{jj}^2)q_{kj}^{-1} \left[x_{ijk}, x_{j}
\right]_c - (1 + q_{jj}^2)(1 + q_{jj}) q_{ij} x_j x_{ijk};
\end{equation}

\noindent $\diamond$ if $i,j,k,l\in\{1,\ldots,\theta \}$ satisfy
$q_{jj}\widetilde{q_{ij}}= q_{jj}\widetilde{q_{jk}}=1$, $q_{kk}=-1$,
$\widetilde{q_{ik}}=\widetilde{q_{il}}=\widetilde{q_{jl}}=1$,
$\widetilde{q_{jk}}^2= \widetilde{q_{lk}}^{-1}= q_{ll}$,
\begin{equation}\label{eqn:relacion super C4}
\left[\left[\left[x_{ijkl},x_k\right]_c, x_j \right]_c, x_k \right]_c;
\end{equation}

\noindent $\diamond$ if $i,j,k,l\in\{1,\ldots,\theta \}$ satisfy
$\widetilde{q_{jk}}= \widetilde{q_{ij}}= q_{jj}^{-1}\in
\G_4'\cup\G_6'$, $q_{ii}=q_{kk}=-1$,
$\widetilde{q_{ik}}=\widetilde{q_{il}}=\widetilde{q_{jl}}=1$,
$\widetilde{q_{jk}}^3= \widetilde{q_{lk}}$,
\begin{equation}\label{eqn:relacion super C4 modificada}
\left[\left[x_{ijk},\left[x_{ijkl}, x_k \right]_c\right]_c, x_{jk}
\right]_c;
\end{equation}

\noindent $\diamond$ if $i,j,k,l\in\{1,\ldots,\theta \}$ satisfy
$q_{ll}=\widetilde{q_{lk}}^{-1}=
q_{kk}=\widetilde{q_{jk}}^{-1}=q^2$, $\widetilde{q_{ij}}=
q_{ii}^{-1}=q^3$ for some $q\in \ku^\times$, $q_{jj}=-1$,
$\widetilde{q_{ik}}=\widetilde{q_{il}}=\widetilde{q_{jl}}=1$,
\begin{equation}\label{eqn:relacion super F4-1}
\left[\left[\left[x_{ijk},x_j\right]_c, \left[x_{ijkl},x_j\right]_c
\right]_c, x_{jk} \right]_c;
\end{equation}

\noindent $\diamond$ if $i,j,k,l\in\{1,\ldots,\theta \}$ satisfy one
of the following situations
\begin{itemize}
  \item[$\circ$] $q_{kk}=-1$, $q_{ii}=\widetilde{q_{ij}}^{-1}= q_{jj}^2$,
$\widetilde{q_{kl}}= q_{ll}^{-1}= q_{jj}^3$, $\widetilde{q_{jk}}=
q_{jj}^{-1}$,
$\widetilde{q_{ik}}=\widetilde{q_{il}}=\widetilde{q_{jl}}=1$, or
  \item[$\circ$] $q_{ii}=\widetilde{q_{ij}}^{-1}= -q_{ll}^{-1}=-\widetilde{q_{kl}}$,
$q_{jj}=\widetilde{q_{jk}}=q_{kk}=-1$,
$\widetilde{q_{ik}}=\widetilde{q_{il}}=\widetilde{q_{jl}}=1$,
\end{itemize}
\begin{equation}\label{eqn:relacion super F4-2}
\left[\left[x_{ijkl}, x_j \right]_c, x_k \right]_c-
q_{jk}(\widetilde{q_{ij}}^{-1}-q_{jj})
\left[\left[x_{ijkl},x_k\right]_c, x_j \right]_c;
\end{equation}

\noindent $\diamond$ if $i,j,k\in\{1,\ldots,\theta\}$ satisfy
$\widetilde{q_{jk}}=1$,
$q_{ii}=\widetilde{q_{ij}}=-\widetilde{q_{ik}}\in \G_3$,
\begin{equation}
\left[x_i, \left[ x_{ij},x_{ik} \right]_c
\right]_c+q_{jk}q_{ik}q_{ji} \left[ x_{iik} ,x_{ij} \right]_c+q_{ij}
\, x_{ij} x_{iik};\label{eqn:relacion mji,mjk=2}
\end{equation}

\noindent $\diamond$ if $i,j,k\in\{1,\ldots,\theta\}$ satisfy
$q_{jj}=q_{kk}=\widetilde{q_{jk}}=-1$,
$q_{ii}=-\widetilde{q_{ij}}\in\G_3$, $\widetilde{q_{ik}}=1$,
\begin{equation}\label{eqn:relacion especial rank3}
\left[x_{iijk}, x_{ijk} \right]_c;
\end{equation}

\noindent $\diamond$ if $i,j\in\{1,\ldots,\theta\}$ satisfy $-q_{ii}, -q_{jj}, q_{ii}\widetilde{q_{ij}}, q_{jj}\widetilde{q_{ij}} \neq 1$,
\begin{equation}\label{eqn:relacion mij,mji mayor que 1}
(1-\widetilde{q_{ij}})q_{jj}q_{ji}\left[x_i, \left[ x_{ij}, x_j \right]_c \right]_c - (1+q_{jj})(1-q_{jj}\widetilde{q_{ij}})x_{ij}^2;
\end{equation}

\noindent $\diamond$ if $i,j\in\{1,\ldots,\theta\}$ satisfy that $m_{ij}\in \{4,5\}$, or $q_{jj}=-1$, $m_{ij}=3$, $q_{ii} \in \G_4$,
\begin{equation}\label{eqn:relacion mij mayor que dos, raiz alta}
\left[x_i,x_{3\alpha_i+2\alpha_j}\right]_c-\frac{1-q_{ii}\widetilde{q_{ij}}-q_{ii}^2\widetilde{q_{ij}}^2q_{jj}}{(1-q_{ii}\widetilde{q_{ij}})q_{ji}}
x_{iij}^2;
\end{equation}

\noindent $\diamond$ if $i,j\in\{1,\ldots,\theta\}$ satisfy $4\alpha_i+3\alpha_j\notin \Delta_+^\chi$, $q_{jj}=-1$ or $m_{ji}\geq2$, and also $m_{ij}\geq 3$,
or $m_{ij}=2$, $q_{ii}\in\G_3$,
\begin{equation}\label{eqn:relacion (m+1)alpha i+m alpha j, caso 2}
x_{4\alpha_i+3\alpha_j}=[x_{3\alpha_i+2\alpha_j}, x_{ij} ]_c;
\end{equation}

\noindent $\diamond$ if $i,j\in\{1,\ldots,\theta\}$ satisfy $3\alpha_i+2\alpha_j\in\Delta_+^\chi$, $5\alpha_i+3\alpha_j\notin\Delta_+^\chi$, and
$q_{ii}^3\widetilde{q_{ij}}, q_{ii}^4\widetilde{q_{ij}}\neq 1$,
\begin{equation}\label{eqn:relacion con 2alpha i+alpha j, caso 2}
[x_{iij}, x_{3\alpha_i+2\alpha_j}]_c;
\end{equation}

\noindent $\diamond$ if $i,j\in\{1,\ldots,\theta\}$ satisfy $4\alpha_i+3\alpha_j\in\Delta_+^\chi$, $5\alpha_i+4\alpha_j\notin\Delta_+^\chi$,
\begin{equation}\label{eqn:relacion (m+1)alpha i+m alpha j, caso 3}
x_{5\alpha_i+4\alpha_j}=[x_{4\alpha_i+3\alpha_j}, x_{ij} ]_c;
\end{equation}

\noindent $\diamond$ if $i,j\in\{1,\ldots,\theta\}$ satisfy $5\alpha_i+2\alpha_j\in\Delta_+^\chi$, $7\alpha_i+3\alpha_j \notin \Delta_+^\chi$,
\begin{equation}\label{eqn:relacion con 2alpha i+alpha j, caso 1}
[[x_{iiij}, x_{iij}], x_{iij} ]_c;
\end{equation}

\noindent $\diamond$ if $i,j\in\{1,\ldots,\theta\}$ satisfy $q_{jj}=-1$, $5\alpha_i+4\alpha_j\in\Delta_+^\chi$,
\begin{equation}\label{eqn:relacion potencia alta}
[x_{iij},x_{4\alpha_i+3\alpha_j}]_c- \frac{b-(1+q_{ii})(1-q_{ii}\zeta)(1+\zeta+q_{ii}\zeta^2)q_{ii}^6\zeta^4} {a\ q_{ii}^3q_{ij}^2q_{ji}^3} x_{3\alpha_i+2\alpha_j}^2,
\end{equation}
where $\zeta=\widetilde{q_{ij}}$, $a=(1-\zeta)(1-q_{ii}^4\zeta^3)-(1-q_{ii}\zeta)(1+q_{ii})q_{ii}\zeta$, $b=(1-\zeta)(1-q_{ii}^6\zeta^5)-a\ q_{ii}\zeta$.
\end{theorem}
\bigskip

In what follows we will use implicitly the isomorphism $\cB(V)\cong\u^+(\chi)$ determined by $x_i \mapsto E_i$; in this way, we identify $\cB(V)$ as a
subalgebra of $\u(\chi)$.

\medskip
For any bicharacter $\chi$ whose root system is finite,
$\cJ^+(\chi)$ denotes the ideal of $\cU^+(\chi)$ generated by all
the relations in Theorem \ref{thm:presentacion minima}, except
\eqref{eqn:potencia raices}, plus the quantum Serre relations
$(\ad_cx_i)^{1-a_{ij}}x_j$ for those vertices such that
$q_{ii}^{a_{ij}}=q_{ij}q_{ji}=q_{ii}$. The last ingredient is to
obtain a quotient of all the algebras $\cU_p(\chi)$, $1\leq p\leq
\theta$.

Call also $\cJ^-(\chi):=\phi_4(\cJ^+(\chi))$,
$\cJ(\chi):=(\cJ^+(\chi)+\cJ^-(\chi))$,
$$ U(\chi):= \cU(\chi)/\cJ(\chi), \qquad U^{\pm}(\chi):= \cU^{\pm}(\chi)/\cJ^{\pm}(\chi).$$

We prove first that $\cJ^+(\chi)$ is contained in the ideal defining the corresponding Nichols algebra. The following Lemma was proved with Agust\'in
Garc\'ia Iglesias and is implicit in other papers.

\begin{lemma}\label{Lema:minimal es primitivo}
Let $I\subset T(V)$ be a braided homogeneous biideal of $T(V)$, so there exists a surjective morphism of braided graded Hopf algebras $\pi:R:=T(V)/I
\to\cB(V)$. Let $\mathbf{x}\in\ker\pi$, $\mathbf{x}\neq 0$ of minimal degree $k\geq2$. Then $\mathbf{x}$ is primitive.
\end{lemma}
\pf As $\pi$ is a morphism of graded braided bialgebras, $\ker\pi$ is a graded biideal:
$$ \Delta(\mathbf{x})= \mathbf{x} \otimes 1+1 \otimes \mathbf{x}+ \sum_{j=1}^n b_j \otimes c_j \in \ker \pi \otimes R + R \otimes \ker \pi, $$
for some homogeneous elements $b_j,c_j\in\bigoplus\limits_{i=1}^{k-1} R_i$, such that $\deg(b_j)+\deg(c_j)=k$. For each $j$ we may assume either
$b_j\in\ker\pi$ or $c_j\in\ker\pi$. If $b_j\in\ker\pi$, then $b_j=0$ by the minimality condition on $k$. Similarly, if $c_j\in\ker\pi$, then $c_j=0$.
Hence $\mathbf{x}$ is primitive in $R$.
\epf

We will work with $\N_0^\theta$-graded ideals, so the following notation will be useful: given $\beta=\sum_i b_i\alpha_i$, $\gamma=\sum_i c_i\alpha_i$, for
some $b_i,c_i\in\N_0$, we say that $\beta\geq\gamma$ (respectively, $\beta>\gamma$) if $b_i\geq c_i$ (respectively, $b_i>c_i$) for all $i\in\{1,\ldots,\theta \}$.

\begin{proposition}\label{prop: proyeccion sobre Nichols}
$\cJ^+(\chi)$ is a braided biideal of $\cU^+(\chi)$, and there exist
a canonical projection of Hopf algebras
$\pi_\chi:U(\chi)\twoheadrightarrow\u(\chi)$ such that
$\pi\left(U^\pm(\chi)\right)=\u^\pm(\chi)$.

Moreover, the multiplication $m:U^+(\chi)\ot U^0(\chi)\ot U^-(\chi)\to U(\chi)$ is an isomorphism of graded vector spaces.
\end{proposition}
\pf We can order the relations according to their $\N$-graduation.
When we quotient by the relations of degree at most $n-1$, the
relations of degree $n$ are primitive by Lemma \ref{Lema:minimal es
primitivo}, because for any of them we can see that the relations in
Theorem \ref{Thm:presentacion} of degree at most $n-1$ are verified.
Moreover, for a relation of degree $\alpha\in\N_0^{\theta}$, it is
enough to verify that the relations of $\N^\theta$-degree lower than
$\alpha$ hold in this partial quotient. For example, each quantum
Serre relation is primitive, and the same holds for $x_i^{N_i}$;
therefore, when we quotient by these relation we have that
$\mathbf{x}= \left[ (\ad_c x_i)^2 x_j , (\ad_c x_i) x_j \right]_c$
is primitive under the conditions for \eqref{eqn:relacion estandar
B2}, because we have quotiented by $x_i^3$, $x_j^2$, so it also
holds that
$$(\ad_c x_i)^3 x_j= (\ad_c x_j)^2 x_i=0.$$
We work in a similar way with the other relations so each partial quotient is a braided bialgebra (and then a Hopf algebra with the induced antipode);
finally, $U^+(\chi)$ is a braided bialgebra, because $\cJ^+(\chi)$ is a braided biideal.

By the definition of Nichols algebra we conclude that $\cJ^+(\chi)
\subseteq I^+(\chi)$. By Lemma \ref{Lemma:caracterizacion ideales
triangulacion}, $\cJ^+(\chi)\cU^0(\chi)\cU^-(\chi)$ is a Hopf ideal
of $\cU(\chi)$, and then the equivalent conditions in Proposition
\ref{prop:desc triangular cocientes} hold. Therefore there exists a
projection of Hopf algebras and a triangular decomposition as in the
Proposition. \epf

Now we prove that the isomorphisms at the beginning of Theorem \ref{thm: iso Lusztig-Heck} induce isomorphisms between the corresponding algebras $U(\chi)$.
The first step is to prove that $T_p(\cJ(\chi))\subset\cJ(\schi)$, which will be proved considering each relation generating the ideal. The following two
Lemmata help us to reduce the number of explicit computations.

\begin{lemma}\label{lema:palabra como combinacion de mayores en dos ideales}
Let $l$ be a Lyndon word such that $[l]_c\cong\sum_{w\in S_{I^+(\chi)},w\succ l}a_w w \ (\mbox{mod }I^+(\chi))$, for some $a_w\in \ku$. Let $I$ be a braided
Hopf ideal $\N^\theta$-graded of $\cU^+(\chi)$ such that the set of good words $S_{I^+(\chi)}$, $S_I$ coincide for those terms $w \succ l$, and assume that
$l$ is written as a linear combination of words greater than $l$ modulo $I$. Then $[l]_c\cong\sum_{w\in S_{I^+(\chi)},w\succ l}a_w w \ (\mbox{mod }I)$.
\end{lemma}
\pf
It is a direct consequence of Corollary \ref{cor:primero}: by this result, $[l]_c$ is written as a linear combination of good hyperwords greater than
$[l]_c$ modulo $I$. Such hyperwords coincide with the corresponding good hyperwords for $I^+(\chi)$ by hypothesis, and also $I\subseteq I^+(\chi)$. Hence
the linear combination should be the same, because the good hyperwords generate a linear complement of the ideal in $\cU^+(\chi)$.
\epf
\medskip

\begin{lemma}\label{lema:dos conjuntos de generadores}
Let $I$ be an $\N_0^\theta$-homogeneous ideal of $T(V)$,
$\theta=\dim V$. Let $S$, $T$ be two minimal sets of homogeneous
generators of $I$. Assume that for each $\alpha \in \N^\theta$ there
exists at most one generator in $S$ (respectively in $T$) of degree
$\alpha$, and denote by $I(S,\alpha)$ (respectively, $I(T,\alpha)$)
the ideal generated by the elements of $S$ (respectively, $T$) of
degree $\beta<\alpha$.

For each $s\in S$ of degree $\alpha\in \N_0^\theta$, there exists $t\in T$ of the same degree, and $c\in\ku^\times$ such that $s\cong c\ t
\left(\mbox{mod }(I(S,\alpha) \right)$, and then $I(S,\alpha)=I(T,\alpha)$.
\end{lemma}
\pf
We prove it by induction on the degree of the generators. Let $s$ be of degree $\alpha$, minimal for the partial order defined on $\N_0^\theta$. Therefore
$\dim I_\alpha=1$, so there exists an element of $T$ which belongs to this subspace of $I$ of dimension $1$.

If the degree of $s$ is not minimal, we apply inductive hypothesis for all the generators in lower degree, so for each $s'$ of lower degree there exists a
corresponding $t'\in T$ of the same degree which satisfies the conditions above, and $I(S,\alpha)=I(T,\alpha)$. Therefore
$$ \dim I(S,\alpha)_\alpha=\dim I(T,\alpha)_\alpha=\dim I_\alpha-1, $$
because $S$ is a minimal set of generators, and by hypothesis there exists a unique generator of degree $\alpha$. As $T$ is also a set of generators of
$I$, there exists $t\in I-I(T,\alpha)=I-I(S,\alpha)$, of degree $\alpha$, so the statement follows.
\epf

\begin{remark}\label{obs:distintos conjuntos de generadores}
This Lemma lets us to identify relations of the same degree for two sets of minimal generators of an ideal, up to relations of lower degree and scalars.
In this way we can consider relations from Theorem \ref{Thm:presentacion} for a fixed order on the letters, and consider relations for another order. If we
have a minimal set and this set contains relations all in different degrees (as we will have for the set of relations of the Nichols algebra or some partial
quotients), then we can find a correspondence as above between the relations of the same $\zt$-degree.

For example, if $q_{ii}^{m_{ij}+1}\neq 1$ for some pair of vertices $i,j$, then the quantum Serre relation $(\ad_cx_i)^{m_{ij}+1}x_j$ is a generator for
the minimal set of generators corresponding to the order $x_i<x_j$, so for the order $x_i>x_j$ we have:
$$ [x_jx_i^{m_{ij}+1}]_c= \left[ \cdots \left[ \left[(\ad_cx_j)x_i, x_i\right]_c,\cdots, \right]_c, x_i \right]_c=a (\ad_cx_i)^{m_{ij}+1}x_j, $$
for some scalar $a\in \ku^\times$.

Also, if $q_{ii}\in \G_3$, $\widetilde{q_{ij}}\in \{\pm q_{ii},-1\}$, $q_{jj}=-1$, we notice that
$$ \left[ (\ad_c x_i)^2x_j, (\ad_cx_i)x_j \right]_c \cong b \left[ (\ad_c x_j)x_i, \left[(\ad_cx_j)x_i, x_i \right]_c \right]_c \ (\mbox{mod }I), $$
for some $b\in\ku^\times$, where $I$ is the ideal generated by $x_i^3$ and $x_j^2$, because such relations correspond to different minimal sets of
generators of the ideal of relations of the Nichols algebra, and these are the generators of degree $3\alpha_i+2\alpha_j$ for each set.
\end{remark}
\medskip

\begin{lemma}\label{lema:corchete con partes triviales}
Let $I$ be a $\zt$-graded ideal of $\cU_p^+(\chi)$. Let $Y,Z \in \cU_p^+(\chi)/I$ be homogeneous elements of degree $\beta, \gamma \in \N_0^\theta$,
respectively, such that $(\ad_c E_p)Y=0$. Then,
\begin{equation}\label{eqn:corchete con elemento que q-conmuta 1}
\left[ (\ad_c E_p)Z, Y \right]_c = (\ad_c E_p) \left[ Z, Y \right]_c.
\end{equation}
If also $\chi(\alpha_p,\beta)\chi(\beta, \alpha_p)=1$, then
\begin{equation}\label{eqn:corchete con elemento que q-conmuta 2}
\chi(\alpha_p,\beta) \left[ Y, (\ad_c E_p)Z \right]_c =  (\ad_c E_p) \left[ Y, Z \right]_c.
\end{equation}
\end{lemma}
\pf
Both identities follow from \eqref{eq:identidad jacobi}. For example, for the second one,
\begin{align*}
(\ad_c E_p) \left[ Y, Z \right]_c &= \left[  E_p,  \left[ Y, Z \right]_c \right]_c= \left[  \left[ E_p,   Y\right]_c , Z  \right]_c +
\chi(\alpha_p, \beta) Y \left[ E_p,Z\right]_c-\chi(\beta,\gamma) \left[ E_p,Z\right]_c Y
\\ &= \chi(\alpha_p, \beta) \left( Y (\ad_c E_p)Z -\chi(\beta,\gamma)\chi(\beta, \alpha_p) (\ad_c E_p)Z Y \right)
\\ &= \chi(\alpha_p,\beta)\left[ Y, (\ad_c E_p)Z \right]_c,
\end{align*}
where we use the condition $\chi(\alpha_p,\beta)\chi(\beta, \alpha_p)=1$.\epf
\medskip

\begin{lemma}\label{lemma:relacion transformada -1}
Let $i,p \in \{ 1, \ldots, \theta \}$ be such that $m_{pi} \geq 2$ and $m_{ip}=1$. Then, in $U(s_p^* \chi)$,
$$ \left[\Eb_{i,m_{pi}}^+, \Eb_{i,m_{pi}-1}^+ \right]_{\cb}= \left[ (\ad_{\cb} \Eb_p)^{m_{pi}}\Eb_i, (\ad_{\cb} \Eb_p)^{m_{pi}-1}\Eb_i \right]_{\cb} =0. $$
\end{lemma}

\begin{remark}
Such relation belongs to the ideal $I^+(\schi)$. In fact, as $2 \alpha_i+\alpha_p \notin \Delta_+^\chi$, we have
$s_p(2 \alpha_i+\alpha_p)=2\alpha_i+(2m_{pi}-1)\alpha_p \notin \Delta_+^\chi$,
so such relation holds in the corresponding Nichols algebra $\u^+(\schi)$.

On the other hand, some of these relations are generators of the ideal $\cJ(\schi)$ by definition, for example \eqref{eqn:relacion estandar B2}. We prove
here that the other relations not in the definition of this ideal are redundant; that is, they are generated by relations of lower degree.
\end{remark}
\pf
We consider the different possible values of $m_{pi}$; we begin with $m_{pi}=2$. Therefore $q_{pp}\in\G_3$ or $q_{pp}^2q_{ip}q_{pi}=1$, and also
$q_{ii}=-1$ or $q_{ii}q_{ip}q_{pi}=1$. If $\underline m_{ip}=1$ for $\schi$, then $p$, $i$ determine a subdiagram of standard type. If
$q_{pp}^3\neq 1$ or $q_{ii}\neq-1$ then $E_p^2E_iE_pE_i$ is written as a linear combination of words greater than $E_p^2E_iE_pE_i$, modulo $\cJ(\schi)$,
using the quantum Serre relations, because in the first case $E_p^2E_iE_p$ appears with non-zero coefficient in $(\ad_c E_p)^3E_i$,
so $E_p^2E_iE_pE_i$ is a linear combination of greater words and $E_p^3E_i^2$, but this last word is in the ideal if $q_{ii}=-1$, or $E_pE_i^2$ appears in
$(\ad_c E_i)^2E_p$ with non-zero coefficient, so in both cases we obtain $E_p^2E_iE_pE_i$ as a combination of greater words. Using Lemma
\ref{lema:palabra como combinacion de mayores en dos ideales}, we conclude that $\left[\Eb_{i,2}^+, \Eb_{i,1}^+ \right]_c=0$. A similar proof in the case
$q_{pp}^3=1$, $q_{ii} \neq -1$ gives us the same conclusion. If $q_{pp}^3=1$, $q_{ii}=-1$, the relation corresponds to \eqref{eqn:relacion estandar B2},
which is by definition a generator of $\cJ(\schi)$.

If $m_{pi}=2$ and $\underline m_{ip}>1$ for $\schi$, then \eqref{eqn:relacion mij,mji mayor que 1} is a generator of $\cJ(\schi)$, and then $E_pE_iE_pE_i$
is a linear combination of $E_p^2E_i^2$ and greater words. Therefore $\left[\Eb_{i,2}^+, \Eb_{i,1}^+ \right]_c \in \cJ(\schi)$, by a similar argument.

If $m_{pi}=3$, then $\underline m_{ip}=1$ for $\schi$, or there exists $\zeta \in \G_{24}$ such that $q_{pp}=\zeta^6$, $q_{ii}^{-1}=q_{pi}q_{ip}=\zeta$.
For the first case, we notice that \eqref{eqn:relacion mij mayor que dos, raiz alta} holds also if $q_{pp}\notin\G_4$, because $\Eb_p^3\Eb_i\Eb_p$ can be written
as a linear combination of other words from the quantum Serre relation $(\ad_{\cb} \Eb_p)^4\Eb_i=0$, and then $\Eb_p^3\Eb_i\Eb_p\Eb_i$ is a linear
combination of greater words multiplying by $\Eb_i$, so we apply Lemma \ref{lema:palabra como combinacion de mayores en dos ideales}; from this relation
we work as above, so we write $E_p^{2}E_iE_p^{2}E_i$ as a linear combination of other words and deduce that $E_p^{3}E_iE_p^{2}E_i$ is a linear combination
of greater words, and we can apply Lemma \ref{lema:palabra como combinacion de mayores en dos ideales} again. For the second case, we write
$\Eb_p^3\Eb_i\Eb_p^2\Eb_i$ as a linear combination of greater words using the quantum Serre relations or the relation
\eqref{eqn:relacion mij mayor que dos, raiz alta}, with the same conclusion.

If $m_{pi}=4,5$, then $\underline m_{ip}=1$ for $\schi$. Therefore we work as before and we obtain the desired relation from
\eqref{eqn:relacion mij mayor que dos, raiz alta} or \eqref{eqn:relacion estandar B2}, according to $3\alpha_p+2\alpha_i$ belongs to $\Delta_+^{\schi}$ or
not. In both situations, we can write $\Eb_p^4\Eb_i\Eb_p^3\Eb_i$ or $\Eb_p^5\Eb_i\Eb_p^4\Eb_i$ as a linear combination of greater words, so we apply Lemma
\ref{lema:palabra como combinacion de mayores en dos ideales} again.
\epf
\medskip

We will prove now that $T_p(x)\in\cJ(\schi)$ for any generator $x$
of the ideal $\cJ^+(\chi)$ so we will have a family of morphisms
between the algebras $U(\chi)$.

\begin{lemma}\label{lemma:potencia raices via Tp}
Let $i$ be a non-Cartan vertex. Then $T_p(E_i^{N_i}) \in \cJ(\schi)$.

\noindent If $i,j$ satisfy $q_{ii}=\widetilde{q_{ij}}=q_{jj}=-1$,
and there exists $k$ such that $\widetilde{q_{ik}}^2\neq-1$ or
$\widetilde{q_{jk}}^2\neq-1$, then $T_p\left(E_{ij}^2 \right) \in
\cJ(\schi)$.
\end{lemma}
\pf Consider the first relation. If $i=p$, then $p$ is not Cartan for $\chi$, so $p$ is not Cartan for $\schi$ too. Therefore, by the definition of the
ideal $\cJ(\schi)$,
$$ T_p(E_p^{N_p})= \Fb_p^{N_p}=\phi_4(\Eb_p^{N_p}) \in \cJ(\schi). $$
We consider then $i \neq p$. In such case, $T_p(E_i^{N_i})=\left( \Eb_{i,m_{pi}}^+ \right)^{N_i}$.

If $m_{pi}=0$, then $\Eb_{i,0}^+=\Eb_i$ and $q_{ip}q_{pi}=1$, so for each $j\neq p$ we have $\qb_{ij}\qb_{ji}=\widetilde{q_{ij}}$. Therefore $i$ is not Cartan
for $\schi$, and $T_p(E_i^{N_i})= \Eb_i^{N_i} \in \cJ(\schi)$.

Consider $m_{pi}\neq0$. As $i$ is not Cartan, there exists $j\neq i$
such that $q_{ii}^{m_{ij}}\widetilde{q_{ij}} \neq 1$.

Assume first that $m_{ip}+1=N_i$. If $m_{ip}=1$, that is $q_{ii}=-1$, there are two possibilities. If $q_{ip}q_{pi}\neq-1$, using Lemma
\ref{lemma:relacion transformada -1}, the identity \eqref{eq:identidad jacobi} and the quantum Serre relation $(\ad_{\cb}\Eb_p)^{m_{pi}+1}\Eb_i=0$,
we compute in $U(s_p^* \chi)$,
\begin{align*}
0&= \left[ \Eb_p, \left[\Eb_{i,m_{pi}}^+, \Eb_{i,m_{pi}-1}^+ \right]_c\right]_{\cb}
\\ &= \left( \schi(\alpha_p,m_{pi}\alpha_p+\alpha_i) - \schi(m_{pi}\alpha_p+\alpha_i,(m_{pi}-1)\alpha_p+\alpha_i) \right)
\left( \Eb_{i,m_{pi}}^+ \right)^2
\\ &= \left( \chi(-\alpha_p,\alpha_i) -\chi(\alpha_i,\alpha_p+\alpha_i) \right) \left( \Eb_{i,m_{pi}}^+ \right)^2=q_{pi}^{-1}(1+q_{ip}q_{pi})
\left( \Eb_{i,m_{pi}}^+ \right)^2,
\end{align*}
so $T_p(E_i^2)=\left(\Eb_{i,m_{pi}}^+ \right)^2\in\cJ(\schi)$. If $q_{pi}q_{ip}=-1$, there are 3 possible subdiagrams determined by $i,p$: it is standard
with $q=-1$, or it is Cartan of type $B_2$ with $q \in \G_4$, or it is Cartan of type $G_2$ with $q \in G_6$. For the first case, if the diagram is of type
$A_2$ associated to $q=-1$, it follows by definition of the ideal $\cJ(\schi)$; for the other cases, we write $\Eb_p^{m_{pi}}\Eb_i\Eb_p^{m_{pi}}\Eb_i$ as
a linear combination of greater words using \eqref{eqn:relacion estandar B2} or the quantum Serre relations, and also the previous Lemmata.

If $m_{ip}>1$ and $m_{pi}=1$, we compute, using \eqref{eq:identidad jacobi} and the relation $(\ad_{\cb} \Eb_p)^2\Eb_i=0$,
\begin{align*}
\ad_{\cb} \Eb_p \left[ \Eb_{i,1}^+, \Eb_i \right]_{\cb}&= \left( \schi(\alpha_p,\alpha_i+\alpha_p) - \schi(\alpha_i+\alpha_p, \alpha_i) \right)
\left(\Eb_{i,1}^+\right)^2
\\ &= q_{pi}^{-1}(1-q_{ii}q_{ip}q_{pi})  \left(\Eb_{i,1}^+\right)^2.
\end{align*}
From this relation and \eqref{eq:identidad jacobi} again, we calculate
\begin{align*}
\ad_{\cb} \Eb_p \left[ \Eb_{i,1}^+, \left[ \Eb_{i,1}^+, \Eb_i \right]_{\cb} \right]_{\cb} = &
\left( \schi(\alpha_p,\alpha_i+\alpha_p) - \schi(\alpha_i+\alpha_p,2\alpha_i+\alpha_p) \right)
\\ &(q_{pi}^{-1}-q_{ii}q_{ip}) \left(\Eb_{i,1}^+\right)^3
\\ = & q_{pi}^{-2}(1-q_{ii}q_{ip}q_{pi})(1-q_{ii}^2q_{ip}q_{pi})  \left(\Eb_{i,1}^+\right)^3.
\end{align*}
So if $m_{ip}=2$ and $m_{pi}=1$, it follows that $\alpha_p+3 \alpha_i \notin \Delta_+^\chi$, and
$$ s_p(\alpha_p+3 \alpha_i)= 3\alpha_i+2 \alpha_p \notin \Delta_+^{\schi}. $$
Using the previous Lemma, $\left[ (\ad_{\cb}\Eb_i)^2\Eb_p,
(\ad_{\cb}\Eb_i)\Eb_p \right]_{\cb}\in \cJ^+(\schi)$, so
$$ \left[ \Eb_{i,1}^+, \left[ \Eb_{i,1}^+, \Eb_i \right]_{\cb} \right]_{\cb} \in \cJ(\schi),$$
because we apply Lemma \ref{lema:dos conjuntos de generadores} if the relation belongs to a minimal set of generators ($q_{ii}^2q_{ip}q_{pi}\neq 1$, so
$q_{ii}\in\G_3$), or we compute it directly for the cases in which it is Cartan of type $B_2$ or standard with $q_{pp}=-1$. Then, by a similar argument,
$$ T_p(E_i^3)= \left(\Eb_{i,1}^+\right)^3 \in \cJ(\schi). $$

If $m_{ip}=3$, $m_{pi}=1$, we have that $s_p(\alpha_p+4\alpha_i)=4\alpha_i+3\alpha_p\notin\Delta_+^{\schi}$, so
$$ \left[ \Eb_{i,1}^+, \left[ \Eb_{i,1}^+, \left[ \Eb_{i,1}^+, \Eb_i \right]_{\cb} \right]_{\cb}\right]_{\cb}\in \cJ(\schi),$$
by a similar argument, using \eqref{eqn:relacion (m+1)alpha i+m alpha j, caso 2}. In this case we deduce that $\left(\Eb_{i,1}^+\right)^4 \in \cJ(\schi)$.

If $m_{ip}=4$, then $q_{ii}^4q_{ip}q_{pi}\neq 1$, so $\left(\Eb_{i,1}^+\right)^5 \in \cJ(\schi)$ in a similar way, using
\eqref{eqn:relacion (m+1)alpha i+m alpha j, caso 3}. We notice that there are no diagrams such that $q_{ii}^{m_{ip}+1}=1$ and $m_{ip}\geq5$.

Now we consider $m_{ip},m_{pi}>1$, so there are 3 possibilities:
\begin{itemize}
  \item $m_{ip}=m_{pi}=2$, so \eqref{eqn:relacion (m+1)alpha i+m alpha j, caso 2} is a generator of $\cJ(\schi)$, and $q_{pp}\in\G_3$.
Therefore we write $\Eb_i\Eb_p^2\Eb_i\Eb_p^2\Eb_i$ as a linear combination of other words, which begin with $\Eb_p$ or they contain $\Eb_p^3$ as a factor.
If we multiply by $\Eb_p^2$ on the left, $\Eb_p^2\Eb_i\Eb_p^2\Eb_i\Eb_p^2\Eb_i$ is a linear combination of greater words modulo $\cJ(\schi)$, because
$\Eb_p^3\in\cJ(\schi)$, so $T_p(E_i^3)=( (\ad_c\Eb_p)^2\Eb_i)^3\in \cJ(\schi)$.
\smallskip
  \item $m_{ip}=3$, $m_{pi}=2$, so $q_{ii}=\zeta^6$, $q_{pp}=\zeta^8$, $q_{ip}q_{pi}=\zeta^{11}$ for some $\zeta\in\G_{24}$, and
\eqref{eqn:relacion potencia alta} is a generator of the ideal $\cJ(\schi)$. Using this relation we write $\Eb_i\Eb_p^2\Eb_i\Eb_p^2\Eb_i\Eb_p^2\Eb_i$ as a
linear combination of words beginning with $\Eb_p$ or words containing $\Eb_p^3$ as a factor. Multiplying by $\Eb_p^2$ on the left, we write
$(\Eb_p^2\Eb_i)^4$ as a linear combination of greater words modulo $\cJ(\schi)$, because $\Eb_p^3\in\cJ(\schi)$, so
$T_p(E_i^4)=( (\ad_c\Eb_p)^2\Eb_i)^4\in \cJ(\schi)$.
\smallskip
  \item $m_{ip}=2$, $m_{pi}=3$, so there are two possible diagrams; in both cases \eqref{eqn:relacion con 2alpha i+alpha j, caso 1} is a generator of
$\cJ(\schi)$. From this relation we write $\Eb_i\Eb_p^3\Eb_i\Eb_p^3\Eb_i$ as a combination of words beginning with $\Eb_p$ o containing $\Eb_p^4$.
Multiplying on the left by $\Eb_p^3$, $\Eb_p^3\Eb_i\Eb_p^3\Eb_i\Eb_p^3\Eb_i$ can be written as a linear combination of greater words, modulo $\cJ(\schi)$,
because $\Eb_p^4\in\cJ(\schi)$ or $(\ad_{\cb}\Eb_p)^4\Eb_i \in \cJ(\schi)$, so, in both cases, $T_p(E_i^3)=( (\ad_c\Eb_p)^3\Eb_i)^3\in \cJ(\schi)$.
\end{itemize}

Finally we consider $q_{ii}^{m_{ip}}q_{ip}q_{pi}=1$, $m_{ip}<N_i-1$, so there exists $j \neq p$ such that $1 \leq m_{ip}<m_{ij}=N_i-1$. In this case, $i$,
$j$, $p$ determine a connected diagram, where $i$ is not Cartan, connected with $j$ and $p$, and also $q_{ii}$ is a root of unity of order $N_i>2$.
We have the following possible diagrams under the previous conditions:
\begin{itemize}
  \item $q_{ii}\in \G_3$, $q_{pp}\in \{q_{ii},-1\}$, $m_{ip}=m_{pi}=1$, $m_{ij}=2$, $m_{pj}=m_{jp}=0$, or
  \item $q_{ii}\in \G_4$, $q_{pp}=-1$, $q_{ip}q_{pi}q_{ii}=1$, $\widetilde{q_{ij}}=q_{ii}$, $m_{pj}=m_{jp}=0$, $m_{ij}=3$ (a diagram of type super $G(3)$,
with $q \in \G_4$).
\end{itemize}
Both possibilities follow in a similar way to the case $m_{pi}=1$.

\bigskip
We analyze now the second relation. As $\widetilde{q_{ij}}=-1$,
$$ \left((\ad_c E_j)E_i \right)^2 = q_{ji}^2\left((\ad_c E_i)E_j \right)^2+2q_{ji}(E_iE_j^2E_i + E_jE_i^2E_j). $$
By the first part $T_p(E_i^2),T_p(E_j^2)\in\cJ(\schi)$, and as $T_p$ is an algebra morphism, it is enough to prove that
$T_p\left(\left((\ad_c E_i)E_j\right)^2\right)\in\cJ(\schi)$, to conclude that also $T_p\left(\left((\ad_c E_j)E_i \right)^2 \right) \in \cJ(\schi)$, and
vice versa. Moreover we need just one of these two relations in order to have a minimal set of relations.

If $p=j$, we have
\begin{align*}
T_p\left(\left((\ad_c E_i)E_p \right)^2 \right) &= \left( \Eb_{i,1}^+\Fb_p\Lb_p^{-1}- q_{ip}\Fb_p\Lb_p^{-1}\Eb_{i,1}^+ \right)^2
\\ &= \left( \left( \Fb_p\Eb_{i,1}^-2\Lb_p\Eb_i \right) \Lb_p^{-1}- q_{ip}\qb_{pp}\qb_{ip}\Fb_p\Eb_{i,1}^+\Lb_p^{-1} \right)^2
\\ &=\left(-2\qb_{ip}^{-1}\Eb_i \right)^2= 4q_{ip}^2 \Eb_i^2 \in \cJ(\schi),
\end{align*}
because $\qb_{pp}=\qb_{pi}\qb_{ip}=\qb_{ii}=-1$.

Now we consider $p\neq i,j$. If $m_{pi},m_{pj}\neq0$, we have two
possibilities:
\begin{itemize}
  \item $q_{pp}=-1$, $q_{ip}q_{pi}q_{jp}q_{pj}=-1$ so it is a diagram of type super $D(2,1;\alpha)$), or
  \item $q_{pp}=q_{ip}^{-1}q_{pi}^{-1}=-q_{jp}q_{pj} \in \G_3 \cup \G_4 \cup \G_6$.
\end{itemize}
For the first case, or the second when $q_{pp}\in \G_4$,
$$ T_p\left(\left((\ad_c E_i)E_j \right)^2 \right)= \left[ (\ad_{\cb} \Eb_p)\Eb_i, (\ad_{\cb} \Eb_p)\Eb_j \right]_{\cb}^2.$$
Using \eqref{eqn:relacion triangulo} and $\Eb_p^2$ if $\qb_{pp}=-1$, or the quantum Serre relations
$$(\ad_{\cb}\Eb_p)^2 \Eb_i=(\ad_{\cb}\Eb_p)^2 \Eb_j=0$$
if $q_{pp}\in \G_4$, $\Eb_i\Eb_p\Eb_j\Eb_p$ is a linear combination of greater words, so $(\Eb_p\Eb_i\Eb_p\Eb_j)^2$ is also a linear combination of greater
words. Then,
$$ \left[ (\ad_{\cb} \Eb_p)\Eb_i, (\ad_{\cb} \Eb_p)\Eb_j \right]_{\cb}^2 \in  \cJ^+(\schi), $$
by an analogous statement to Lemma \ref{lema:palabra como combinacion de mayores en dos ideales} but for powers of hyperwords, and such relation is in
$I^+(\schi)$.

For the remaining cases, $q_{pp}\in \G_3 \cup \G_6$ and
$$ T_p\left(\left((\ad_c E_i)E_j \right)^2 \right)= \left[ (\ad_{\cb} \Eb_p)\Eb_i, (\ad_{\cb} \Eb_p)^2\Eb_j \right]_{\cb}^2.$$
We write $(\Eb_p\Eb_i\Eb_p^2\Eb_j)^2$ as a linear combination of greater words using the quantum Serre relations or \eqref{eqn:relacion triangulo},
so
$$T_p\left(\left((\ad_c E_i)E_j \right)^2 \right) \in  \cJ^+(\schi)$$
by an analogous argument.

If $m_{pi}=1$, $m_{pj}=0$, we have two possibilities. If
$q_{pp}=-1$, then $\widetilde{q_{pi}}\neq -1$ and
$\qb_{ii}=\widetilde{q_{pi}}=\widetilde{\qb_{pi}}^{-1}$, so
\eqref{eqn:relacion parecida a super C3} is a generator of
$\cJ(\schi)$, as well as $\Eb_p^2$ and $\Eb_{pj}$. By
\eqref{eq:identidad jacobi}, we have that
$$ q_{pi}q_{pj}(1+\widetilde{q_{pi}})\Eb_{pij}^2 =[\Eb_{ppij},\Eb_{ij}]_{\cb}- [\Eb_p,[\Eb_{pij},\Eb_{ij}]_{\cb}]_{\cb}\in\cJ(\schi),$$
so $T_p(E_{ij}^2)= \Eb_{pij}^2\in\cJ(\schi)$. If not, then
$q_{pp}^{-1}=\widetilde{q_{pi}}\neq-1$, and we write
$\Eb_p\Eb_i\Eb_j\Eb_p\Eb_i\Eb_j$ as a linear combination of greater
words modulo $\cJ(\schi)$, using the quantum Serre relations
(observe that $(\qb_{ij})$ is twist equivalent to the original
braiding). Therefore $T_p(E_{ij}^2)= \Eb_{pij}^2\in\cJ(\schi)$.

If $m_{pi}>1$, $m_{pj}=0$, then $q_{pp}=-\widetilde{q_{pi}}\in\G_3$,
and the proof follows in a similar way to the case $q_{pp}=-1$, but
using the relation \eqref{eqn:relacion especial rank3}.

If $m_{pi}=m_{pj}=0$, the proof follows easily, because
$\qb_{ii}=\qb_{ij}\qb_{ji}=\qb_{jj}=-1$, and $T_p\left(E_{ij}^2
\right)=\Eb_{ij}^2 \in \cJ(\schi)$ by definition of the ideal.\epf
\medskip

\begin{lemma}\label{lemma:quantum serre via Tp}
Let $i,j\in\{1,\ldots,\theta\}$ be such that
$q_{ii}^{m_{ij}}\widetilde{q_{ij}}=1$. Then
$$T_p\left( (\ad_c E_i)^{m_{ij}+1}E_j \right) \in \cJ(\schi).$$
\end{lemma}
\pf
\vi The case $p=i$ was considered in the first part of Theorem \ref{thm: iso Lusztig-Heck}.
\medskip

\noindent \vii Let $p=j$: we analyze all the possible values of $m_{ip}$. If $m_{ip}=0$, then $q_{ip}q_{pi}=1$, and
$$ T_p\left((\ad_c E_i)E_p\right)=\Eb_i\Fb_p\Lb_p^{-1}-q_{ip}\Fb_p\Lb_p^{-1}\Eb_i=(\Eb_i\Fb_p-\Fb_p\Eb_i)\Lb_p^{-1}\in\cJ(\schi).$$

\textsc{Consider $m_{ip}=1$}; by \eqref{eqn:corchete E+ con Fp} we have
\begin{align}
T_p\left( (\ad_c E_i)E_p \right) =& \Eb_{i,m_{pi}}^+ \Fb_p\Lb_p^{-1}-q_{ip}\Fb_p\Lb_p^{-1}\Eb_{i,m_{pi}}^+ \nonumber
\\ =& \left(\Fb_p\Eb_{i,m_{pi}}^+ + (m_{pi})_{\qb_{pp}} (\qb_{pp}^{m_{pi}-1}\qb_{pi}\qb_{ip}-1)\Lb_p \Eb_{i,m_{pi}-1}^+ \right)\Lb_p^{-1} \nonumber
\\ &- q_{ip} s_p^+\chi(m_{pi}\alpha_p+\alpha_i, \alpha_p )\Fb_p\Eb_{i,m_{pi}}^+\Lb_p^{-1} \nonumber
\\ =& (m_{pi})_{\qb_{pp}} (\qb_{pp}^{m_{pi}-1}\qb_{pi}\qb_{ip}-1)\schi\left((m_{pi}-1)\alpha_p+\alpha_i,\alpha_p \right)^{-1} \Eb_{i,m_{pi}-1}^+  \nonumber
\\ & + \Fb_p\Eb_{i,m_{pi}}^+\Lb_p^{-1} - q_{ip} \chi(\alpha_i, -\alpha_p )\Fb_p\Eb_{i,m_{pi}}^+\Lb_p^{-1} \nonumber
\\ \label{eqn:Tp adEi Ep} =& (m_{pi})_{q_{pp}} (q_{pp}^{-1-m_{pi}}q_{pi}^{-1}q_{ip}^{-1}-1)q_{ip}q_{pp} \Eb_{i,m_{pi}-1}^+.
\end{align}
If $m_{pi}=1$, we have by this identity and Remark \ref{rem:corchete via Tp}:
$$ T_p\left( (\ad_c E_i)^2 E_p \right) = T_p\left( \left[ E_i, (\ad_c E_i)E_p \right]_c \right)= \left[ (\ad_{\cb} \Eb_p)\Eb_i, a_1\Eb_i \right]_{\cb}, $$
where $a_{m_{pi}}:=(m_{pi})_{q_{pp}}(q_{pp}^{-1-m_{pi}}q_{pi}^{-1}q_{ip}^{-1}-1)q_{ip}q_{pp}\in\ku^\times$. This element is in $\cJ(\schi)$
because $\underline m_{ip}=1$, so $(\ad_{\cb}\Eb_i)^2\Eb_p\in\cJ(\schi)$. We consider now $m_{pi}\geq2$; by Lemma \ref{lemma:relacion transformada -1},
$$ T_p\left( (\ad_c E_i)^2 E_p \right) = \left[ \Eb_{i,m_{pi}}^+, a_{m_{pi}}\Eb_{i,m_{pi}-1}^+ \right]_{\cb} \in \cJ(\schi). $$

\textsc{Consider now $m_{ip}=2$}. We look at all the possible
diagrams with two vertices and note that $m_{pi}=1$, or there exists
$\zeta \in \G_9$ such that $q_{ii}=-\zeta$, $q_{ip}q_{pi}=\zeta^7$,
$q_{pp}=\zeta^3$. In the first case, $q_{pp}\in \{ -1,q_{ii}^2 \}$,
so this diagram is standard of type $B_2$, and \eqref{eqn:relacion
estandar B2} belongs to $\cJ(\schi)$ by Lemma \ref{lemma:relacion
transformada -1}. Therefore
$$ T_p\left( (\ad_c E_i)^3 E_p \right) = a_1\left[ (\ad_{\cb}\Eb_p)\Eb_i, \left[ (\ad_{\cb}\Eb_p)\Eb_i, \Eb_i \right]_{\cb} \right]_{\cb} \in \cJ(\schi).$$
For the second case, the braiding matrix of $\schi$ is $\qb_{ii}=-1$, $\qb_{ip}\qb_{pi}=\zeta^8$, $\qb_{pp}=\zeta^3$. Then
$$ T_p\left((\ad_c E_i)^3E_p\right)=\left[(\ad_{\cb}\Eb_p)^2\Eb_i,\left[(\ad_{\cb}\Eb_p)^2\Eb_i,(\ad_{\cb}\Eb_p)\Eb_i\right]_{\cb}\right]_{\cb} \in \cJ(\schi), $$
because \eqref{eqn:relacion con 2alpha i+alpha j, caso 2} is a generator of this ideal.

\textsc{It remains to consider $m_{ip}\in\{3,4,5\}$}. The unique diagram with $m_{pi}>1$ is
$$ \xymatrix{ \circ^{-\zeta} \ar@{-}[rr]^{-\zeta^{12}} & & \circ^{\zeta^5} }, $$
where $\zeta \in \G_{15}$, $q_{ii}=-\zeta$ and $m_{ip}=3$, $m_{pi}=2$; applying $s_p$ we obtain
$$ \xymatrix{ \circ^{-1} \ar@{-}[rr]^{-\zeta^{13}} & & \circ^{\zeta^5} }. $$
By \eqref{eqn:relacion (m+1)alpha i+m alpha j, caso 3} we write $\Eb_i\Eb_p^2\Eb_i\Eb_p^2\Eb_i\Eb_p\Eb_i$ as a linear combination of words beginning with
$\Eb_p$, or containing $\Eb_p^3$ as a factor, or greater than this word for the order $p<i$. Multiplying on the left by $\Eb_p^2$ and using that
$\Eb_p^3\in\cJ(\schi)$, $(\Eb_p^2\Eb_i)^3\Eb_p\Eb_i$ can be written as a linear combination of greater words, modulo $\cJ(\schi)$. By Lemma
\ref{lema:palabra como combinacion de mayores en dos ideales} we conclude that
$$ T_p\left( (\ad_c E_i)^4 E_p \right) = \left[  \Eb_{i,2}^+, \left[  \Eb_{i,2}^+, \left[ \Eb_{i,2}^+ , \Eb_{i,1}^+ \right]_{\cb} \right]_{\cb} \right]_{\cb}
=[(\Eb_p^2\Eb_i)^3\Eb_p\Eb_i]_{\cb}\in \cJ(\schi). $$
Finally we consider $m_{pi}=1$, so we have
\begin{align*}
s_p\left( m_{ip}\alpha_i+\alpha_p \right)&= m_{ip}\alpha_i+(m_{pi}-1)\alpha_p \in \Delta_+^{\schi},
\\  s_p\left( (m_{ip}+1)\alpha_i+\alpha_p \right)&= (m_{ip}+1)\alpha_i+m_{pi}\alpha_p \notin \Delta_+^{\schi}.
\end{align*}
If $q_{pp}=\qb_{pp}\neq-1$ then $m_{pi}=3$ and $(\Eb_p\Eb_i)^2\Eb_p\Eb_i^2$ is a linear combination of greater words modulo $\cJ(\schi)$, where we use
first the quantum Serre relation $(\ad_c\Eb_p)^2\Eb_i=0$ to write $\Eb_p\Eb_i\Eb_p$ as a combination of the words $\Eb_p^2\Eb_i$, $\Eb_i\Eb_p^2$ and then
\eqref{eqn:relacion mij mayor que dos, raiz alta}, which also holds in $U(\schi)$. By this relation,
$$ T_p\left((\ad_c E_i)^4 E_p\right)=\left[\Eb_{i,1}^+,\left[\Eb_{i,1}^+,\left[\Eb_{i,1}^+,\Eb_i\right]_{\cb}\right]_{\cb}\right]_{\cb}\in\cJ(\schi).$$
In other case, $q_{pp}=-1$ and $m_{ip}\in\{3,4,5\}$, so we also have that
$$ T_p\left( (\ad_c E_i)^{m_{ip}+1}E_p \right)= [\Eb_{m_{pi}\alpha_i+(m_{pi}-1)\alpha_p}, (\ad_{\cb} \Eb_p) \Eb_i ]_c \in \cJ(\schi),$$
by \eqref{eqn:relacion (m+1)alpha i+m alpha j, caso 2}, \eqref{eqn:relacion (m+1)alpha i+m alpha j, caso 3} or \eqref{eqn:relacion potencia alta},
depending on the value of $m_{ip}$.
\medskip

\viii Let $p\neq j$: if $m_{pi}=0$ (i.e. $\widetilde{q_{ip}}=1$),
then $\qb_{ii}=q_{ii}$, $\qb_{ij}\qb_{ji}=\widetilde{q_{ij}}$, so
$\underline m_{ij}=m_{ij}$, and $(\ad_{\cb} \Eb_i)^{m_{ij}+1}
\Eb_j=0$ holds in $U(\schi)$. Moreover, in $U(\schi)$ we have
$\Eb_p\Eb_i=\qb_{pi} \Eb_p \Eb_i$, so
$$ (\ad_{\cb} \Eb_i)(\ad_{\cb} \Eb_p)X= \qb_{ip} (\ad_{\cb} \Eb_p)(\ad_{\cb} \Eb_i)X ,$$
for any $X\in U(\schi)$, by the second part of Lemma
\ref{lema:corchete con partes triviales}. By Remark
\ref{rem:corchete via Tp} and the previous results, in $U(\schi)$ we
have
\begin{align*}
T_p \left( (\ad_c E_i)^{m_{ij}+1} E_j \right) & = (\ad_{\cb} \Eb_i)^{m_{ij}+1}(\ad_{\cb} \Eb_p)^{m_{pj}} \Eb_j
\\ & = \qb_{ip}^{m_{pj}(m_{ij}+1)} (\ad_{\cb} \Eb_p)^{m_{pj}} (\ad_{\cb} \Eb_i)^{m_{ij}+1} \Eb_j =0.
\end{align*}

Consider now $m_{pi}\neq 0$. If $m_{ij}=m_{pj}=0$, we apply Lemma \ref{lema:corchete con partes triviales} to obtain
$$ T_p\left((\ad_c E_i)E_j\right)=\left[(\ad_{\cb}\Eb_p)^{m_{pi}}\Eb_i,\Eb_j\right]_{\cb}=(\ad_{\cb}\Eb_p)^{m_{pi}}\left(\left[\Eb_i,\Eb_j\right]_{\cb}\right) =0.$$
where we use that $\qb_{ij}\qb_{ji}=\widetilde{q_{ij}}=1$, so in
$U(\schi)$ it holds that $\left[\Eb_i,\Eb_j\right]_{\cb}=0$. It
remains to consider the case in which $i$, $j$ and $p$ determine a
connected diagram, and $m_{pi}\neq 0$.

\textsc{First we analyze the case $m_{ij}=0$, $m_{pj} \neq 0$}. If
$q_{pp}=-1$ it follows that $m_{pi}=m_{pj}=1$. Then
$\qb_{ij}\qb_{ji}=q_{ip}q_{pi}q_{jp}q_{pj}$, and
$\Eb_p\Eb_i\Eb_p\Eb_j$ is a linear combination of greater words for
the order $p<i<j$, modulo $\cJ(\schi)$:
\begin{itemize}
  \item if $\qb_{ij}\qb_{ji}=1$, it follows from \eqref{eqn:relacion vertice -1},
  \item if $\qb_{ij}\qb_{ji}\neq 1$, we write $\Eb_i\Eb_p\Eb_j$ as a linear combination of other words by \eqref{eqn:relacion triangulo}, where those words
are greater than $\Eb_i\Eb_p\Eb_j$ or begin with $\Eb_p$, so we multiply on the left by $\Eb_p$ and use that $\Eb_p^2\in\cJ(\schi)$.
\end{itemize}
In this way, $T_p\left((\ad_c E_i)E_j
\right)=\left[(\ad_{\cb}\Eb_p)\Eb_i,(\ad_{\cb}\Eb_p)\Eb_j\right]_{\cb}\in\cJ(\schi)$.
If $m_{pi}=m_{pj}=1$ and $q_{pp}\neq-1$, then
$(\ad_{\cb}\Eb_p)^2\Eb_i,(\ad_{\cb}\Eb_p)^2\Eb_j \in\cJ(\schi)$; by
these relations and $(\ad_{\cb}\Eb_i)\Eb_j$, $\Eb_p\Eb_i\Eb_p\Eb_j$
can be written as a linear combination of greater words for the
order $p<i<j$, modulo $\cJ(\schi)$, and also
$T_p\left((\ad_cE_i)E_j\right)\in\cJ(\schi)$ in this case.

If $m_{pj}=1$ and $m_{pi}>1$ (or analogously, $m_{pj}>1$, $m_{pi}=1$), then $q_{pp}q_{pj}q_{jp}=1$, and $q_{pp}\neq-1$. Note that $$\widetilde{q_{ij}}=\schi(\alpha_i,\alpha_j)\schi(\alpha_j,\alpha_i)=q_{pp}^{m_{pi}}q_{pi}q_{ip}. $$
If $q_{pp}^{m_{pi}}q_{pi}q_{ip}\neq 1$, then \eqref{eqn:relacion triangulo} holds in $U(\schi)$, so we can write $\Eb_i\Eb_p\Eb_j$ as a linear combination of other words, greater than $\Eb_i\Eb_p\Eb_j$ for the order $p<i<j$, or beginning with $\Eb_p$. Multiplying on the left by $\Eb_p^{m_{pi}}$, we express  $\Eb_p^{m_{ij}}\Eb_i\Eb_p\Eb_j$ as a linear combination of greater words, using that $\Eb_p^{m_{pi}+1}\in\cJ(\schi)$, or $(\ad_{\cb}\Eb_p)^{m_{pi}+1}\Eb_i\in \cJ(\schi)$, so
$$ T_p\left((\ad_cE_i)E_j\right)= \left[(\ad_{\cb}\Eb_p)^{m_{pi}}\Eb_i,(\ad_{\cb}\Eb_p)\Eb_j\right]_{\cb}\in\cJ(\schi). $$
If $q_{pp}^{m_{pi}}q_{pi}q_{ip}=1$ and $q_{pp}^{m_{pi}+1}\neq1$, $\Eb_p^{m_{ij}}\Eb_i\Eb_p\Eb_j$ is written as a linear combination of greater words for the same order using $(\ad_{\cb}\Eb_p)^{m_{pi}+1}\Eb_i$, $(\ad_{\cb}\Eb_p)^2\Eb_j$ and $(\ad_{\cb}\Eb_i)\Eb_j$, so we obtain the same conclusion. If $q_{pp}^{m_{pi}}q_{pi}q_{ip}=1$ and $q_{pp}^{m_{pi}+1}=1$, then $m_{pi}=2$ or $m_{pi}=3$, and the conclusion follows from \eqref{eqn:relacion super C3 raiz de orden 3} or \eqref{eqn:relacion super G3 raiz de orden 4}, respectively.

If $m_{pi},m_{pj}>1$, there is only one possibility:
$m_{pi}=m_{pj}=2$. The proof is as above, expressing
$\Eb_p^2\Eb_i\Eb_p^2\Eb_j$ as a linear combination of greater words
in the two possible cases: if $q_{pp}\notin\G_3$, using the quantum
Serre relations; if $q_{pp}\in\G_3$, by \eqref{eqn:relacion
mji,mjk=2} and $\Eb_p^3$.
\smallskip

\textsc{We consider now $m_{pj}=0$, $m_{ij} \neq 0$}. Note that $m_{ij} \leq 3$, because we have a connected diagram with three vertices and $q_{ii} \neq -1$:
$q_{ii}^{m_{ij}+1}\neq 1$. If $m_{ij}=3$, it corresponds to a diagram of type super $G(3)$:
$$ \chi: \xymatrix{ \circ^{-1} \ar@{-}[r]^{q^{-1}} & \circ^{q} \ar@{-}[r]^{q^{-3}} & \circ^{q^3} }\quad \leftrightsquigarrow_{s_p} \quad
\schi:\xymatrix{ \circ^{-1} \ar@{-}[r]^{q} & \circ^{-1}
\ar@{-}[r]^{q^{-3}} & \circ^{q^3} }.$$ Using \eqref{eqn:relacion
super G3}, $\Eb_i(\Eb_p\Eb_i)^3\Eb_j$ can be written as a linear
combination of other words modulo $\cJ(\schi)$, which are greater
than this word for the order $p<i<j$, or begin with $\Eb_p$ (recall
that $\Eb_i^2\in\cJ(\schi)$). Multiplying on the left by $\Eb_p$,
$(\Eb_p\Eb_i)^4\Eb_j$ is expressed as a linear combination of
greater words, modulo $\cJ(\schi)$, using that
$\Eb_p^2\in\cJ(\schi)$. By Lemma \ref{lema:palabra como combinacion
de mayores en dos ideales},
\begin{align*}
T_p\left((\ad_c E_i)^4 E_j\right)=& \left[(\ad_{\cb} \Eb_p)\Eb_i,\left[(\ad_{\cb}\Eb_p)\Eb_i,\left[(\ad_{\cb} \Eb_p)\Eb_i,
(\ad_{\cb} \Eb_p)(\ad_{\cb} \Eb_i) \Eb_j \right]_{\cb} \right]_{\cb}  \right]_{\cb}
\\ =& [(\Eb_p\Eb_i)^4\Eb_j]_{\cb}\in \cJ(\schi).
\end{align*}
If $m_{ij}=2$, then $m_{pi}=2$ for the diagram
$$\chi: \xymatrix{ \circ^{\zeta^6} \ar@{-}[r]^{\zeta^4} & \circ^{\zeta^5}\ar@{-}[r]^{\zeta^8} & \circ^{\zeta} } \quad \leftrightsquigarrow_{s_p}
\quad \schi:
\xymatrix{\circ^{\zeta^6}\ar@{-}[r]^{\zeta^8}&\circ^{\zeta}
\ar@{-}[r]^{\zeta^8}&\circ^{\zeta}},$$ where $\zeta\in \G_9$, or
$m_{pi}=1$. In the first case, we use \eqref{eqn:relacion fila 18
rango 3}, $\Eb_p^3\in\cJ(\schi)$ and the quantum Serre relations,
and call $c=\schi(4\alpha_p+2\alpha_i+\alpha_j,\alpha_p+\alpha_i)=
\chi(2\alpha_i+\alpha_j,\alpha_p+\alpha_i)$, so
\begin{align*}
 T_p\big((\ad_cE_i)^3 & E_j\big)=\left[(\ad_{\cb}\Eb_p)^2\Eb_i,\left[(\ad_{\cb}\Eb_p)^2\Eb_i,(\ad_{\cb}\Eb_p)^2(\ad_{\cb}\Eb_i)\Eb_j\right]_{\cb}\right]_{\cb}
\\&=(\ad_{\cb}\Eb_p)\left(\left[(\ad_{\cb}\Eb_p)\Eb_i,\left[(\ad_{\cb}\Eb_p)^2\Eb_i,(\ad_{\cb}\Eb_p)^2(\ad_{\cb}\Eb_i)\Eb_j\right]_{\cb}\right]_{\cb}\right)
\\&= -c \, (\ad_{\cb}\Eb_p) \left( \big[[\Eb_{ppi},\Eb_{ppij}], \Eb_{pi} \big]_{\cb} \right)
\in\cJ(\schi).
\end{align*}
If $m_{pi}=1$ and $q_{pp}\neq -1$, using $(\ad_{\cb}\Eb_p)^2\Eb_i$
we write $\Eb_p\Eb_i\Eb_p\Eb_i\Eb_j\Eb_i$ as a linear combination of
greater words and $\Eb_p^2\Eb_i^2\Eb_j\Eb_i$, for the order induced
by $p<i<j$, modulo $\cJ(\schi)$. Using now $(\ad_{\cb}\Eb_p)\Eb_j$,
and $(\ad_{\cb}\Eb_i)^3\Eb_j$ when $\widetilde{q_{ip}}^3\neq 1$, or
\eqref{eqn:relacion super C3 raiz de orden 3} in other case,
$\Eb_p\Eb_i\Eb_p\Eb_i\Eb_j\Eb_i$ is expressed as a linear
combination of greater words modulo $\cJ(\schi)$. If $m_{pi}=1$ and
$q_{pp}=-1$, then \eqref{eqn:relacion super C3} is a generator of
$\cJ(\schi)$, so
$$ \left[\left[(\ad_{\cb}\Eb_p)\Eb_i,(\ad_{\cb}\Eb_p)(\ad_{\cb}\Eb_i)\Eb_j\right]_{\cb}, \Eb_i \right]_{\cb}\in \cJ(\schi),$$
or \eqref{eqn:relacion estandar B2} (considered for the pair $p,i$) is a generator of the
ideal. In any case we have that
\begin{align*}
T_p\left((\ad_c E_i)^3E_j\right)&=\left[(\ad_{\cb}\Eb_p)\Eb_i,\left[(\ad_{\cb}\Eb_p)\Eb_i,(\ad_{\cb}\Eb_p)(\ad_{\cb}\Eb_i)\Eb_j\right]_{\cb}\right]_{\cb}
\\ &=(\ad_{\cb}\Eb_p)\left[\Eb_i,\left[(\ad_{\cb}\Eb_p)\Eb_i,(\ad_{\cb}\Eb_p)(\ad_{\cb}\Eb_i)\Eb_j\right]_{\cb}\right]_{\cb}\in \cJ(\schi).
\end{align*}
Now we fix $m_{ij}=1$. If $m_{pi}=1$, we analyze each different
possible diagram.
\begin{itemize}
 \item If $q_{pp}\neq-1$, then $\schi$ is twist equivalent to $\chi$ (restricted to the vertices $p,i,j$), and $\Eb_p\Eb_i\Eb_p\Eb_i\Eb_j$
can be expressed as a linear combination of greater words modulo $\cJ(\schi)$, using the quantum Serre relations
$$ (\ad_{\cb} \Eb_p)^2 \Eb_i = (\ad_{\cb} \Eb_i)^2 \Eb_j =(\ad_{\cb} \Eb_p)\Eb_j =0. $$
 \item If $q_{pp}= -1$ and $q_{ii}q_{ip}q_{pi}=1$, then $\qb_{ii}=-1$ and $\qb_{ip}\qb_{pi}\qb_{ij}\qb_{ji}=1$. In this way \eqref{eqn:relacion vertice -1} is a generator of the ideal, and by Lemma \ref{lema:corchete con partes triviales},
\begin{align*}
T_p \left( (\ad_c E_i)^2 E_j \right) &= \left[ (\ad_{\cb} \Eb_p)\Eb_i , (\ad_{\cb} \Eb_p)(\ad_{\cb} \Eb_i) \Eb_j \right]_{\cb}
\\ &= (\ad_{\cb} \Eb_p) \left( \left[ \Eb_i , (\ad_{\cb} \Eb_p)(\ad_{\cb} \Eb_i) \Eb_j \right]_{\cb} \right) \in \cJ(\schi).
\end{align*}
 \item If $q_{pp}=-1$ and $q_{ii}q_{ip}q_{pi} \neq  1$, then $\qb_{pj}\qb_{jp}=q_{pj}q_{jp}=1$, $\qb_{ij}\qb_{ji}=\widetilde{q_{ij}}\neq -1$ and
$$ \qb_{ii}^{-1} = q_{pp}^{-1}q_{ip}^{-1}q_{pi}^{-1}q_{ii}^{-1}=- \qb_{ip}\qb_{pi} \qb_{ij}\qb_{ji} \neq -1,$$
so \eqref{eqn:relacion parecida a super C3} or \eqref{eqn:relacion
parecida a super C3-bis} are generators of $\cJ(\schi)$, and then
$T_p \left( (\ad_c E_i)^2 E_j \right) \in \cJ(\schi)$.
\end{itemize}
Now we fix $m_{pi} \neq 1$. The possible connected diagrams of rank
three with these conditions must verify $m_{pi}=2$, $m_{ip}=1$.
Using the quantum Serre relations
$(\ad_{\cb}\Eb_p)\Eb_j=(\ad_{\cb}\Eb_p)^3\Eb_i=0$ if $q_{pp}^3 \neq
-1$, or \eqref{eqn:relacion estandar B3}, \eqref{eqn:relacion fila
18 rango 3, caso 2} depending on the case,
$\Eb_p^2\Eb_i\Eb_j\Eb_p\Eb_i$ can be expressed as a linear
combination of greater words for the order $p<i<j$, and by Lemma
\ref{lema:palabra como combinacion de mayores en dos ideales},
$$ \left[(\ad_{\cb}\Eb_p)^2(\ad_{\cb}\Eb_i)\Eb_j,(\ad_{\cb}\Eb_p)\Eb_i\right]_{\cb} \in \cJ(\schi). $$
By Lemma \eqref{lema:corchete con partes triviales} we conclude that
\begin{align*}
T_p \left( (\ad_c E_i)^2 E_j \right) &= \left[ (\ad_{\cb} \Eb_p)^2\Eb_i , (\ad_{\cb} \Eb_p)^2(\ad_{\cb} \Eb_i) \Eb_j \right]_{\cb}
\\ &= (\ad_{\cb} \Eb_p) \left( \left[ (\ad_{\cb} \Eb_p)\Eb_i , (\ad_{\cb} \Eb_p)^2(\ad_{\cb} \Eb_i) \Eb_j \right]_{\cb} \right) \in \cJ(\schi).
\end{align*}

\smallskip
\textsc{Finally we consider $m_{ij},m_{pj} \neq 0$}, so each pair of vertices is connected. If $m_{ij}=2$, there is just one possibility,
$$ \chi:  \xymatrix{ & \circ^{-1} \ar@{-}[rd]^{q^3} & \\ \circ^{q} \ar@{-}[ru]^{q^{-1}} \ar@{-}[rr]^{q^{-2}} &  & \circ^{-1} } \quad \leftrightsquigarrow_{s_p} \quad \schi: \xymatrix{ \circ^{-1} \ar@{-}[r]^{q} & \circ^{-1} \ar@{-}[r]^{q^{-3}} & \circ^{q^3} }$$
which is a diagram of type super $G(3)$. By \eqref{eqn:relacion super G3} and Lemma \ref{lema:dos conjuntos de generadores} we have that
$$T_p\left((\ad_{\cb}E_i)^3E_j\right)=\left[(\ad_{\cb}\Eb_p)\Eb_i,\left[(\ad_{\cb}\Eb_p)\Eb_i,\left[(\ad_{\cb} \Eb_p)\Eb_i,
(\ad_{\cb}\Eb_p)\Eb_j\right]_{\cb}\right]_{\cb} \right]_{\cb} \in \cJ(\schi).$$
The remaining case is $m_{ij}=1$. If $m_{pi}=m_{pj}=1$, there are two possible cases:
\begin{itemize}
  \item $q_{pp}=-1$; in this case \eqref{eqn:relacion super C3} is a generator of the ideal $\cJ(\schi)$ by definition, and by Lemma
  \ref{lema:dos conjuntos de generadores} we have that
$$T_p\left((\ad_cE_i)^2E_j\right)=\left[(\ad_{\cb}\Eb_p)\Eb_i,\left[(\ad_{\cb}\Eb_p)\Eb_i,(\ad_{\cb}\Eb_p)\Eb_j\right]_{\cb}\right]_{\cb}\in\cJ(\schi).$$
  \item $q_{pp}\neq-1$, $q_{pp}q_{pi}q_{ip}=q_{pp}q_{pj}q_{jp}=1$; $\schi$ is twist equivalent to $\chi$, so \eqref{eqn:relacion triangulo} is a generator of $\cJ(\schi)$. Using also the quantum Serre relations $(\ad_{\cb}\Eb_p)^2\Eb_i$, $(\ad_{\cb}\Eb_p)^2\Eb_j$, $(\ad_{\cb}\Eb_i)^2\Eb_j$, $\Eb_p\Eb_i\Eb_p\Eb_i\Eb_p\Eb_j$ is written as a linear combination of greater words, modulo $\cJ(\schi)$, so as before $T_p\left((\ad_cE_i)^2E_j\right)\in\cJ(\schi)$.
\end{itemize}
It remains to consider the following braiding:
$$\xymatrix{ & \circ^{\zeta} \ar@{-}[rd]^{\zeta^{-1} }& \\ \circ^{-\zeta} \ar@{-}[rr]^{-\zeta^{-1}} \ar@{-}[ru]^{-\zeta^{-1}} & & \circ^{-1} },
\qquad q_{pp}=\zeta \in \G_3, \ m_{ij}=m_{pj}=1, \ m_{pi}=2.$$
The diagram of $\schi$ is $\xymatrix{ \circ^{-1} \ar@{-}[r]^{-1} & \circ^{\zeta} \ar@{-}[r]^{\zeta^{-1}} & \circ^{-1} }$. Then \eqref{eqn:relacion estandar B2} holds for $p$, $i$, and so $\Eb_i\Eb_p^2\Eb_i\Eb_p$ is expressed as a linear combination of other words of the same $\zt$-degree. Multiplying on the left by $\Eb_p^2$, on the right by $\Eb_j$, and using that $\Eb_p^3=0$, $\Eb_p^2\Eb_i\Eb_p^2\Eb_i\Eb_p\Eb_j$ can be written as a linear combination of greater words, so
$$  T_p \left( (\ad_{\cb} E_i)^2 E_j \right) = \left[ (\ad_{\cb} \Eb_p)^2\Eb_i , \left[ (\ad_{\cb} \Eb_p)^2\Eb_i ,(\ad_{\cb} \Eb_p)\Eb_j \right]_{\cb}
\right]_{\cb} \in \cJ(\schi).$$
Therefore we analyze all the cases and the proof is completed.
\epf
\medskip

\begin{lemma}\label{lema:relacion vertice -1}
Let $i,j,k \in \{1,\ldots,\theta\}$ be such that $q_{jj}=-1$, $\widetilde{q_{ik}}=\widetilde{q_{ij}}\widetilde{q_{jk}}=1$. Then,
$$ T_p \left( \left[ E_{ijk} , E_j \right]_c \right) \in \cJ(\schi).$$
\end{lemma}
\pf
Note that in $U(\chi)$ we have the following identity, using the condition on the scalars, \eqref{eq:identidad jacobi} and $(\ad_c E_i)E_k=E_j^2=0$:
$$ \left[E_{ijk},E_j\right]_c = q_{ij}q_{kj} \left[
E_j,E_{ijk}\right]_c = q_{ij}q_{kj} \left[ E_{ji}, E_{jk}
\right]_c,$$
so it is enough to prove that one of these relations is
applied in $\cJ(\schi)$ by $T_p$ for each possible diagram. Let
$p=j$. Note that $\qb_{pp}=-1$,
$\qb_{ip}\qb_{pi}\qb_{pk}\qb_{kp}=\qb_{ik}\qb_{ki}=1$, so
$(\ad_{\cb}\Eb_i)\Eb_k$ and \eqref{eqn:relacion vertice -1} are
generators of $\cJ(\schi)$. By \eqref{eqn:corchete E+ con Fp} we
have:
\begin{align*}
T_p\left(\left[E_{ipk},E_p\right]_c\right)=&\left[a_1\Eb_i,(\ad_{\cb}\Eb_p)\Eb_k\right]_{\cb}\Fb_p\Lb_p^{-1}+
q_{ip}q_{kp}\Fb_p\Lb_p^{-1}\left[a_1\Eb_i,(\ad_{\cb}\Eb_p)\Eb_k\right]_{\cb}
\\ =& a_1\left( \qb_{ip}^{-1}\qb_{pp}^{-1}\qb_{kp}^{-1}+q_{ip}q_{kp}\right)\Fb_p\Lb_p^{-1}\left[\Eb_i,(\ad_{\cb}\Eb_p)\Eb_k\right]_{\cb}
\\ &+ a_1q_{pk}^{-1}(1-q_{pk}q_{kp})\left( \Eb_i\Eb_k+q_{ip}q_{pk}q_{ik}\Eb_k\Eb_i\right)
\\ =& a_1q_{pk}^{-1}(1-q_{pk}q_{kp})\left(\ad_{\cb}\Eb_i\right)\Eb_k \in \cJ(\schi).
\end{align*}
Let $p=i$, which is analogous to the case $p=k$. By \eqref{eqn:corchete con elemento que q-conmuta 1} and \eqref{eqn:Tp adEi Ep},
$$  T_p\left(\left[E_{jp},E_{jk}\right]_c\right)=\left[a_{m_{pj}}(\ad_{\cb}\Eb_p)^{m_{pj}-1}\Eb_j,
(\ad_{\cb} \Eb_p)^{m_{pj}}(\ad_{\cb} \Eb_j)\Eb_k \right]_{\cb}. $$
Note that $m_{pj}=1,2$. If $m_{pj}=1$, $q_{pp}\neq-1$ or $m_{pj}=2$, $q_{pp}\notin\G_3$, then $q_{pp}^{m_{pj}}q_{pj}q_{jp}=1$. In this case,
$-\qb_{jj}=\qb_{pk}\qb_{kp}=\qb_{pj}\qb_{jp}\qb_{kj}\qb_{jk}=1$, so $(\ad_{\cb}\Eb_p)^{m_{pj}+1}\Eb_j=0=(\ad_{\cb}\Eb_p)\Eb_k$. Note also that
\eqref{eqn:relacion vertice -1} is a generator of $\cJ(\schi)$, hence $T_p\left(\left[(\ad_c E_j)E_i,(\ad_c E_j)E_k \right]_c\right)\in\cJ(\schi)$. If
$q_{pp}=-1$, then
$$ \qb_{jj}\qb_{jp}\qb_{pj}=\qb_{jj}\qb_{jk}\qb_{kj}=1, \quad \qb_{pk}\qb_{kp}=1, $$
hence in $U(\schi)$, $(\ad_{\cb} \Eb_j)^2\Eb_p=(\ad_{\cb}\Eb_j)^2 \Eb_k=0$ if $\qb_{jj}\neq-1$, or \eqref{eqn:relacion vertice -1} if $\qb_{jj}=-1$: in
this way,
$$  T_p\left(\left[E_{jp},E_{jk}\right]_c\right)=a_1\left[\Eb_j,(\ad_{\cb}\Eb_p)(\ad_{\cb}\Eb_j)\Eb_k\right]_{\cb}\in\cJ(\schi).$$
The remaining case is $q_{pp}\in\G_3$: by \eqref{eqn:relacion estandar B3} we obtain that
$$T_p\left(\left[E_{jp},E_{jk}\right]_c\right)=a_2\left[(\ad_{\cb}\Eb_p)\Eb_j,(\ad_{\cb}\Eb_p)^2(\ad_{\cb} \Eb_j)\Eb_k \right]_{\cb}
\in \cJ(\schi). $$

Finally take $p\neq i,j,k$. First, the proof is trivial if $p$ is not connected with $i$, $j$, $k$, because in such case $\schi$ is twist equivalent to
$\chi$, and then
$$ T_p\left(\left[E_{ijk},E_j\right]_c\right)=\left[\Eb_{ijk},\Eb_j\right]_{\cb}\in \cJ(\schi).$$
Now, if $p$ is connected just with $i$ (or analogously, just with
$k$), we have
$$ \qb_{jj}=-1, \quad \qb_{ji}\qb_{ij}\qb_{jk}\qb_{kj}=1, \quad \qb_{ik}\qb_{ki}=1, $$
and by Lemma \ref{lema:corchete con partes triviales}:
$$ T_p \left( \left[E_{ijk},E_j\right]_c\right)=\left[(\ad_{\cb}\Eb_p)^{m_{pi}}(\Eb_{ijk}),\Eb_j\right]_{\cb}
= (\ad_{\cb} \Eb_p)^{m_{pi}} \left(\left[\Eb_{ijk},
\Eb_j\right]_{\cb} \right) \in \cJ(\schi).$$ If $p$ is connected
just with $j$, then
$\widetilde{q_{pj}}\in\{\widetilde{q_{ij}},\widetilde{q_{jk}}\}$,
and $m_{pj}=1$. We assume that
$\widetilde{q_{pj}}=\widetilde{q_{ij}}= \widetilde{q_{kj}}^{-1}$. If
$q_{pp}\widetilde{q_{ip}}=1$, $\schi$ is twist equivalent to $\chi$;
in other case, $q_{pp}=-1$, and then
$$ \qb_{jj}= q_{pj}q_{jp}, \qquad \qb_{pj}\qb_{jp}=q_{pj}^{-1}q_{jp}^{-1}=q_{kj}q_{jk}=\qb_{kj}\qb_{jk},$$
so in both cases $\left[(\ad_{\cb}\Eb_p)(\ad_{\cb}\Eb_j)\Eb_k,\Eb_j\right]_{\cb}\in\cJ(\schi)$. Therefore $\Eb_p\Eb_j\Eb_k\Eb_p\Eb_j\Eb_i$ is a linear
combination of greater words (for the order $p<j<k<i$), so we have that
$$ T_p\left(\left[(\ad_c E_j)E_k,(\ad_cE_j)E_i\right]_c\right)=\left[(\ad_{\cb}\Eb_p)(\ad_{\cb}\Eb_j)\Eb_k,(\ad_{\cb} \Eb_p)(\ad_{\cb} \Eb_j)\Eb_i\right]
_{\cb}\in\cJ(\schi). $$

The remaining case is that $p$ is connected with two consecutive
vertices. We can assume that $p$ is connected with $i$ and $j$.
There exist six possible diagrams satisfying these conditions. For
two of them, the diagram of $\schi$ is the same.
$$ \xymatrix{ \circ^{-1} \ar@{-}[r]^{-\zeta} \ar@{-}[rd]_{-1} & \circ^{-1} \ar@{-}[r]^{-\zeta^2} \ar@{-}[d]^{\zeta^2} &
\circ^{-\zeta} \\ & \circ^{q} & }, \, \zeta\in\G_3, \qquad \qquad
\xymatrix{ \circ^{-1} \ar@{-}[r]^{q^2} \ar@{-}[rd]_{q^{-1}} &
\circ^{-1} \ar@{-}[r]^{q^{-2}}\ar@{-}[d]^{q^{-1}} & \circ^{q_{kk}}
\\ & \circ^{q} & }, \, q\neq -1. $$
In the second case, $m_{pi}=m_{pj}=1$. Consider the order $p<j<i<k$:
$$ T_p \left( \left[E_{ji},E_{jk}\right]_c \right) = \left[ \left[
\Eb_{pj}, \Eb_{pi} \right]_{\cb} , \left[\Eb_{pj}, \Eb_k
\right]_{\cb} \right]_{\cb}=
[\Eb_p\Eb_j\Eb_p\Eb_i\Eb_p\Eb_j\Eb_k].$$ Note that $\Eb_j\Eb_p\Eb_i$
is a linear combination of other words of the same degree by
\eqref{eqn:relacion triangulo}, where those words are greater than
$\Eb_j\Eb_p\Eb_i$ or begin with $\Eb_p$. In all the cases we
conclude that $T_p \left( \left[ E_{ji},E_{jk}\right]_c \right) \in
\cJ(\schi)$. The proof for the first case is analogous.

For the remaining four diagrams, we write also the diagram
corresponding to $\schi$:
\begin{align*}
& \xymatrix{ \circ^{-1} \ar@{-}[r]^{q^{-1}} \ar@{-}[rd]_{q^2} & \circ^{-1} \ar@{-}[r]^{q}\ar@{-}[d]^{q^{-1}} & \circ^{q_{kk}} & & \ar@{~>}[rr]^{s_p}
& & & \circ^{q^2}\ar@{-}[rd]_{q^{-2}} & \circ^{q^{-1}} \ar@{-}[r]^{q} \ar@{-}[d]^{q} & \circ^{q_{kk}}\\ & \circ^{-1} &  & &  & & &  & \circ^{-1} & },
\\ & \xymatrix{ \circ^{-1} \ar@{-}[r]^{q} \ar@{-}[rd]_{-q^{-1}} & \circ^{-1} \ar@{-}[r]^{q^{-1}}\ar@{-}[d]^{-1} & \circ^{q_{kk}} & & \ar@{~>}[rr]^{s_p}
& & & \circ^{-q^{-1}}\ar@{-}[rd]_{-q} & \circ^{-1} \ar@{-}[r]^{q^{-1}} \ar@{-}[d]^{-1} & \circ^{q_{kk}}
\\ & \circ^{-1} &  & &  & & &  & \circ^{-1} & },
\\ & \xymatrix{ \circ^{-1} \ar@{-}[r]^{q^2} \ar@{-}[rd]_{q^{-3}} & \circ^{-1} \ar@{-}[r]^{q^{-2}}\ar@{-}[d]^{q} & \circ^{q^2} & &
\ar@{~>}[rr]^{s_p} & & & \circ^{q^{-3}}\ar@{-}[rd]_{q^3} & \circ^{q}
\ar@{-}[r]^{q^{-2}} \ar@{-}[d]^{q^{-1}} & \circ^{q^2}\\ & \circ^{-1}
& & & & & &  & \circ^{-1} & },
\\ & \xymatrix{ \circ^{-1} \ar@{-}[r]^{-\zeta} \ar@{-}[rd]_{-\zeta} & \circ^{-1} \ar@{-}[r]^{-\zeta^2} \ar@{-}[d]^{\zeta} &
\circ^{-\zeta} & & \ar@{~>}[rr]^{s_p} & & &
\circ^{-\zeta}\ar@{-}[rd]_{-\zeta^2} & \circ^{\zeta}
\ar@{-}[r]^{-\zeta^2} \ar@{-}[d]^{\zeta^2} & \circ^{-\zeta}
\\ & \circ^{-1} &  & &  & & &  & \circ^{-1} & }.
\end{align*}
Note that $m_{pi}=m_{pj}=1$. If we fix the order $p<j<i<k$, and
obtain that $ T_p \left( \left[E_{ji},E_{jk}\right]_c \right) =
[\Eb_p\Eb_j\Eb_p\Eb_i\Eb_p\Eb_j\Eb_k]$. We can write
$\Eb_p\Eb_j\Eb_p\Eb_i\Eb_p\Eb_j\Eb_k$ as a linear combination of
greater words modulo $\cJ(\schi)$: in the first case, using
\eqref{eqn:relacion super C3}; for the remaining cases, we use
\eqref{eqn:relacion super F4-2}.
\epf
\medskip

\begin{lemma} \label{lemma: relacion estandar B2}
Let $i,j\in\{1,\ldots,\theta \}$ be such that $q_{jj}=-1$, and also
$q_{ii}=\pm \widetilde{q_{ij}}\in \G_3$, or
$q_{ii}\widetilde{q_{ij}}\in\G_6$. Then,
$T_p\left(\left[(E_{iij},E_{ij}\right]_c\right)\in\cJ(\schi)$, for
any $p\in\{1,\ldots,\theta\}$.
\end{lemma}
\pf We denote $\xb:=\left[(E_{iij},E_{ij}\right]_c$. We begin with
the case $p=j$. Note that $m_{pi}=1$ (because $q_{pp}=-1$,
$\widetilde{q_{pi}}\neq1$),
$3\alpha_i+2\alpha_p\notin\Delta^\chi_+$, so
$$s_p(3\alpha_i+2 \alpha_p) = 3\alpha_i+\alpha_p \notin \Delta^{\schi}_+.$$
Using \eqref{eqn:Tp adEi Ep} and $(\ad_{\cb}\Eb_i)^3\Eb_p\in\cJ(\schi)$, we obtain that
$$ T_p(\xb)=a_2a_1\ \left[\left[\Eb_{pi},\Eb_i\right]_{\cb},\Eb_i\right]_{\cb}\in\cJ(\schi).$$

Now let $p=i$. By Lemma \ref{lema:palabra como combinacion de mayores en dos ideales} it is equivalent to prove that
$$ T_p(\xb')\in \cJ(\schi),\qquad\mbox{where }\xb':=\left[\Eb_{jp},\left[\Eb_{jp},\Eb_p\right]_c\right]_c, $$
because we have proved that $T_p$ apply the generating relations of
degree less than $\xb$ in elements of $\cJ(\schi)$. By
\eqref{eqn:corchete E+ con Fp},
$$T_p(\xb')=a_1^2a_2 \ \left[(\ad_{\cb}\Eb_p)\Eb_j,\Eb_j\right]_{\cb}\in\cJ(\schi),$$
because it holds that $\qb_{jj}=-1$, or $\qb_{jj}\qb_{jp}\qb_{pj}=1$.

Finally, let $p\neq i,j$; the case $m_{pi}=m_{pj}=0$ follows easily as in the previous Lemmata, so consider the case in which $p$, $i$, $j$ determine a
connected subdiagram of rank three. We note that $q_{ii}\in\G_3$.

We take first $m_{pi}\neq 0$, $m_{pj}=0$. The possible braidings verify that $m_{pi}=1$, so for the order $p<i<j$,
$$ T_p(\xb)= \left[ \left[ \Eb_{pi}, \left[ \Eb_{pi}, \Eb_j \right]_{\cb} \right]_{\cb} ,
\left[ \Eb_{pi}, \Eb_j \right]_{\cb} \right]_{\cb}= \left[
\Eb_p\Eb_i(\Eb_p\Eb_i\Eb_j)^2 \right]_{\cb}, $$ where we use that
$(\ad_c\Eb_p)\Eb_j=0$ in $U(\schi)$. As $q_{pi}q_{ip}q_{ii}=1$ or
$q_{pi}q_{ip}= \pm q_{ii}$, we have that $\qb_{ii}=-1$, or
$\qb_{ii}\qb_{ip}\qb_{pi}=1$, or $\qb_{ii}^2\qb_{ip}\qb_{pi}=1$, or
$\qb_{ii}=-\qb_{ip}\qb_{pi}\in \G_3$, so
$\Eb_p\Eb_i\Eb_p\Eb_i\Eb_j\Eb_p\Eb_i\Eb_j$ can be expressed as a
linear combination of greater words modulo $\cJ(\schi)$, using the
quantum Serre relations or \eqref{eqn:relacion estandar B2}. We
deduce that $T_p(\xb) \in \cJ(\schi)$, using the Lemma
\ref{lema:palabra como combinacion de mayores en dos ideales}.

Now let $m_{pi}=0$, $m_{pj}\neq 0$. We note that $m_{pj}=1$ for any possible diagram, and in $U(\schi)$ we have that:
\begin{align*}
T_p(\xb) &= \left[ (\ad_{\cb} \Eb_i)^2(\ad_{\cb} \Eb_p)\Eb_j, (\ad_{\cb} \Eb_i)(\ad_{\cb} \Eb_p)\Eb_j \right]_{\cb}
\\ &= \qb_{ip}^3\left[ (\ad_{\cb} \Eb_p)(\ad_{\cb} \Eb_i)^2\Eb_j, (\ad_{\cb} \Eb_p)(\ad_{\cb} \Eb_i)\Eb_j \right]_{\cb}
\\ &= \qb_{ip}^3(\ad_{\cb} \Eb_p)\left[ (\ad_{\cb} \Eb_i)^2\Eb_j, (\ad_{\cb} \Eb_p)(\ad_{\cb} \Eb_i)\Eb_j \right]_{\cb}
\\ &= \qb_{ip}^3(\ad_{\cb} \Eb_p)\left[ (\ad_{\cb} \Eb_i)^2(\ad_{\cb} \Eb_j)\Eb_p, (\ad_{\cb} \Eb_i)\Eb_j \right]_{\cb} =0,
\end{align*}
by applying \eqref{eqn:corchete con elemento que q-conmuta 1}
(because $(\ad_{\cb}\Eb_p)\Eb_i=0$), \eqref{eq:identidad jacobi},
and also that \eqref{eqn:relacion estandar B3} is a generator of
$\cJ(\schi)$, because $\qb_{ii}=q_{ii}\in\G_3$, $\underline
m_{ji}=\underline m_{jp}=1$, and \eqref{eqn:relacion estandar B2} is
another generator by Lemma \ref{lemma:relacion transformada -1}.

Finally, consider $m_{pi},m_{pj}\neq0$. There exists just one
possible braiding: $q_{pp}=-1=q_{pj}q_{jp}$,
$q_{ii}=-\widetilde{q_{ij}}=q_{pi}^{-1}q_{ip}^{-1}$. The diagram of
$\schi$ is $\xymatrix{ \circ^{-1}  \ar@{-}[r]^{-1} & \circ^{-1}
\ar@{-}[r]^{q_{ii}} &\circ^{-1} }$, and the solution is analogous to
the previous case, but now we use \eqref{eqn:relacion super C3}.
\epf
\medskip

\begin{lemma}\label{lemma: relacion estandar B3}
Let $i,j,k\in\{1,\ldots,\theta\}$ be such that $q_{ii}=\pm
\widetilde{q_{ij}}\in\G_3$, $\widetilde{q_{ik}}=1$, and
$q_{jj}\widetilde{q_{ij}}=q_{jj}\widetilde{q_{jk}}=1$ or
$q_{jj}=-1$, $\widetilde{q_{ij}}\widetilde{q_{jk}}=1$. Then, for any
$p \in \{1, \ldots, \theta \}$,
$$T_p\left(\left[ E_{iijk} , E_{ij} \right]_c \right) \in \cJ(\schi).$$
\end{lemma}
\pf We denote $\xb=\left[ E_{iijk} , E_{ij} \right]_c$. We begin
with the case $p=k$. In all the cases we have that $m_{kj}=1$, and
$\schi$ satisfies the same conditions, so \eqref{eqn:relacion
estandar B3} is a generator of $\cJ(\schi)$. Then,
\begin{align*}
T_p(\xb) &\cong \left[ (\ad_{\cb} \Eb_i)^2(a_1 \Eb_j), (\ad_{\cb}
\Eb_i)(\ad_{\cb}\Eb_k)\Eb_j \right]_{\cb}\cong a_1q_{ik}\left[
\Eb_{iij}, \Eb_{kij} \right]_{\cb}
\\&\cong a_1q_{ik} \left[ \Eb_{iijk}, \Eb_{ij} \right]_{\cb} \cong 0 \quad (\mbox{mod }\cJ(\schi)),
\end{align*}
where we apply first \eqref{eqn:Tp adEi Ep}, then \eqref{eqn:corchete con elemento que q-conmuta 1}, $(\ad_{cb} \Eb_k)\Eb_i=0$ for the second line, y finally
\eqref{eq:identidad jacobi} plus the fact that \eqref{eqn:relacion estandar B2} is a generator of $\cJ(\schi)$. The cases $p=i$, $p=j$ are proved in a similar way to the case $p=i$ of previous Lemma.

Finally take $p\neq i,j,k$, and assume that $p$ is connected with at least one of the other vertices; in other case the proof is easy as above. We have two possible cases: $m_{pi}=1$, $m_{pj}=m_{pk}=0$, or $m_{pk}=1$, $m_{pj}=m_{pi}=0$. For the first one,
$$ T_p(\xb) = \left[ \left[ \Eb_{pi},  \left[ \Eb_{pi}, \Eb_{jk} \right]_{\cb}\right]_{\cb},
[\Eb_{pi},\Eb_j]_{\cb} \right]_{\cb},$$ and we have two
possibilities:
\begin{itemize}
 \item if $q_{pp}=-1$, then $q_{ii}\widetilde{q_{ip}}=1$, and $\qb_{ii}=-1$, so \eqref{eqn:relacion super C3}
 is a generator of $\cJ(\schi)$ for the subdiagram determined by $p,i,j$. Therefore $T_p(\xb)\in\cJ(\schi)$.
 \item if $q_{pp}\neq-1$, then $q_{pp}=q_{ip}^{-1}q_{pi}^{-1}=q_{ii}$, so $\schi$ is twist equivalent to $\chi$ and \eqref{eqn:relacion estandar B3} is a generator of $\cJ(\schi)$. Then $T_p(\xb)\in\cJ(\schi)$, because it is obtained after to apply $(\ad_{cp}\Eb_p)^3$ to \eqref{eqn:relacion estandar B3}
and multiply by a non-zero scalar, where we use also the quantum Serre relations involving $\Eb_p$.
\end{itemize}
For the second case, we use $(\ad_{\cb} \Eb_p)\Eb_i,(\ad_{\cb}\Eb_p)\Eb_j\in\cJ(\schi)$ to obtain that
\begin{align*}
T_p(\xb) &\cong \left[ (\ad_{\cb} \Eb_i)^2( \ad_{\cb} \Eb_j) ( \ad_{\cb} \Eb_p)\Eb_k, (\ad_{\cb} \Eb_i)\Eb_j \right]_{\cb}
\\ &\cong \qb_{ip}^2\qb_{jp}\left[ ( \ad_{\cb} \Eb_p)(\ad_{\cb} \Eb_i)^2( \ad_{\cb} \Eb_j) \Eb_k, (\ad_{\cb} \Eb_i)\Eb_j \right]_{\cb}
\\ &\cong \qb_{ip}^2\qb_{jp} (\ad_{\cb}\Eb_p)\left(\left[(\ad_{\cb}\Eb_i)^2(\ad_{\cb}\Eb_j)\Eb_k,(\ad_{\cb}\Eb_i)\Eb_j\right]_{\cb}\right)\cong0
\quad(\mbox{mod }\cJ(\schi)),
\end{align*}
by \eqref{eqn:corchete con elemento que q-conmuta 1} and the fact that \eqref{eqn:relacion estandar B3} is a generator of $\cJ(\schi)$.
\epf
\medskip

\begin{lemma}\label{lema:relacion triangulo}
Let $i,j,k \in \{1, \ldots, \theta \}$ be such that
$\widetilde{q_{ik}}, \widetilde{q_{ij}}, \widetilde{q_{jk}} \neq 1$.
Then, for any $p$,
$$ T_p\left(E_{ijk}-\frac{1-\widetilde{q_{jk}}}{q_{kj}(1-\widetilde{q_{ik}})}\left[E_{ik},E_j\right]_c-q_{ij}
(1-\widetilde{q_{jk}})E_jE_{ik}\right)\in\cJ(\schi).$$
\end{lemma}
\pf Let
$\xb=E_{ijk}-\frac{1-\widetilde{q_{jk}}}{q_{kj}(1-\widetilde{q_{ik}})}\left[E_{ik},E_j\right]_c-q_{ij}
(1-\widetilde{q_{jk}})E_jE_{ik}$. By a direct computation we obtain
the same relation, up to an scalar, if we permute the vertices $i$,
$j$, $k$, where we use that
$\widetilde{q_{ik}}\widetilde{q_{ij}}\widetilde{q_{jk}}=1$, so it is
enough to consider one of these permutations for each $p$.

Consider then $p=k$, which is analogous to take $p=i$ or $p=j$. Note that $\{ m_{pi},m_{pj}\}= \{1,1\}$, or $\{ m_{pi},m_{pj}\}=\{1,2\}$, so we fix $m_{pj}=1$, $m_{pi} \in \{1,2\}$. By \eqref{eqn:Tp adEi Ep},
\begin{align*}
 T_p \left(E_{ijk}\right) &=(q_{pp}^{-1}q_{pj}^{-1}q_{jp}^{-1}-1)q_{jp}q_{pp}
\left[(\ad_{\cb}\Eb_p)^{m_{pi}}\Eb_i,\Eb_j\right]_{\cb},
\\ T_p \left(\left[E_{ik},E_j\right]_c\right) &=(q_{pp}^{-1-m_{pi}}q_{pi}^{-1}q_{ip}^{-1}-1)q_{ip}q_{pp}
\left[(\ad_{\cb} \Eb_p)^{m_{pi}-1}\Eb_i,(\ad_{\cb}
\Eb_p)\Eb_j\right]_{\cb},
\\ T_p \left( E_jE_{ik} \right) &=(q_{pp}^{-1-m_{pi}}q_{pi}^{-1}q_{ip}^{-1}-1)q_{ip}q_{pp}(\ad_{\cb}\Eb_p)\Eb_j(\ad_{\cb}\Eb_p)^{m_{pi}-1}\Eb_i.
\end{align*}
If $m_{pi}=2$, or $m_{pi}=1$, $q_{pp}\neq-1$, then $\qb_{ik}\qb_{ki}, \qb_{ij}\qb_{ji}, \qb_{jk}\qb_{kj} \neq 1$ and we deduce that
$T_p(\xb)\in\cJ(\schi)$ from the fact that \eqref{eqn:relacion triangulo} is a generator of $\cJ(\schi)$, because we can write then $\Eb_p^{m_{pi}}\Eb_i\Eb_j$ as a linear combination of greater words (for the order on the letters $p<i<j$), modulo $\cJ(\schi)$, y apply then Lemma
\ref{lema:palabra como combinacion de mayores en dos ideales}. If $q_{pp}=-1$ then $\qb_{ij}\qb_{ji}=1$, so $(\ad_{\cb}\Eb_i)\Eb_j\in\cJ(\schi)$.
By a direct computation, there exists $a\in\ku^\times$ such that
$$ T_p(\xb)= a(\ad_{\cb}\Eb_p)(\ad_{\cb}\Eb_i)\Eb_j\in \cJ(\schi). $$

Let $p \neq i,j,k$. We note that $p$ is not connected with any of the other vertices (so the proof follows easily as in the previous Lemmata), or $p$
is connected just with one of these vertices. For the last case we can assume that $m_{pi}\neq0$, so the unique possibility is $m_{pi}=m_{ip}=1$. Then $\widetilde{q_{ik}}=\qb_{ik}\qb_{ki}$, $\widetilde{q_{ij}}=\qb_{ij}\qb_{ji}$, $\widetilde{q_{jk}}=\qb_{kj}\qb_{jk}\neq1$, so \eqref{eqn:relacion triangulo} is a generator of  $\cJ(\schi)$. By Lemma \ref{lema:corchete con partes triviales} and the relations $(\ad_{\cb}\Eb_p)\Eb_j=(\ad_{\cb}\Eb_p)\Eb_k=0$, we deduce that $T_p(\xb)$ is obtained, up to a non-zero scalar, after to apply $(\ad_{\cb}\Eb_p)$ to \eqref{eqn:relacion triangulo}, modulo $\cJ(\schi)$, so $T_p(\xb) \in \cJ(\schi)$.
\epf
\medskip

\begin{lemma}\label{lema: relacion superC3}
Let $i,j,k \in \{1, \ldots, \theta \}$ be such that

\vi $q_{ii}=q_{jj}=-1$, $\widetilde{q_{ij}}^2=
\widetilde{q_{jk}}^{-1}$, $\widetilde{q_{ik}}=1$, or

\vii $\widetilde{q_{ij}}=q_{jj}=-1$, $q_{ii}=
-\widetilde{q_{jk}}^2\in\G_3$, $\widetilde{q_{ik}}=1$, or

\viii $q_{kk}=\widetilde{q_{jk}}=q_{jj}=-1$, $q_{ii}=
-\widetilde{q_{ij}}\in\G_3$, $\widetilde{q_{ik}}=1$, or

\viv $q_{jj}=-1$, $\widetilde{q_{ij}}=q_{ii}^{-2}$,
$\widetilde{q_{jk}}=-q_{ii}^{-3}$, $\widetilde{q_{ik}}=1$, or

\vv $q_{ii}=q_{jj}=q_{kk}=-1$,
$\pm\widetilde{q_{ij}}=\widetilde{q_{jk}}\in\G_3$,
 $\widetilde{q_{ik}}=1$,

\noindent Then, for any $p$,
$T_p\left(\left[\left[E_{ij},E_{ijk}\right]_c,E_j\right]_c\right)\in\cJ(\schi)$.
\end{lemma}
\pf Denote $\xb=\left[\left[E_{ij},E_{ijk}\right]_c,E_j\right]_c$;
we analyze each case.

\noindent\vi We begin with the case $p=k$; by \eqref{eqn:Tp adEi Ep} and as $\Eb_i^2$, $\Eb_j^2$, $(\ad_{\cb}\Eb_i)\Eb_p$ are generators of $\cJ(\schi)$ (note that $\schi$ is twist equivalent to $\chi$), we have that
\begin{align*}
T_p(\xb)\cong &
a_1\left[\left[\Eb_{ipj},\Eb_{ij}\right]_{\cb},\Eb_{pj}
\right]_{\cb} \cong
a_1\qb_{ip}\left[\left[\Eb_{pij},\Eb_{ij}\right]_{\cb},\Eb_{pj}
\right]_{\cb}
\\ \cong & a \left[\left[\Eb_{pji},\Eb_{ji}\right]_{\cb},\Eb_{pj} \right]_{\cb}\cong a \left[\Eb_p\Eb_j\Eb_i\Eb_j\Eb_i\Eb_p\Eb_j \right]_{\cb} \quad \left( \mbox{mod }\cJ(\schi) \right),
\end{align*}
for some $a\in\ku^\times$, where we use the order on the letters $p<j<i$. As \eqref{eqn:relacion super C3} is also a generator of $\cJ(\schi)$,
we can write $\Eb_j\Eb_i\Eb_j\Eb_i\Eb_p\Eb_j$ as a linear combination of other words, greater than this word or beginning with $\Eb_p$. Multiplying on the left by  $\Eb_p$ and using the quantum Serre relations $\Eb_p\Eb_j\Eb_i\Eb_j\Eb_i\Eb_p\Eb_j$ is expressed as a linear combination of greater words modulo $\cJ(\schi)$, so by Lemma \ref{lema:palabra como combinacion de mayores en dos ideales}, $T_p(\xb)\in\cJ(\schi)$.

Let $p=j$; note that $m_{pi}=m_{pk}=1$. Also, $\qb_{ii}^{-1}=\qb_{ik}\qb_{ki}$, so $(\ad_{\cb}\Eb_i)^2\Eb_k\in\cJ(\schi)$; use \eqref{eqn:Tp adEi Ep} and work as in the case $p=i$ of Lemma \ref{lemma: relacion estandar B2} to obtain that
\begin{align*}
T_p(\xb)=&a_1^2\left[\Eb_i,\Eb_{ipk}\right]_{\cb}\Fb_p\Lb_p^{-1}-a_1^2q_{ip}^2q_{kp}\Fb_p\Lb_p^{-1}
\left[\Eb_i,\Eb_{ipk}\right]_{\cb}= b\ (\ad_{\cb}\Eb_i)^2\Eb_k \in
\cJ(\schi),
\end{align*}
for some $b\in\ku^\times$.

Let now $p=i$. As in the previous Lemmata, it is enough to prove the statement for
$$ \xb':= \left[\left[ E_{kjp},E_{jp} \right]_c,E_j\right]_c. $$
We apply \eqref{eqn:Tp adEi Ep} to obtain, for the order on the letters $k<i<j$,
$$T_p(\xb')=\left[\left[\left[\Eb_k,a_1\Eb_j\right]_{\cb},\Eb_j\right]_{\cb},(\ad_{\cb}\Eb_p)\Eb_j\right]_{\cb}=a_1^2\left[\Eb_k\Eb_j^2\Eb_p\Eb_j\right]_{\cb}.$$
As $\qb_{jj}\qb_{ji}\qb_{ij}=\qb_{jj}^2\qb_{jk}\qb_{kj}=1$, we deduce by \eqref{eqn:relacion super C3 raiz de orden 3} if $\qb_{jj}\in\G_3$, or by
$(\ad_{\cb} \Eb_j)^2\Eb_p=(\ad_{\cb} \Eb_j)^3\Eb_k=0$, if $\qb_{jj} \notin \G_3$, that $\Eb_k\Eb_j^2\Eb_p\Eb_j$ is a linear combination of greater words, so $T_p(\xb') \in \cJ(\schi)$, and then $T_p(\xb) \in \cJ(\schi)$.

If $p\neq i,j,k$, then $p$ is not connected with these three
vertices, or $p$ is connected just with $i$, or $p$ is connected
with two vertices. For the second case we have that
$q_{pp}\widetilde{q_{pi}}=1$, or $q_{pp}=-1$,
$\widetilde{q_{ip}}\widetilde{q_{ij}}=1$, so
\begin{align*}
T_p(\xb)\cong & \left[\left[ \Eb_{pij},\Eb_{pijk}
\right]_{\cb},\Eb_j\right]_{\cb}\cong (\ad_{\cb}
\Eb_p)\left(\left[\left[\Eb_{pij},\Eb_{ijk}\right]_{\cb},\Eb_j\right]_{\cb}\right)
\\ \cong & (\ad_{\cb} \Eb_p) \left( \left[\Eb_p\Eb_i\Eb_j\Eb_i\Eb_j\Eb_k\Eb_j\right]_{\cb} \right) \quad \mbox{mod }\cJ(\schi) \end{align*}
by using first Lemma \ref{lema:corchete con partes triviales}, then $(\ad_{\cb}\Eb_p)\Eb_j$, $(\ad_{\cb}\Eb_p)\Eb_j$, $(\ad_{\cb}\Eb_p)\Eb_j\in\cJ(\schi)$, and fixing the order $p<i<j<k$.
We conclude that $T_p(\xb)\in\cJ(\schi)$ by using \eqref{eqn:relacion super C4} if $q_{pp}=-1$, or using the quantum Serre relations corresponding to $\ad_{\cb}\Eb_p$ to write $\Eb_p\Eb_i\Eb_j\Eb_i\Eb_j\Eb_k\Eb_j$ as a linear combination of greater words and apply Lemma \ref{lema:palabra como combinacion de mayores en dos ideales} to deduce that
$$\left[\Eb_p\Eb_i\Eb_j\Eb_i\Eb_j\Eb_k\Eb_j\right]_{\cb}\in\cJ(\schi).$$
For the last case, we have two possibilities:
\begin{itemize}
  \item $p$ is connected with $i$ and $j$, in which case the proof
  follows by the fact that \eqref{eqn:relacion super F4-1} is a
  generator of the ideal, or
  \item $p$ is connected with $j$ and $k$, in which case it follows
  because \eqref{eqn:relacion super C4 modificada} is a generator
  of $\cJ(\schi)$.
\end{itemize}
\smallskip

\noindent\vii, \viii, \viv, \vv \, If $p\in\{i,j,k\}$ the proof is
completely analogous to the previous case.

Let $p\neq i,j,k$, so $p$ is not connected with any of these
vertices, or it connected only with $i$, or only with $k$. The first
case is easy. For the second case, $m_{pi}=1$, because $q_{pp}=1$ or
$q_{pp}^{-1}=\widetilde{q_{ip}}\neq-1$, and the solution follows as
in the previous case. For the last case, $m_{pk}=1$ and the proof is
also analogous, considering the previous $\xb'$. \epf
\medskip

\begin{lemma}\label{lema:relaciones super C4, G3}
\noindent \vi Let $i,j,k,l\in\{1,\ldots,\theta\}$ be such that
$q_{kk}=-1$, $\widetilde{q_{jk}}^2= \widetilde{q_{lk}}^{-1}=
q_{ll}$, $q_{jj}\widetilde{q_{ij}}=
q_{jj}\widetilde{q_{jk}}=\widetilde{q_{ik}}=\widetilde{q_{il}}=\widetilde{q_{jl}}=1$.
Then, for all $p$, $ T_p
\left(\left[\left[\left[E_{ijkl},E_k\right]_c,E_j\right]_c,E_k\right]_c
\right) \in \cJ(\schi)$.

\noindent \vii Let $i,j,k\in\{1,\ldots,\theta\}$ be such that
$q_{ii}=q_{jj}=-1$,
$\widetilde{q_{ij}}^3=\widetilde{q_{jk}}^{-1}=q_{kk}\neq\pm1$,
$\widetilde{q_{ik}}=1$. Then, for all $p$,
$T_p\left(\left[\left[E_{ij},\left[E_{ij},E_{ijk}\right]_c\right]_c,E_j
\right]_c\right)\in\cJ(\schi)$.
\end{lemma}
\pf \vi The proof is analogous to \vi of the previous Lemma, because
if $p\neq i,j,k,l$ is connected with some of them, then $p$ is
connected only with $i$ with the same conditions.
\medskip

\noindent \vii If $p\in\{i,j,k\}$ the proof is completely analogous
to the previous Lemma. If $p\neq i,j,k$ is connected with some of
them, then $p$ is connected only with $i$ and $q_{pp}=-1$,
$\widetilde{q_{pi}}=-\widetilde{q_{ij}}\in\G_4$. Anyway, the proof
is analogous to the previous Lemma. \epf
\medskip

\begin{lemma}\label{lemma:casos F4 y C4 modificado}
\vi Let $i,j,k,l\in\{1,\ldots,\theta \}$ be such that
$q_{ll}=\widetilde{q_{lk}}^{-1}=
q_{kk}=\widetilde{q_{jk}}^{-1}=q^2$, $\widetilde{q_{ij}}=
q_{ii}^{-1}=q^3$ for some $q\in \ku^\times$, $q_{jj}=-1$,
$\widetilde{q_{ik}}=\widetilde{q_{il}}=\widetilde{q_{jl}}=1$. Then,
for all $p$, $T_p \left( \left[\left[\left[E_{ijk},E_j\right]_c,
\left[E_{ijkl},E_j\right]_c \right]_c, E_{jk} \right]_c
\right)\in\cJ(\schi)$.
\smallskip

\noindent \vii Let $i,j,k,l\in\{1,\ldots,\theta \}$ be such that
$\widetilde{q_{jk}}= \widetilde{q_{ij}}= q_{jj}^{-1}\in
\G_4'\cup\G_6'$, $q_{ii}=q_{kk}=-1$,
$\widetilde{q_{ik}}=\widetilde{q_{il}}=\widetilde{q_{jl}}=1$,
$\widetilde{q_{jk}}^3= \widetilde{q_{lk}}$. Then, for all $p$, $T_p
\left( \left[\left[E_{ijk},\left[E_{ijkl}, E_k \right]_c\right]_c,
E_{jk} \right]_c \right)\in\cJ(\schi)$.
\end{lemma}
\pf \vi Let $\xb=\left[\left[E_{ijk},\left[E_{ijkl}, E_k
\right]_c\right]_c, E_{jk} \right]_c$. If $p\notin\{i,j,k,l\}$, then
$p$ is not connected with $i,j,k,l$, so we consider the case
$p\in\{i,j,k,l\}$. If $p=i$, the diagram for $\schi$ is the same,
and $T_p(\xb)$ corresponds to a relation of degree
$3\alpha_p+5\alpha_j+3\alpha_k+\alpha_l$, obtained after to apply
$\ad_c\Eb_p$ to \eqref{eqn:relacion super F4-1}. Similar situations
hold when $p=k$ and $p=l$. If $p=j$, then $T_p(\xb)$ corresponds to
the relation \eqref{eqn:relacion super C3}, and the proof follows.

\medskip
\noindent \vii The proof is very similar to previous item, because
$p\in\{i,j,k,l\}$, or $p$ is not connected with $i,j,k,l$. \epf

\medskip

\begin{lemma}\label{lemma: varios casos rango 4}
Let $i,j,k,l\in\{1,\ldots,\theta \}$ be such that one of the
following situations hold:
\begin{itemize}
  \item[$\circ$] $q_{kk}=-1$, $q_{ii}=\widetilde{q_{ij}}^{-1}= q_{jj}^2$,
$\widetilde{q_{kl}}= q_{ll}^{-1}= q_{jj}^3$, $\widetilde{q_{jk}}=
q_{jj}^{-1}$,
$\widetilde{q_{ik}}=\widetilde{q_{il}}=\widetilde{q_{jl}}=1$,
  \item[$\circ$] $q_{ii}=\widetilde{q_{ij}}^{-1}= -q_{ll}^{-1}=-\widetilde{q_{kl}}$,
$q_{jj}=\widetilde{q_{jk}}=q_{kk}=-1$,
$\widetilde{q_{ik}}=\widetilde{q_{il}}=\widetilde{q_{jl}}=1$, or
  \item[$\circ$] $q_{jj}=\widetilde{q_{jk}}^{-1}\in\G_3$, $q_{ii}=\widetilde{q_{ij}}^{-1}=q_{ll}=\widetilde{q_{kl}}^{-1}=-q_{jj}$,
$q_{kk}=-1$,
$\widetilde{q_{ik}}=\widetilde{q_{il}}=\widetilde{q_{jl}}=1$.
\end{itemize}
Then, $T_p \left( \left[\left[E_{ijkl}, E_j \right]_c, E_k
\right]_c- q_{jk}(\widetilde{q_{ij}}^{-1}-q_{jj})
\left[\left[E_{ijkl},E_k\right]_c, E_j \right]_c
\right)\in\cJ(\schi)$, for all $p$.
\end{lemma}
\pf No one of these diagrams can be extended in order to have a
connected diagram with finite root system. In consequence, it is
enough to consider (as in the previous Lemmata) the cases
$p\in\{i,j,k,l\}$. The proof for these cases is similar to the Lemma
\ref{lema:relacion triangulo}, up to consider the necessary
relations under the conditions of these new situations. \epf

\medskip

\begin{lemma}\label{lemma: relaciones fila 18, rango3}
\vi Let $i,j,k\in\{1,\ldots,\theta\}$ be such that
$q_{kk}=q_{jj}=\widetilde{q_{ij}}^{-1}=\widetilde{q_{jk}}^{-1}\in
\G_9$, $\widetilde{q_{ik}}=1$, $q_{ii}=q_{kk}^6$. Then, for all $p$,
$T_p \left( \left[ \left[ E_{iij} , E_{iijk} \right]_c, E_{ij}
\right]_c \right)\in\cJ(\schi)$.
\smallskip

\noindent \vii Let $i,j,k\in\{1,\ldots,\theta\}$ be such that
$q_{ii}=\widetilde{q_{ij}}^{-1}\in \G_9$,
$q_{jj}=\widetilde{q_{jk}}^{-1}=q_{ii}^5$, $\widetilde{q_{ik}}=1$,
$q_{kk}=q_{ii}^6$. Then, for all $p$, $T_p \left( [\left[E_{ijk},
E_{j} \right]_c, E_k]_c -(1 + \widetilde{q_{jk}})^{-1}q_{jk} \left[
\left[E_{ijk}, E_{k} \right]_c , E_{j} \right]_c
\right)\in\cJ(\schi)$.
\end{lemma}
\pf For both cases, there exist no extensions of these diagrams. In
consequence, it is enough to consider the case $p\in\{i,j,k\}$.
\smallskip

\noindent \vi The cases $p=k$ and $p=j$ follows easily because the
diagram for $\schi$ in these cases coincide with the one for $\chi$.

If $p=i$, then $T_p \left( \left[ \left[ E_{ppj} , E_{ppjk}
\right]_c, E_{pj} \right]_c \right)$ is, up to an scalar,
$[\Eb_{pjk}, \Eb_{j}]_{\cb}$, which belongs to $\cJ(\schi)$ because
$(\ad_{\cb}\Eb_j)^2\Eb_p, (\ad_{\cb}\Eb_j)^2\Eb_k \in \cJ(\schi)$
(we have that $\qb_{jj}\widetilde{\qb_{jk}}=\qb_{jj}
\widetilde{\qb_{jp}}=1$).
\smallskip

\noindent \vii The proof is similar to the one for Lemma
\ref{lema:relacion triangulo}.

\epf

\medskip

\begin{lemma}\label{lema: relacion C3 con raiz en G3, G3 raiz en G4}
\vi Let $i,j,k\in\{1,\ldots,\theta\}$ be such that
$q_{jj}=\widetilde{q_{ij}}^2=\widetilde{q_{jk}}\in \G_3$,
$\widetilde{q_{ik}}=1$. Then, for all $p$, $T_p \left( \left[ \left[
E_{ijk} , E_j \right]_c, E_j \right]_c \right) \in \cJ(\schi)$.

\noindent \vii Let $i,j,k\in\{1,\ldots,\theta\}$ be such that
$q_{jj}=\widetilde{q_{ij}}^3=\widetilde{q_{jk}}\in\G_4$,
$\widetilde{q_{ik}}=1$. Then, for all $p$,
$T_p\left(\left[\left[\left[E_{ijk},E_j\right]_c,E_j\right]_c,E_j\right]_c\right)\in
\cJ(\schi)$.
\end{lemma}
\pf \vi Let $\xb=\left[ \left[ E_{ijk} , E_j \right]_c, E_j
\right]_c$. For the case $p=k$, note that $m_{pj}=1$ in all the
cases, so for the order $i<p<j$ on the letters we have by
\eqref{eqn:Tp adEi Ep}:
$$ T_p(\xb)=a_1\left[\left[\Eb_{ij},\Eb_{pj}\right]_{\cb},\Eb_{pj}\right]_{\cb}
=\left[\Eb_i\Eb_j\Eb_p\Eb_j\Eb_p\Eb_j\right]_{\cb}. $$
As $(\ad_{\cb}\Eb_p)^2\Eb_j\in\cJ(\schi)$, or \eqref{eqn:relacion super C3} is a generator of $\cJ(\schi)$, and also $\Eb_j^3,(\ad_{\cb}\Eb_k)\Eb_i\in\cJ(\schi)$, $\Eb_i\Eb_j\Eb_k\Eb_j\Eb_k\Eb_j$ can be written as a linear combination of greater words modulo $\cJ(\schi)$, so $T_p(\xb)\in\cJ(\schi)$.

Consider now $p=j$. By \eqref{eqn:Tp adEi Ep} and the relations defining $\cJ(\schi)$,
$$ T_p(\xb)=\left[\left[\left[\Eb_i,(\ad_{\cb}\Eb_p)^2\Eb_k\right]_{\cb},\Fb_p\Lb_p^{-1}\right]_{\cb},\Fb_p\Lb_p^{-1}\right]_{\cb}=a(\ad_{\cb}\Eb_i)\Eb_k,$$
for some $a\in\ku^\times$, so $T_p(\xb)\in\cJ(\schi)$.

Let $p=i$. It is equivalent to prove that $T_p(\xb')\in\cJ(\schi)$, where
$\xb'= \left[ \left[E_{kjp}, E_j \right]_c , E_j \right]_c$.
We note that $m_{ij}=1$ for all the possible diagrams, so $ T_p(\xb)=a_1\left[\left[\Eb_{kj},\Eb_{pj}\right]_{\cb},\Eb_{pj}\right]_{\cb}$,
and a proof similar to the case $p=k$ tells us that $T_p(\xb')\in\cJ(\schi)$, so $T_p(\xb)\in\cJ(\schi)$.

Finally, if $p\neq i,j,k$, then $p$ is not connected to any of these vertices, or $p$ is connected only with $i$, or it is connected only with $k$. The proof of the first case is again trivial, and for the other two cases $T_p(\xb)\in\cJ(\schi)$, using Lemma \ref{lema:corchete con partes triviales} and the fact that $\schi$ is twist equivalent to $\chi$.
\medskip

\vii The proof is analogous to \vi.
\epf
\medskip

\begin{lemma}\label{lema:relacion similar a super C3}
\vi Let $i,j,k\in\{1,\ldots,\theta\}$ be such that
$q_{ii}=\widetilde{q_{ij}}=-1$, $q_{jj}=\widetilde{q_{jk}}^{-1}\neq
-1$, $\widetilde{q_{ik}}=1$. Then, for all $p$, $T_p \left(
\left[E_{ij}, E_{ijk} \right]_c \right) \in \cJ(\schi)$.
\smallskip

\vii Let $i,j,k\in\{1,\ldots,\theta\}$ be such that
$q_{jj}=q_{kk}=\widetilde{q_{jk}}=-1$,
$q_{ii}=-\widetilde{q_{ij}}\in \G_3$, $\widetilde{q_{ik}}=1$. Then,
for all $p$, $T_p \left( \left[E_{iijk}, E_{ijk} \right]_c \right)
\in \cJ(\schi)$.
\end{lemma}

\pf \vi Let $\xb=\left[E_{ij}, E_{ijk}\right]_c$. If $p=k$, $m_{pj}=1$ in all the possible diagrams and then:
$$ T_p(\xb)=a_1[\Eb_{ikj},\Eb_{ij}]_{\cb}=a_1\qb_{ip}[\Eb_{kij},\Eb_{ij}]_{\cb}. $$
We consider the two possible values of $q_{pp}$. When $q_{pp}\neq-1$, the diagram is of Cartan $C_3$ type, associated to a root of order $4$, and $\schi$ is twist equivalent to $\chi$. Therefore $T_p(\xb)\in\cJ(\schi)$, using that \eqref{eqn:relacion parecida a super C3} is again a generator of $\cJ(\schi)$. If $q_{pp}=-1$, then $q_{jj}\in\G_3\cup\G_4\cup\G_6$, and $\qb_{jj}=\widetilde{\qb_{ij}}=\qb_{ii}=-1$, so \eqref{eqn:relacion dos vertices con -1} is a generator of $\cJ(\schi)$ and
$$ T_p(\xb)=a_1\qb_{ip}[\Eb_{kij},\Eb_{ij}]_{\cb}=a_1\qb_{ip} \ad_{\cb}\Eb_k\left(E_{ij}^2\right)\in\cJ(\schi).$$

If $p=j$, we consider the different possible orders of $q_{pp}$. If $q_{pp}\in\G_4\cup\G_6$, then $p$ is a Cartan vertex and $\schi$ is twist equivalent to $\schi$, and $m_{pi}=2,3$
$$ T_p(\xb)=a_{m_{pi}} \left[ (\ad_{\cb}\Eb_p)^{m_{pi}-1}\Eb_i,[(\ad_{\cb}\Eb_p)^{m_{pi}-1}\Eb_i,\Eb_{pk}]_{\cb} \right]_{\cb}= [\Eb_p^{m_{pi-1}}\Eb_i\Eb_p^{m_{pi-1}}\Eb_i\Eb_p\Eb_j]_{\cb}. $$
We use the quantum Serre relations and \eqref{eqn:relacion parecida a super C3} to write $\Eb_p^{m_{pi-1}}\Eb_i\Eb_p^{m_{pi-1}}\Eb_i\Eb_p\Eb_j$ as a linear combination of greater words modulo $\cJ(\schi)$. If $q_{pp}\in\G_3$, then
$$ \xymatrix{ \circ^{q_{ii}=-1} \ar@{-}[r]_{-1} &  \circ^{q_{pp}} \ar@{-}[r]_{q_{pp}^2}  & \circ^{q_{kk}=-1} } \quad  \leftrightsquigarrow_{s_p}  \quad
\xymatrix{ & \circ^{q_{pp}} \ar@{-}[rd]^{q_{pp}^2} & \\ \circ^{-q_{pp}} \ar@{-}[ru]^{-q_{pp}^2} \ar@{-}[rr]_{-q_{pp}^2} &  & \circ^{q_{kk}=-1}. } $$
The result follows in a similar way, but using that \eqref{eqn:relacion estandar B2} and \eqref{eqn:relacion triangulo} are generators of $\cJ(\schi)$ in this case.

If $p=i$, by Lemma \ref{lema:dos conjuntos de generadores}, it is equivalent to prove that
$$ T_p(\xb')\in \cJ(\schi), \qquad\xb':= \left[E_{kjp},E_{jp}\right]_c.$$
Note that $(\ad_{\cb}\Eb_j)^2\Eb_k\in\cJ(\schi)$, because for all the possible diagrams $\qb_{jj}^{-1}=\qb_{jk}\qb_{kj}\neq-1$. By \eqref{eqn:Tp adEi Ep}, we have
$$ T_p(\xb')=a_1^2 \ \left[ (\ad_{\cb}\Eb_k)\Eb_j,\Eb_j \right]_{\cb}\in\cJ(\schi). $$

Finally, if $p\neq i,j,k$, then $p$ is not connected with any of these vertices, or $p$ is connected just with $i$ and $m_{pi}=1$, or $p$ is connected just with $k$ and $m_{pk}=1$. The proof is analogous to the corresponding case in previous Lemmata.
\medskip

\vii If $p\neq i,j,k$, then $p$ is not connected with them and the result follows easily. In consequence, we consider the case $p\in\{i,j,k\}$. If $p=k$,
$\schi$ is twist equivalent to $\chi$ and \eqref{eqn:relacion estandar B2} is a generator of $\cJ(\schi)$. Therefore \eqref{eqn:Tp adEi Ep} implies that
$$ T_p \left( \left[E_{iijk}, E_{ijk} \right]_c \right) = a_2^2\left[ \Eb_{iij},\Eb_{ij}\right]_{\cb} \in \cJ(\schi).$$
If $p=j$, we have that
$$ \xymatrix{ \circ^{q_{ii}} \ar@{-}[r]_{-q_{ii}} &  \circ^{-1} \ar@{-}[r]_{-1}  & \circ^{-1} } \quad  \leftrightsquigarrow_{s_p}  \quad
\xymatrix{ & \circ^{-1} \ar@{-}[rd]^{-1} & \\ \circ^{q_{ii}^2} \ar@{-}[ru]^{-q_{ii}^2} \ar@{-}[rr]_{q_{ii}} &  & \circ^{-1}. } $$
By \eqref{eqn:Tp adEi Ep} and fixing the order $p<i<k$, we have that
$$ T_p \left( \left[E_{iijk}, E_{ijk} \right]_c \right) = a_2^2\left[ \left[ \Eb_{pi},\Eb_{pik}\right]_{\cb},\Eb_{pik}\right]_{\cb} =[\Eb_p\Eb_i\Eb_p\Eb_i\Eb_k\Eb_p\Eb_i\Eb_k]_{\cb} $$
We write $\Eb_p\Eb_i\Eb_p\Eb_i\Eb_k\Eb_p\Eb_i\Eb_k$ as a linear combination of greater words modulo $\cJ(\schi)$ using that \eqref{eqn:relacion triangulo}, \eqref{eqn:relacion estandar B2}, $\Eb_p^2$, $\Eb_k^2$ and $\Eb_{iik}$ are generators of $\cJ(\schi)$, so Lemma \ref{lema:palabra como combinacion de mayores en dos ideales} implies that $T_p \left( \left[E_{iijk}, E_{ijk} \right]_c \right)\in \cJ(\schi)$.

Finally let $p=i$. By Lemma \ref{lema:dos conjuntos de generadores}, it is equivalent to prove that
$$ T_p(\xb')\in \cJ(\schi), \qquad\xb':= \left[E_{kjp},E_{kjpp}\right]_c.$$
Note that $\schi$ is twist equivalent to $\chi$, so $x_{jp}^2$ is a generator of $\cJ(\schi)$. Applying \eqref{eqn:Tp adEi Ep},
$$ T_p(\xb')= a_2^2a_1\left[\Eb_{kjp},\Eb_{kj}\right]_c = a_2^2a_1\left[\Eb_{kj}^2,\Eb_p \right]_c \in\cJ(\schi). $$
Therefore $T_p \left( \left[E_{iijk}, E_{ijk} \right]_c \right)\in \cJ(\schi)$.
\epf
\medskip

\begin{lemma}
\vi Let $i,j,k\in\{1,\ldots,\theta\}$ be such that $q_{ii}=\widetilde{q_{ij}}=-\widetilde{q_{ik}}\in\G_3$, $\widetilde{q_{jk}}=1$, $q_{jj}=-1$,
$q_{kk}\in\{-1,\widetilde{q_{ik}}^{-1}\}$. Then, for all $p$,
$$ T_p \left( \left[E_i, \left[ E_{ij},E_{ik} \right]_c \right]_c+q_{jk}q_{ik}q_{ji} \left[ E_{iik} ,E_{ij} \right]_c
+q_{ij}E_{ij}E_{iik} \right) \in\cJ(\schi).$$
\smallskip

\vii Let $i,j,k\in\{1,\ldots,\theta\}$ be such that $q_{ii}=q_{kk}=-1$, $\widetilde{q_{ik}}=1$, $\widetilde{q_{ij}}\in\G_3$,
$q_{jj}=-\widetilde{q_{jk}}= \pm\widetilde{q_{ij}}$. Then, for all $p$,
$$ T_p \left( \left[E_i, E_{jjk}\right]_c-q_{kj}^{-1}(1+q_{jj}^2) \left[ E_{ijk} ,E_{j} \right]_c+(1+q_{jj})(1+q_{jj}^2)E_{j}E_{ijk} \right) \in\cJ(\schi).$$
\end{lemma}
\pf
If $p\neq i,j,k$, then $p$ is not connected with these three vertices and the proof follows easily for both relations. We consider then $p\in\{i,j,k\}$ for each item.
\smallskip

\vi If $p=k$, $\Eb_p\Eb_i\Eb_p\Eb_i\Eb_j\Eb_i$ is a linear combination of greater words by \eqref{eqn:relacion estandar B2} and the quantum Serre relations, depending on the value of $q_{kk}$. In this way, there exist $a,b\in\ku$ such that
$$ \xb':=[\Eb_p\Eb_i\Eb_p\Eb_i\Eb_j\Eb_i]_{\cb}+a[\Eb_p\Eb_i^2\Eb_p\Eb_i\Eb_j]_{\cb}+b[\Eb_p\Eb_i\Eb_j]_{\cb}[\Eb_p\Eb_i^2]_{\cb}\in\cJ(\schi). $$
On the other hand, by \eqref{eqn:Tp adEi Ep} and for the order on the letters $p<i<j$,
\begin{align*}
T_p\left(\left[E_i, \left[ (\ad_c E_i)E_j,(\ad_c E_i)E_p \right]_c \right]_c\right) &= a_1 \left[ (\ad_{\cb} \Eb_p) \Eb_i,
\left[ (\ad_{\cb} \Eb_p)(\ad_{\cb} \Eb_i) \Eb_j,\Eb_i \right]_{\cb} \right]_{\cb}
\\&= a_1[\Eb_p\Eb_i\Eb_p\Eb_i\Eb_j\Eb_i]_{\cb},
\\ T_p\left(\left[ (\ad_c E_i)^2E_p ,(\ad_c E_i)E_j \right]_c\right) &= a_1\left[ \left[ (\ad_{\cb} \Eb_p) \Eb_i, \Eb_i \right]_{\cb},
(\ad_{\cb}\Eb_p)(\ad_{\cb}\Eb_i)\Eb_j \right]_{\cb}
\\&= a_1[\Eb_p\Eb_i^2\Eb_p\Eb_i\Eb_j]_{\cb};
\\ T_p\left((\ad_c E_i)E_j (\ad_c E_i)^2 E_p\right) &= a_1 (\ad_{\cb}\Eb_p)(\ad_{\cb}\Eb_i)\Eb_j\left[ (\ad_{\cb} \Eb_p) \Eb_i, \Eb_i \right]_{\cb}
\\&=a_1 [\Eb_p\Eb_i\Eb_j]_{\cb}[\Eb_p\Eb_i^2]_{\cb}.
\end{align*}
Calculating explicitly the scalars $a,b$, we notice that $T_p\left(\xb\right)=a_1\xb'\in\cJ(\schi)$. The case $p=j$ is analogous.

Finally the case $p=i$ follows as the corresponding case in Lemma \ref{lema:relacion triangulo}.
\medskip

\vii The proof is analogous to the previous case.
\epf
\medskip

\begin{lemma}\label{lema:relaciones iguales a cuadrados de raices}
\vi Let $i,j\in\{1,\ldots,\theta\}$ be such that $m_{ij}, m_{ji} >1$. Then, for all $p$,
$$T_p \left((1-\widetilde{q_{ij}})q_{jj}q_{ji}\left[E_i,\left[E_{ij},E_j\right]_c\right]_c-(1+q_{jj})(1-q_{jj}\widetilde{q_{ij}})E_{ij}^2\right)\in \cJ(\schi).$$
\smallskip

\noindent \vii Let $i,j\in\{1,\ldots,\theta\}$ be such that $q_{jj}=-1$, $q_{ii}\widetilde{q_{ij}}\notin\G_6$ or $m_{jj}=2$, and $q_{ii}\in \G_4$, $m_{ij}=4$, or
$m_{ij}\in \{4,5\}$. Then, for all $p$,
$$T_p\left( \left[E_i,\left[E_{iij},E_{ij}\right]_c\right]_c-\frac{1-q_{ii}\widetilde{q_{ij}}-q_{ii}^2\widetilde{q_{ij}}^2q_{jj}}
{(1-q_{ii}\widetilde{q_{ij}})q_{ji}} E_{iij}^2\right)\in \cJ(\schi).$$
\smallskip

\noindent \viii Let $i,j\in\{1,\ldots,\theta\}$ be such that $q_{jj}=-1$, $5\alpha_i+4\alpha_j\in\Delta_+^\chi$. Let
\begin{align*}
& \upsilon=\widetilde{q_{ij}}, \qquad a=(1-\upsilon)(1-q_{ii}^4\upsilon^3)-(1-q_{ii}\upsilon)(1+q_{ii})q_{ii}\upsilon
\\ & b=(1-\upsilon)(1-q_{ii}^6\upsilon^5)-a\ q_{ii}\upsilon,
\\ & d= \frac{b-(1+q_{ii})(1-q_{ii}\upsilon)(1+\upsilon+q_{ii}\upsilon^2)q_{ii}^6\upsilon^4}{a\ q_{ii}^3q_{ij}^2q_{ji}^3}.
\end{align*}
Then, for all $p$, $T_p\left([E_{2\alpha_i+\alpha_j},E_{4\alpha_i+3\alpha_j}]_c-d\ E_{3\alpha_i+2\alpha_j}^2\right)\in\cJ(\schi)$.
\end{lemma}
\pf
\vi Let $\xb$ be the relation we are considering here. We note that if $p\neq i,j$ then $m_{pi}=m_{pj}=0$, so the proof follows easily. Moreover the conditions about $i,j$ are the same but one relation implies the other holds in $\cU(\chi)$ too by Lemma \ref{lema:palabra como combinacion de mayores en dos ideales}. Therefore it is enough to consider one of cases $p=i$ or $p=j$; consider $p=j$, in order to apply \eqref{eqn:Tp adEi Ep}. Note that $m_{pi}=2,3$.

If $m_{pi}=3$, then $m_{ip}=2$ and we have that
$$ T_p \left( \left[E_i, \left[ E_{ip}, E_p \right]_c \right]_c \right) = \left[ \Eb_{pppi}, \Eb_{pi} \right]_{\cb}.$$
By \eqref{eq:identidad jacobi} we can write $T_p(\xb)$ as a linear combination of
$$ \left[\Eb_p,\left[\Eb_{ppi},\Eb_{pi}\right]_{\cb}\right]_{\cb}, \qquad \Eb_{ppi}^2;$$
note that $\left[\Eb_p,\left[\Eb_{ppi},\Eb_{pi}\right]_{\cb}\right]_{\cb}=[\Eb_p^3\Eb_i\Eb_p\Eb_i]_{\cb}$ if we consider
$p<i$. Using the quantum Serre relations or \eqref{eqn:relacion mij mayor que dos, raiz alta}, depending on the case (there exist two possible diagrams), $\Eb_p^3\Eb_i\Eb_p\Eb_i$ is expressed as a linear combination of greater words modulo $\cJ(\schi)$, so $T_p(\xb)\in\cJ(\schi)$ by Lemma \ref{lema:palabra como combinacion de mayores en dos ideales}.

If $m_{pi}=2$, there exist three diagrams such that $m_{ip}=2$, and two such that $m_{ip}=3$. In all the cases,
\begin{align*}
T_p\left(E_{ip}^2\right)=&a_2^2 \ \Eb_{pi}^2,
\\ T_p \left(\left[E_i,\left[E_{ip},E_p\right]_c\right]_c\right)=&a_2a_1\, \left[\Eb_{ppi}, \Eb_i\right]_{\cb} =a_2a_1\, \left[\Eb_p,\left[\Eb_{pi},\Eb_i\right]_{\cb}\right]_{\cb}+a_2a_1\ \qb_{pi}(\qb_{pp}-\qb_{ii})\Eb_{pi}^2,
\end{align*}
where we use \eqref{eq:identidad jacobi} for the last equality. If $q_{pp}=-\zeta$, $\widetilde{q_{ip}}=\zeta^7$, $q_{ii}=\zeta^3$, for some primitive root $\zeta\in\G_9$, then $\schi$ is twist equivalent to $\chi$ and \eqref{eqn:relacion mij,mji mayor que 1} is a generator of $\cJ(\schi)$, so $T_p(\xb)\in\cJ(\schi)$ by this relation and Lemma \ref{lema:palabra como combinacion de mayores en dos ideales}. For the other braidings $\qb_{ii}=-1$, so $\left[\Eb_{pi},\Eb_i\right]_{\cb}\in\cJ(\schi)$ and the coefficient of $\Eb_{pi}^2$ in the expression of $T_p(\xb)$ is zero. Then $T_p(\xb)\in\cJ(\schi)$.

\medskip
\vii Let $\xb$ be the relation we are considering in this item. First we consider $p=j$; if $q_{pp}=-1$, then $m_{pi}=1$, so
$$s_p(3\alpha_i+\alpha_p)=3\alpha_i+2\alpha_p, \ s_p(3\alpha_i+2\alpha_p)=3\alpha_i+\alpha_p \in \Delta_+^{\schi},$$
so $\underline{m}_{ip}\geq3$. Applying \eqref{eqn:Tp adEi Ep} we have that:
\begin{align*}
T_p\left(E_{iip}^2\right)&=a_1^2 \ \left[\Eb_{pi}, \Eb_i\right]_{\cb}^2,
\\ T_p \left(\left[E_i,\left[E_{iij},E_{ij}\right]_c\right]_c\right)&=a_1^2\ \left[\Eb_{pi},\left[\left[\Eb_{pi}, \Eb_i\right]_{\cb},\Eb_i\right]_{\cb}\right]_{\cb}.
\end{align*}
As $\underline{m}_{ip}\geq3$, \eqref{eqn:relacion mij mayor que dos, raiz alta} is a generator of $\cJ(\schi)$, or $\underline{m}_{ip}=3$, $\qb_{ii}\notin\G_4$, so $\Eb_p\Eb_i\Eb_p\Eb_i^3$ can be written as a linear combination of greater words modulo $\cJ(\schi)$, for the order $p<i$, using the corresponding quantum Serre relation and $\Eb_p^2$. In both cases we apply Lemma \ref{lema:palabra como combinacion de mayores en dos ideales} to deduce that $T_p(\yb)\in\cJ(\schi)$.

If $m_{pi}=2$, then $m_{ip}=3$; in this case,
\begin{align*}
T_p\left(E_{iip}^2\right)&=a_2^2 \ \left[\Eb_{ppi},\Eb_{pi}\right]_{\cb}^2,
\\ T_p \left(\left[E_i,\left[E_{iip},E_{ip}\right]_c\right]_c\right)&=a_2^2\ \left[\Eb_{ppi},\Eb_{4\alpha_p+3\alpha_i} \right]_{\cb}.
\end{align*}
We have two possibilities for $\schi$:
\begin{itemize}
  \item $\xymatrix{\circ^{\zeta^8} \ar@{-}[r]^{\zeta^5} & \circ^{-1}}$, $\zeta\in\G_{24}$, so \eqref{eqn:relacion potencia alta} is a generator of $\cJ(\schi)$,
  \item $\xymatrix{\circ^{\zeta^5}\ar@{-}[r]^{-\zeta^{13}}&\circ^{-1}}$, $\zeta\in\G_{15}$, so \eqref{eqn:relacion (m+1)alpha i+m alpha j, caso 3} and $\Eb_p^3$ are generators of $\cJ(\schi)$.
\end{itemize}
Then $\Eb_p^2\Eb_i\Eb_p^2\Eb_i\Eb_p\Eb_i\Eb_p\Eb_i$ is written as a linear combination of greater words modulo $\cJ(\schi)$ in both cases, so by Lemma \ref{lema:palabra como combinacion de mayores en dos ideales} we have that $T_p(\yb)\in\cJ(\schi)$.

Let $p=i$; by Lemma \ref{lema:dos conjuntos de generadores}, it is equivalent to prove that $T_p(\yb')\in\cJ(\schi)$, where
$$ \yb':=\left[\left[E_{jp},\left[E_{jp},E_p\right]_c \right]_c,E_p\right]_c-a \ \left(\left[E_{jp},E_p\right]_c\right)^2,$$
and $a\in\ku^\times$ is fixed. Note that
\begin{align*}
T_p & \left(\left[\left[(\ad_cE_j)E_p,\left[(\ad_cE_j)E_p,E_p\right]_c \right]_c,E_p\right]_c\right)
\\ &= a_{m_{pi}}^2a_{m_{pi}-1}\left(\left[ (\ad_{\cb}\Eb_p)^{m_{pi}-1}\Eb_j, (\ad_{\cb}\Eb_p)^{m_{pi}-2}\Eb_j\right]_{\cb}\Fb_p\Lb_p^{-1} \right.
\\ &\left.\quad -q_{jp}^2q_{pp}^3\Fb_p\Lb_p^{-1}\left[(\ad_{\cb}\Eb_p)^{m_{pi}-1}\Eb_j, (\ad_{\cb}\Eb_p)^{m_{pi}-2}\Eb_j\right]_{\cb} \right)
\\T_p&\left(\left[E_{jp},E_p\right]_c\right)= a_{m_{pi}}a_{m_{pi}-1} (\ad_{\cb}\Eb_p)^{m_{pi}-2}\Eb_i
\end{align*}
In any case, $T_p(\yb')\in\ker \pi_{\schi}$ is a linear combination of
$$[\Eb_p^{m_{ip}-1}\Eb_i\Eb_p^{m_{ip}-3}\Eb_i]_{\cb}, \qquad [\Eb_p^{m_{ip}-2}\Eb_i]_{\cb}^2,$$
so by Lemma \ref{lema:palabra como combinacion de mayores en dos ideales}, $T_p(\yb')\in\cJ(\schi)$, because \eqref{eqn:relacion mij,mji mayor que 1}, (respectively, \eqref{eqn:relacion mij mayor que dos, raiz alta}, \eqref{eqn:relacion potencia alta}) is a generator of $\cJ(\schi)$ if $m_{pi}=3$, (respectively, $m_{pi}=4$, $m_{pi}=5$).

Finally we take $p\neq i,j$, so $p$ is not connected with $i$ and $j$ (and the proof follows easily by Lemma \ref{lema:corchete con partes triviales}), or $p$ is connected only with $i$, and $q_{ii}=\widetilde{q_{ij}}=\widetilde{q_{pi}}^{-1}\in\G_4$, $q_{pp}=-1$. Consider the order $p<i<j$, so
\begin{align*}
T_p\left(E_{iij}\right) &=\left[\Eb_{pi},\Eb_{pij}\right]_{\cb},
\\ T_p \left(\left[E_i,\left[E_{iij},E_{ij}\right]_c\right]_c\right) &=\left[\Eb_{pi},\left[\left[\Eb_{pi}, \Eb_{pij} \right]_{\cb},\Eb_{pij} \right]_{\cb} \right]_{\cb}= [\Eb_p\Eb_i\Eb_p\Eb_i\Eb_p\Eb_i\Eb_j\Eb_p\Eb_i\Eb_j]_{\cb}.
\end{align*}
By \eqref{eqn:relacion super G3 raiz de orden 4}, $\Eb_i\Eb_p\Eb_i\Eb_p\Eb_i\Eb_j\Eb_p\Eb_i$ can be written as a linear combination of other words modulo $\cJ(\schi)$, which are greater than it or they begin with $\Eb_p$; multiplying on the left by $\Eb_p$, on the right by $\Eb_j$, and using that $\Eb_p^2\in\cJ(\schi)$, $\Eb_p\Eb_i\Eb_p\Eb_i\Eb_p\Eb_i\Eb_j\Eb_p\Eb_i\Eb_j$ is a linear combination of greater words modulo $\cJ(\schi)$, so $T_p(\yb)\in\cJ(\schi)$ by a similar argument to the previous steps.
\medskip

\viii The proof is analogous to the previous items, where we note that in the two possible cases $q_{jj}=-1$, and if $p\neq i,j$, then $p$ is not connected with them.
\epf
\medskip

\begin{lemma}
\vi Let $i,j \in \{1, \ldots, \theta \}$ be such that $4\alpha_i+3\alpha_j\notin \Delta_+^\chi$, $q_{jj}=-1$ or $m_{ji}=2$, and also $m_{ij}\geq 3$, or $m_{ij}=2$, $q_{ii}\in\G_3$. Then, $T_p\left([E_{3\alpha_i+2\alpha_j},E_{ij} ]_c\right)\in\cJ(\schi)$, for all $p$.
\smallskip

\noindent \vii Let $i,j\in\{1,\ldots,\theta\}$ be such that $4\alpha_i+3\alpha_j\in\Delta_+^\chi$, $5\alpha_i+4\alpha_j\notin\Delta_+^\chi$. Then, for all $p$, $T_p\left([E_{4\alpha_i+3\alpha_j}, E_{ij}]_c\right)\in\cJ(\schi)$.
\smallskip

\noindent \viii Let $i,j\in\{1,\ldots,\theta\}$ be such that $3\alpha_i+2\alpha_j\in\Delta_+^\chi$, $5\alpha_i+3\alpha_j\notin\Delta_+^\chi$, and
$q_{ii}^3\widetilde{q_{ij}}$, $q_{ii}^4\widetilde{q_{ij}} \neq 1$. Then, $T_p\left([E_{iij},E_{3\alpha_i+2\alpha_j}]_c\right)\in\cJ(\schi)$ for all $p$.
\smallskip

\noindent \viv Let $i,j\in\{1,\ldots,\theta\}$ be such that $5\alpha_i+2\alpha_j\in\Delta_+^\chi$, $7\alpha_i+3\alpha_j\notin\Delta_+^\chi$. Then, for all $p$, $T_p\left([[E_{iiij},E_{iij}],E_{iij}]_c\right)\in\cJ(\schi)$.
\end{lemma}
\pf For these four sets of conditions, if $p\neq i,j$ then $p$ is not connected with $i$ and $j$, so the proof follows easily using Lemma \ref{lema:corchete con partes triviales}, or we have a diagram as in Lemma \ref{lema:relaciones iguales a cuadrados de raices}, \vii, and the proof is analogous to this one. In consequence we will consider $p=i$ and $p=j$ for each one of these cases.

\noindent\vi Let $\xb=[E_{3\alpha_i+2\alpha_j},E_{ij}]_c$, and take $p=j$. If $m_{pi}=1$ (that is, $q_{pp}=-1$ or $q_{pp}\widetilde{q_{ip}}=1$), we have that
$$ s_p(3\alpha_i+2\alpha_p)=3\alpha_i+\alpha_p\in\Delta_+^{\schi}, \quad s_p(4\alpha_i+3\alpha_p)=4\alpha_i+\alpha_p\notin\Delta_+^{\schi}. $$
Therefore $\underline{m}_{ip}=3$, so $\Eb_i^4$ (respectively, $(\ad_{\cb}\Eb_i)^4\Eb_p$) is a generator of $\cJ(\schi)$, if $\qb_{ii}$ belongs (respectively, does not belong) to $\G_4$. By \eqref{eqn:Tp adEi Ep} and the previous relations, depending on the case,
$$T_p(\xb)= a_1^4 \left[ \left[ \left[\Eb_{pi} ,\Eb_i\right]_{\cb},\Eb_i\right]_{\cb},\Eb_i\right]_{\cb}\in\cJ(\schi).$$
The remaining case is $m_{pi}=2$, for which there exist two possible diagrams:
$$ \xymatrix{\circ^{-\zeta} \ar@{-}[r]^{\zeta^7} & \circ^{\zeta^3}},\zeta\in\G_{9};\qquad \xymatrix{\circ^{-\zeta} \ar@{-}[r]^{-\zeta^{12}} & \circ^{\zeta^5}},\zeta\in\G_{15}.$$
In both cases $q_{pp}\in\G_3$, and also
$$ s_p(3\alpha_i+2\alpha_p)=3\alpha_i+4\alpha_p\in\Delta_+^{\schi}, \quad s_p(4\alpha_i+3\alpha_p)=4\alpha_i+5\alpha_p\notin\Delta_+^{\schi}. $$
Then \eqref{eqn:relacion (m+1)alpha i+m alpha j, caso 3} is a generator of $\cJ(\schi)$ if $3\alpha_i+5\alpha_p\in\Delta_+^{\schi}$, or \eqref{eqn:relacion con 2alpha i+alpha j, caso 1} is a generator of $\cJ(\schi)$ in other case, so for both braidings $E_p^2E_iE_pE_iE_pE_iE_pE_i$ is a linear combination of greater words modulo $\cJ(\schi)$, and \eqref{eqn:relacion (m+1)alpha i+m alpha j, caso 3} belongs to $\cJ(\schi)$. Therefore
$$T_p(\xb)= a_2^4 \left[ \left[ \left[\Eb_{ppi},\Eb_{pi}\right]_{\cb},\Eb_{pi}\right]_{\cb}, \Eb_{pi}\right]_{\cb}\in\cJ(\schi).$$

Consider now $p=i$, so by Lemma \ref{lema:dos conjuntos de generadores} it is enough to prove that
$$ T_p(\xb')\in\cJ(\schi), \qquad \xb':=\left[E_{jp},\left[E_{jp},\left[E_{jp},E_p\right]_c\right]_c\right]_c.$$

If $m_{pj}=2$, then $s_p(3\alpha_p+2\alpha_j)=\alpha_p+2\alpha_j\in\Delta_+^{\schi}$, $s_p(4\alpha_p+3\alpha_j)=2\alpha_p+3\alpha_j\notin\Delta_+^{\schi}$, so $\underline m_{jp}=2$, and \eqref{eqn:relacion estandar B2} is a generator of $\cJ(\schi)$; then
$$T_p(\xb')= a_2^3a_1 \left[\Eb_{pj},\left[\Eb_{pj},\Eb_j\right]_{\cb}\right]_{\cb}\in\cJ(\schi).$$
If $m_{pj}=3$, then $s_p(3\alpha_p+2\alpha_j)=3\alpha_p+2\alpha_j\in\Delta_+^{\schi}$, $s_p(4\alpha_p+3\alpha_j)=5\alpha_p+3\alpha_j\notin\Delta_+^{\schi}$, so \eqref{eqn:relacion con 2alpha i+alpha j, caso 1} is a generator of $\cJ(\schi)$. By \eqref{eqn:Tp adEi Ep},
$$T_p(\xb')= a_3^3a_2 \left[\Eb_{ppj},\left[\Eb_{ppj},\Eb_{pj}\right]_{\cb}\right]_{\cb}\in\cJ(\schi).$$
If $m_{pj}=4$, then $s_p(3\alpha_p+2\alpha_j)=5\alpha_p+2\alpha_j\in\Delta_+^{\schi}$; moreover, we note that $7\alpha_p+3\alpha_j\notin\Delta_+^{\schi}$ in both cases, and \eqref{eqn:relacion con 2alpha i+alpha j, caso 2} is a generator of $\cJ(\schi)$. By this relation and the quantum Serre relations, $\Eb_p^3\Eb_j\Eb_p^3\Eb_j\Eb_p^2\Eb_j$ is written as a linear combination of greater words modulo $\cJ(\schi)$, so
$$T_p(\xb')= a_4^3a_3 \left[\Eb_{pppj},\left[\Eb_{pppj},\Eb_{ppj}\right]_{\cb}\right]_{\cb}\in\cJ(\schi).$$
\medskip

\noindent\vii The proof is similar to \vi, but more simple: for $p=j$ we have just one possibility, $m_{pi}=1$.
\medskip

\noindent\viii Let $\xb=[E_{iij},E_{3\alpha_i+2\alpha_j}]_c$. Consider $p=j$; we consider $i$ a non-Cartan vertex, because in other case  $\Eb_i\Eb_p\Eb_i^2$ or $\Eb_i^2\Eb_p\Eb_i^2$ can be written as a linear combination of other words using the corresponding quantum Serre relation, and finally $\Eb_i^2\Eb_p\Eb_i^2\Eb_p\Eb_i$ is a linear combination of greater words modulo $\cJ(\schi)$, so the relation is redundant in this case. In consequence we consider $m_{pi}=1$, and for this case $\Eb_p\Eb_i^2\Eb_p\Eb_i^3$ is a linear combination of greater words modulo $\cJ(\schi)$, so
$$ T_p(\xb)=\left[ \left[\Eb_{pi},\Eb_i\right]_{\cb}, \left[\left[\Eb_{pi},\Eb_i\right]_{\cb},\Eb_i \right]_{\cb} \right]_{\cb}\in\cJ(\schi). $$

Let $p=i$; as above, it is enough to prove that
$$ T_p(\xb')\in\cJ(\schi), \qquad \xb':=\left[\left[E_{jp},\left[E_{jp},E_p\right]_c \right]_c, \left[E_{jp},E_p\right]_c \right]_c.$$

If $m_{pj}=2$, then $s_p(3\alpha_p+2\alpha_j)=\alpha_p+2\alpha_j\in\Delta_+^{\schi}$, $s_p(5\alpha_p+3\alpha_j)=\alpha_p+3\alpha_j\notin\Delta_+^{\schi}$, so $(\ad_{\cb}\Eb_j)^3\Eb_p$ (or $\Eb_i^3$) is a generator of $\cJ(\schi)$; therefore
$$T_p(\xb')= a_2^3a_1^2 \left[\left[\Eb_{pj},\Eb_j\right]_{\cb},\Eb_j\right]_{\cb}\in\cJ(\schi).$$
If $m_{pj}=3$, then $s_p(3\alpha_p+2\alpha_j)=3\alpha_p+2\alpha_j\in\Delta_+^{\schi}$, $s_p(5\alpha_p+3\alpha_j)=4\alpha_p+3\alpha_j\notin\Delta_+^{\schi}$, so \eqref{eqn:relacion (m+1)alpha i+m alpha j, caso 3} is a generator of $\cJ(\schi)$. Therefore
$$T_p(\xb')= a_3^3a_2^2 \left[\left[\Eb_{ppj},\Eb_{pj}\right]_{\cb},\Eb_{pj}\right]_{\cb}\in\cJ(\schi).$$
If $m_{pj}=4$, then $s_p(3\alpha_p+2\alpha_j)=5\alpha_p+2\alpha_j\in\Delta_+^{\schi}$, $s_p(5\alpha_p+3\alpha_j)=7\alpha_p+3\alpha_j\notin\Delta_+^{\schi}$, so \eqref{eqn:relacion con 2alpha i+alpha j, caso 2} is a generator of $\cJ(\schi)$. In consequence,
$$T_p(\xb')= a_4^3a_3^2 \left[\left[\Eb_{pppj},\Eb_{ppj}\right]_{\cb},\Eb_{ppj}\right]_{\cb}\in\cJ(\schi).$$
\medskip

\noindent\viv The proof is analogous to the previous one.
\epf
\medskip

Now we are ready to prove that the Lusztig isomorphisms descend to the family of algebras $U(\chi)$, so we will look at the root system of this family of algebras. As we consider finite root systems, they are univocally determined as sets of real roots, and using this result we will obtain the desired Theorem of presentation by generators and relations of Nichols algebras.

\begin{proposition}\label{prop:iso Lusztig para las U}
The morphisms \eqref{eqn:iso lusztig para cUp} induce algebra isomorphisms
$$ T_p, T_p^-: U(\chi) \to U(\schi), $$
such that $T_p T_p^-=T_p^- T_p= \id_{U(\chi)}$.
\end{proposition}
\pf By the definition of the ideals $\cJ(\chi)$ and the previous
Lemmata, $T_p(\cJ(\chi)) \subset \cJ(\schi)$, so there exists an
algebra morphism $T_p: U(\chi) \to U(\schi)$. By $\phi_4^2=\id$ and
the definition of the ideal, $\phi_4(\cJ(\chi))=\cJ(\chi)$, and also
$\varphi_{\lambda}(\cJ(\chi))=\cJ(\chi)$ for any
$\lambda\in(\ku^\times)^\theta$, because the ideal is
$\zt$-homogeneous. By \eqref{eqn:Tp con phi4} we have that
$T_p^-(\cJ(\chi))\subset \cJ(\schi)$, so there exists also an
algebra morphism $T_p^-: U(\chi) \to U(\schi)$, induced by the
corresponding morphism.

These algebras are generated by $E_i$, $F_i$, $L_i$, $K_i$, and the identities $T_p T_p^-=T_p^- T_p=\id$ hold for each one of these elements, so these identities hold for all the elements of these algebras, and these morphisms are isomorphisms.
\epf

This result lets us to prove the main result of this Section. The proof is similar to the one for \cite[Thm.5.25]{A-standard}.

\medskip
\textbf{Proof of Theorem \ref{thm:presentacion minima}.}

Set $\underline\Delta_+^\chi:=\Delta^+(U^+(\chi))\setminus\{N_\alpha\alpha: \ \alpha\in\Delta_+^\chi\}$. By the triangular decomposition, Lemma \ref{lema:kernel derivaciones generado por ad}, Theorem \ref{thm: iso Lusztig-Heck} and Proposition \ref{prop:iso Lusztig para las U}, we have that
\begin{equation}\label{eqn:serie Hilbert Uchi}
 \cH_{U^+(\chi)}= \cH_{U^+_{+p}(\chi)}\bq_{h(E_p)}=s_p(\cH_{U^+_{+p}(\schi)})\bq_{h(E_p)},
\end{equation}
for all $p\in\{ 1,\ldots,\theta\}$, because $\deg(T_p(X))=s_p(\deg X)$ for each homogeneous element $X\in\cU(\chi)$. Recall that $h(E_p)\in\{\ord q_{pp},\infty\}$, so
\begin{equation}\label{eqn:sistema de raices Uchi}
\Delta^+\left(U^+(\chi)\right)= s_p\left(\Delta^+ \left(U^+(\schi) \right)\setminus\{\alpha_p, N_p\alpha_p \} \right) \cup S_p,
\end{equation}
where $S_p=\{\alpha_p\}$, or $S_p=\{\alpha_p, N_p\alpha_p\}$, so $\underline\Delta_+^\chi=s_p\left(\underline\Delta_+^{\schi}\setminus\{\alpha_p\}\right) \cup \{\alpha_p\}$.

In this way, if we consider the sets $\underline\Delta_+^\chi$, for each $\chi$ in a Weyl equivalence class of a fixed braiding with finite root system, then $R=\{\underline\Delta_+^\chi\}$ is a root system for our Weyl groupoid, according with the Definition \ref{def:sistema de raices}. As we have a finite root system, it follows that $\underline\Delta_+^\chi=\Delta_+^\chi$, for all $\chi$, because by Proposition \ref{Prop:maxlongitudCH} all the roots are real. In this way, $\Delta^+(U^+(\chi))$ is obtained from $\Delta_+^\chi$ adding $N_\alpha\alpha$, for some $\alpha \in \Delta_+^\chi$. Fix an order on the letters $x_i$ and consider the corresponding PBW basis. We have a projection $\pi_{\chi}: U(\chi)\twoheadrightarrow\u(\chi)$ of graded braided Hopf algebras, so the corresponding $x_{\alpha}$ of the PBW basis of $\u(\chi)$ are generators of the PBW basis of $U(\chi)$, by the definition of Kharchenko's PBW basis. On the other hand, each simple root of a non-Cartan vertex satisfies $E_i^{N_i}=0$ in $U(\chi)$, so $N_i\alpha_i\notin\Delta^+\left(U^+(\chi)\right)$. Therefore \eqref{eqn:serie Hilbert Uchi} implies that
$$ N_\alpha \alpha \notin \Delta^+\left(U^+(\chi)\right), \qquad \mbox{for all }\alpha \in \Delta_+^{\chi}\setminus\cO(\chi),$$
because $\alpha$ is of the way $\alpha=w(\alpha_i)$ for some $w\in\cW$ and $i\in\{1,\ldots,\theta\}$, $i$ a non-Cartan vertex in the corresponding $\chi'$. Analogously, for each Cartan vertex $i$, $N_i\alpha_i\in\Delta^+\left(U^+(\chi)\right)$, because $E_i^{N_i}\neq0$ in $U(\chi)$, so
$$ N_\alpha \alpha \in \Delta^+\left(U^+(\chi)\right), \qquad \mbox{for all }\alpha \in \cO(\chi).$$
Therefore $\Delta^+(U^+(\chi))= \Delta_+^{\chi} \cup \{ N_\alpha \alpha: \ \alpha \in \cO(\chi) \}$.

Suppose that the degree $N_\alpha \alpha$ in $\Delta^+(U^+(\chi))$ corresponds to a Lyndon word of this degree: we can assume that it is of minimal length,
and we denote it by $u$; set $(v,w)=\Sh(u)$. In this way, $\deg v=\beta$, $\deg w=\gamma$, for some $\beta$, $\gamma\in\Delta_+^{\chi}$, and
$\beta+\gamma=N_\alpha \alpha$. As all the roots are real, we deduce that if $\beta<\gamma$, then $\beta<\alpha<\gamma$, by a similar argument to the
convexity properties in \cite{A-presentation}. We can consider then the case $\beta=\alpha_i$, because if $\beta=s_{i_1}\cdots s_{i_k}(\alpha_{i_{k+1}})$,
where $w=s_{i_1}\cdots s_{i_k}$ is the beginning of the expression of the element of maximal length, we apply $w^{-1}$ to obtain that
$\alpha_{i_{k+1}}+\gamma'=N_\alpha \alpha'$ for some $\alpha', \gamma'\in \Delta_+^{\chi}$. Note also that $N_\alpha>2$, because if we suppose
$N_\alpha=2$, then $\alpha$ is applied in a simple root $\alpha_i$ corresponding a Cartan vertex by some element of the Weyl groupoid, and as
$N_{\alpha}$ is invariant by the action of the Weyl groupoid, it should be $q_{ii}=-1$, but it corresponds to an isolated vertex or a non-Cartan vertex,
which is a contradiction. Set
$$\alpha=\sum_{j=1}^\theta n_j\alpha_j, \quad \gamma=\sum_{j=1}^\theta m_j\alpha_j, \qquad \mbox{for some }n_j,m_j\in\N_0.$$
Note that $m_i=N_{\alpha}n_i-1\geq 2$, and for $j\neq i$, $m_j=N_{\alpha}n_j\geq 3$, so $\sop\gamma=\sop\alpha$. By simplicity assume that $\sop\alpha=\{1, \ldots, \theta\}$; note that the vertices of $\sop\beta$ corresponds to a connected subdiagram, for any positive root $\beta$.

By these considerations we reduce the problem to an analysis case by case of the positive root systems of connected diagrams, and we do it with the help of the program \emph{SARNA} \cite{GHV}. We look for the possible $\gamma$ such that all the coordinates, except at most one, are divisible by an integer $\geq3$, and the remaining coordinate is $\geq 2$, so we just have a few 3-uples in rank two or three. Analyzing each of these 3-uples
$$(\alpha,\gamma,i)\in \Delta_+^\chi\times\Delta_+^\chi\times\{1, \ldots, \theta\} \mbox{ such that there exists }N\in\N: \ \alpha_i+\gamma=N\alpha, $$
we note that $N\neq N_\alpha$ for all of them. Therefore, there are no Lyndon words of degree $N_\alpha \alpha$, so the generators of degree $N_\alpha \alpha$ are $x_\alpha^{N_\alpha}$, and then the elements
$$ x_{\beta_1}^{n_1} \cdots x_{\beta_k}^{n_k}, \qquad \beta_i \in \Delta_+^{\chi}, \quad \left\{\begin{array}{ll} 0\leq n_j <N_j , &
\hbox{if }\beta_j \notin \cO(\chi) \\0\leq n_j <\infty , & \hbox{si
}\beta_j \in \cO(\chi)\end{array}\right. $$ are a PBW basis of
$U(\chi)$. As $x_\alpha^{N_\alpha}=0$ in $\u^+(\chi)$, $\pi_\chi$
induces a surjective morphism
$$ \pi_\chi ': U(\chi) / \langle x_\alpha^{N_\alpha}: \ \alpha \in \cO(\chi) \rangle \ \longrightarrow \ \u^+(\chi),$$
which applies the set
$$ \{ x_{\beta_1}^{n_1} \cdots x_{\beta_k}^{n_k}, \qquad \beta_i \in \Delta_+^{\chi}, 0\leq n_j <N_j \}, $$
generating linearly the quotient, to the corresponding PBW basis of $\u^+(\chi)$. Therefore $\pi_\chi'$ is an isomorphism.
\qed

\section{Generation in degree one}\label{section:conjetura AS}

Now we answer positively the Conjecture \ref{conj:AS, generacion en
grado 1}, formulated by Andruskiewitsch and Schneider, but
restricting to the case in which $G(H)$ is abelian. The technique of
the proof is the same that these authors use in \cite{AS Class},
extended in some works to other families. In particular, the first
Lemmata of this Section correspond to relations generating the ideal
for standard braidings as in \cite{AnGa}, but the proof is made in a
general context.
\medskip

In what follows $\Gamma$ denotes a finite abelian group, and $S= \bigoplus_{n \geq 0} S(n)$ is a graded braided Hopf algebra in $\gyd$ such that
$S(0)=\ku 1$, generated as an algebra by $V:=S(1)$. Fix a basis $\{x_1,\ldots, x_{\theta}\}$ of $V$, so $V$ has a braiding of diagonal type: we can assume
that $x_i\in S(1)^{\chi_i}_{g_i}$ for some $g_i\in\Gamma$ and $\chi_i\in\widehat{\Gamma}$. Set then $q_{ij}:=\chi_j(g_i)\in\ku^\times$.

We will prove that if $S$ is finite dimensional, then $S$ is the Nichols algebra $\cB(V)$ associated to $V$. We will obtain then the main Theorem of this
Section, answering this Conjecture.
\smallskip

We begin by extending \cite[Lemma 5.4]{AS Class} for a general quantum Serre relation, proving that they hold in $S$, or $S$ is infinite-dimensional.

\begin{proposition}\label{prop:gen-relacion qSerre} Let $S$ be as above. If there exist $i,j\in\{1,\ldots,\theta\}$ such that $q_{ii}^{m_{ij}+1}\neq1$, and also $\ad_c(x_i)^{1+m_{ij}}(x_j)\neq0$, then $S$ is infinite-dimensional.
\end{proposition}
\pf
By definition of $m_{ij}$, we have that $q_{ii}^{m_{ij}}\widetilde{q_{ij}}=1$. We begin the proof as in \cite[Lemma 5.4]{AS Class}. To simplify the notation, call $m=m_{ij}$, $q=q_{ii}$, $y_1:=x_i$, $y_2:=x_j$, $y_3:= \ad_c(x_i)^{1+m}(x_j)$. Set also
\begin{align*}
&h_1=g_i, &&h_2=g_j, &&h_3=g_i^{m+1}g_j, \\ &\eta_1=\chi_i, &&\eta_2=\chi_j, &&\eta_3= \chi_i^{m+1}\chi_j,
\end{align*}
so $y_k\in S^{\eta_k}_{h_k}$, $1\leq k\leq3$. If $W=\ku y_1+\ku y_2+\ku y_3$, then $W\subset\cP(S)$ (because $y_3$ is also primitive), so there exists a monomorphism $\cB(W)\hookrightarrow S$. We compute the corresponding braiding matrix $\left(Q_{kl}=\eta_l(h_k)\right)_{1\leq k,l\leq3}$, and consider the corresponding generalized Dynkin diagram:
\begin{equation}\label{diagram:Qserre}
\xymatrix{ & \circ^{q_{jj}} \ar@{-}[rd]^{q^{-m(m+1)}q_{jj}^{2}} &  \\  \circ^{q}\ar@{-}[ru]^{q^{-m}} \ar@{-}[rr]_{q^{m+2}} & & \circ^{q^{m+1}q_{jj}.} }
\end{equation}

We will consider the different possible cases and prove that no one of them are in \cite{H-classif RS}, so $\cB(W)$, and in consequence $S$, is
infinite-dimensional. Suppose that the diagram \eqref{diagram:Qserre} is in Heckenberger's list:
\smallskip

\noindent\textbf{\emph{Case I:}} $Q_{kl}Q_{lk}\neq1$ for all $1\leq k<l\leq 3$. By \cite[Lemma 9]{H-classif RS}, $1=\prod_{k<l}Q_{kl}Q_{lk}=q^{2-m(m+1)}q_{jj}^2$, and at least one of the vertices is labeled with $-1$. Note that $q\neq-1$, because in such case $m=1$ (and we assume $q^{m+1}\neq1$). Also $q_{jj}\neq q^{m+1}q_{jj}$ by hypothesis, so there is only one vertex labeled with $-1$.
\begin{itemize}
  \item If $q_{jj}=-1$, then $1=(q^{m+1}q_{jj})(q^{-m(m+1)}q_{jj}^2)=-q^{1-m^2}$ and $m=1$ by the same Heckenberger's Lemma, which is a contradiction.
  \item If $q^{m+1}q_{jj}=-1$, then $1=qq^{m+2}=q^{m+3}$ by the same result, and also
  $$1=q_{jj}(q^{-m(m+1)}q_{jj}^{2})=q_{jj}^3q^{-m(m+3)+2m}=q_{jj}^3q^{2m}, $$
so we deduce that
$$ -1= (-1)^3= q_{jj}^{3}q^{3m+3}= (q_{jj}^3q^{2m})q^{m+3}=q^{m+3}, $$
which is also a contradiction. Therefore \eqref{diagram:Qserre} does not belong to Heckenberger list for this case.
\end{itemize}
\medskip

\noindent\textbf{\emph{Case II}:} $Q_{12}Q_{21}=q^{-m}=1$. In this case $m=0$, so \eqref{diagram:Qserre} becomes:
\begin{equation}\label{diagram:symmetriccase}
 \xymatrix{ \circ^{q} \ar@{-}[r]_{q^2} & \circ^{qq_{jj}} \ar@{-}[r]_{q_{jj}^2} & \circ^{q_{jj}}. }
\end{equation}
If $q_{jj}=-1$ we obtain the connected subdiagram $\xymatrix{\circ^{q}\ar@{-}[r]_{q^2}&\circ^{-q}}$, which has no vertices labeled with $-1$, and these labels are different. Such diagram is not of finite Cartan type and moreover it does not correspond to any diagram without $-1$ in the vertices in rows 5, 9, 11, 12, 15 of \cite[Table 1]{H-classif RS}, so $\cB(W)$ is infinite-dimensional.

If $q_{jj}\neq-1$ and $q=-1$, we have an analogous situation, so $q\neq-1$ and \eqref{diagram:symmetriccase} is a connected diagram of three vertices. If $qq_{jj}\neq-1$, then \cite[Lemma 9]{H-classif RS} implies that one of the following situations holds:
\begin{itemize}
  \item the diagram is of finite Cartan type, so it contains a subdiagram of Cartan type $A_2$. Then $1=qq^2=(qq_{jj})q^2$, or $1=q_{jj}q_{jj}^2=(qq_{jj})q_{jj}^2$, so $q=1$ or $q_{jj}=1$;
  \item $q^3=1$, $q_{jj}$, $q_{jj}q\in\G_6\cup\G_9$, and $q_{jj}q_{jj}^2=1$ or $q_{jj}^3=1$, $q,q_{jj}q\in\G_6\cup\G_9$, $qq^2=1$.
\end{itemize}
No one of these situations hold, so $qq_{jj}=-1$. If this diagram is in \cite[Table 2]{H-classif RS}, it follows that $Q_{ii}Q_{i2}Q_{2i}=1$ for some $i\in\{1,3\}$ in all the possible cases. We can assume then $i=1$, $q^3=1$. By \cite[Lemma 9]{H-classif RS}, one of the following situations holds:
\begin{itemize}
  \item $q_{jj}^3=1$, but also $q_{jj}^3=-q^{-3}=-1$,
  \item $q_{jj}^4=1$,
  \item $q_{jj}=-q$.
\end{itemize}
No one of these situations are possible, so we obtain a contradiction in this case too.
\medskip

\noindent\textbf{\emph{Case III:}} $Q_{13}Q_{31}=q^{m+2}=1$. We obtain the diagram:
$$ \xymatrix{ \circ^{q} \ar@{-}[r]_{q^2} & \circ^{q_{jj}} \ar@{-}[r]_{q^{-2}q_{jj}^2} & \circ^{q^{-1}q_{jj}}. } $$
Such diagram is the corresponding to \eqref{diagram:symmetriccase}, but changing $q_{jj}$ by $q_{jj}q^{-1}$, so it does not belong to \cite[Table 2]{H-classif RS}. Then $q^{m+2} \neq 1$.
\medskip

\noindent\textbf{\emph{Case IV:}} $Q_{23}Q_{32}=1$. This means $q_{jj}^2= q^{m(m+1)}$, so we have the diagram:
\begin{equation}\label{diagram:lastcase}
 \xymatrix{ \circ^{q_{jj}} \ar@{-}[r]_{q^{-m}} & \circ^{q} \ar@{-}[r]_{q^{m+2}} & \circ^{q^{m+1}q_{jj}}. }
\end{equation}
This diagram is connected by the previous cases. As $m\neq 0$, $q^{m+1}\neq1$, it follows that $q\neq{-1}$. Consider the different possible values of the labels of the vertices:

$\mathbf{q_{jj}=q^{m+1}q_{jj}=-1}$: that is, $q^{m+1}=1$ and we have the diagram:
$$ \xymatrix{ \circ^{-1} \ar@{-}[r]_{q} & \circ^{q} \ar@{-}[r]_{q} & \circ^{-1}, } $$
which is not in Heckenberger's list.
\smallskip

$\mathbf{q_{jj}=-1, q^{m+1}q_{jj}\neq -1}$: By \cite[Table 2]{H-classif RS}, it should be $1=Q_{22}Q_{23}Q_{32}=q^{m+3}$, and we should have the diagram
  $$ \xymatrix{ \circ^{-q^{-2}} \ar@{-}[r]_{q^{-1}} & \circ^{q} \ar@{-}[r]_{q^3} & \circ^{-1}. } $$
Moreover, $1=q_{jj}^2= q^{m(m+1)}=q^{2m}=q^{-6}$. Note that $q^3\neq1$ because $q^m\neq1$, so $q\in\G_6$. But this diagram is not in Heckenberger's list.
\smallskip

$\mathbf{q_{jj}\neq -1, q^{m+1}q_{jj}=-1}$: as above, $1=Q_{22}Q_{21}Q_{12}=q^{1-m}$. By definition it should be $m=1$, with the same diagram of the previous case and $q\in \G_6$, so we obtain the same contradiction.
\smallskip

$\mathbf{q_{jj}, q^{m+1}q_{jj} \neq -1}$: By \cite[Lemma 9]{H-classif RS}, one of the following situations holds:
\begin{itemize}
  \item the diagram is of Cartan type. Then, $q=q_{jj}$ and $m=1$, or $q=q^{m+1}q_{jj}=q^{-m-2}$. In both cases we obtain the same diagram,
  $$ \xymatrix{ \circ^{q}   \ar@{-}[r]_{q^{-1}} & \circ^{q} \ar@{-}[r]_{q^3} & \circ^{q^3}. } $$
We discard easily the cases $A_3,C_3$, because $q,q^2\neq q^3$. If it is of type $B_3$, $q=(q^3)^2=q^{-3}$, which is a contradiction.
  \item $q_{jj}\in\G_3$, $q\in\G_6\cup\G_9$ and $1=q^{1-m}=q_{jj}q^{2m+3}$. Then $m=1$ and $q^5=q_{jj}^{-1}$, so $q^{15}=1$, but we obtain then a contradiction with the fact that $q\in\G_6 \cup \G_9$ is primitive.
  \item $q^{m+1}q_{jj}\in\G_3$, $q\in\G_6\cup\G_9$ y $1=q_{jj}q^{-m}=q^{m+3}$. Again $q^{15}=1$, and we obtain the same contradiction.
\end{itemize}
In consequence, \eqref{diagram:Qserre} is not in Heckenberger's list, and $S$ is infinite-dimensional.
\epf

Now we continue with another Lemmata from \cite{AnGa}, just adapted to this general context.

\begin{lemma}\label{Lema:gen-relacion vertice -1} Let $j,k,l\in\{1,\ldots,\theta\}$ be such that $q_{kk}=-1$, $\widetilde{q_{kj}}=\widetilde{q_{kl}}^{-1}\neq 1$, $\widetilde{q_{jl}}=1$. If $\left[x_{jkl},x_k\right]_c\neq0$ is a primitive element of $S$, then $S$ is infinite-dimensional.
\end{lemma}
\pf
Set $u:=\left[x_{jkl},x_k\right]_c$, $g_u:=g_jg_k^{2}g_l \in \Gamma$, $\chi_u:=\chi_j\chi_k^{2}\chi_l\in\widehat{\Gamma}$, $q:=\widetilde{q_{lk}}$; we work then as in the previous Lemma.

We compute the braiding corresponding to the primitive elements $y_1=x_j$, $y_2=x_k$, $y_3=x_l$ and $y_4=u$, with the corresponding elements $h_i \in \Gamma$, $\eta_i \in \widehat{\Gamma}$; we will prove that such braiding has an associated Nichols algebra of infinite dimension, and so $S$ has infinite dimension. The corresponding generalized Dynkin diagram to $(Q_{rs}=\eta_s(h_r))_{1 \leq r,s \leq 4}$ is:
\begin{equation}\label{diagram:gen-relacion vertice -1}
    \xymatrix{ \circ^{q_{jj}} \ar@{-}[r]_{q^{-1}} \ar@{-}[d]_{q_{jj}^{2}q^{-2}} & \circ^{-1} \ar@{-}[d]_{q} \\
    \circ^{q_{jj}q_{ll}} \ar@{-}[r]_{q_{ll}^{2}q^2} & \circ^{q_{ll}}.  }
\end{equation}
Suppose that such diagram is in Heckenberger's list. If $q=-1$, then \eqref{diagram:gen-relacion vertice -1} contains \eqref{diagram:symmetriccase} as a subdiagram, so it does not appear in the list. Therefore $q\neq-1$. As each diagram in \cite[Table 3]{H-classif RS} does not contain a 4-cycle, it follows that $q_{jj}^2q^{-2}=1$, or $q_{ll}^2q^2=1$. As the conditions are symmetric, it is enough to consider the case $q_{jj}=\pm q$.

If we also have $q_{ll}=\pm q^{-1}$, and as $Q_{44}=q_{jj}q_{ll}\neq1$, the diagram contains the following
$$\xymatrix{\circ^{q} \ar@{-}[r]_{q^{-1}} & \circ^{-1} \ar@{-}[r]_{q} & \circ^{-q^{-1}} },$$
which is a contradiction with \cite[Lemma 9]{H-classif RS}. In consequence we have:
$$ \xymatrix{\circ^{\pm q} \ar@{-}[r]_{q^{-1}} & \circ^{-1} \ar@{-}[r]_{q} & \circ^{q_{ll}} \ar@{-}[r]_{q_{ll}^2q^2} & \circ^{q_{jj}q_{ll}}. } $$
Suppose that $q_{jj}=-q$. As $Q_{11}Q_{12}Q_{21}\neq 1$, we deduce from \cite[Table 3]{H-classif RS} that $m_{12}=2$; that is, $$0=(1-Q_{11}^3)(Q_{11}^2Q_{12}Q_{21}-1)=(1+q^3)(q-1),$$
which gives conditions about $q$, but each diagram in \cite[Lemma 9]{H-classif RS} does not satisfy this condition.

Therefore $q_{jj}=q$. We look at \cite[Table 3]{H-classif RS} but a diagram in such list does not satisfy $Q_{22}=-1$, $Q_{11}=Q_{44}Q_{33}^{-1}=q\neq\pm1$, so \eqref{diagram:gen-relacion vertice -1} is not in the list. In consequence, $S$ has infinite dimension.
\epf
\medskip

\begin{lemma}\label{Lema:gen-relaciones tipo B}
\vi Let $i,j\in\{1,\ldots,\theta\}$ be such that $q_{jj}=-1$, $q_{ii}\widetilde{q_{ij}}\in\G_6$, and also $q_{ii}\in\G_3$ or $m_{ij}\geq 3$. If
$\left[x_{iij},x_{ij}\right]_c\in\cP(S)\setminus\{0\}$, then $S$ is infinite-dimensional.

\noindent \vii Let $i,j,k\in\{1,\ldots,\theta\}$ be such that $q_{ii}=\pm \widetilde{q_{ij}}\in\G_3$, $\widetilde{q_{ik}}=1$, and also
$-q_{jj}=\widetilde{q_{ij}}\widetilde{q_{jk}}=1$ or $q_{jj}^{-1}=\widetilde{q_{ij}}=\widetilde{q_{jk}}\neq -1$. If $\left[x_{iijk},x_{ij}\right]_c
\in\cP(S)\setminus\{0\}$, then $S$ is infinite-dimensional.
\end{lemma}
\pf
\vi We follow the same scheme of proof. Set
$$y_1=x_i, \quad y_2=x_j, \quad y_3=\left[x_{iij},x_{ij}\right]_c,$$
and $h_i\in\Gamma$, $\eta_i\in\widehat{\Gamma}$, $i=1,2,3$ the corresponding elements. Suppose that the braiding matrix
$(Q_{rs}=\eta_s(h_r))_{1\leq r,s\leq 3}$ appears in Heckenberger's list. The associated generalized Dynkin diagram is
$$ \xymatrix{ \circ^{q_{ii}} \ar@{-}[rr]_{q} \ar@{-}[rd]_{q^2} &  & \circ^{-1} \\ & \circ^{q_{ii}^3} \ar@{-}[ru]_{q^3} & }, \qquad q:=\widetilde{q_{ij}}. $$
Then $Q_{33}=q_{ii}^3\neq 1$, so $m_{ij}\geq3$. Moreover the diagram is connected, so it is of type super $G(3)$, the unique diagram of rank three such
that some $m_{rs}$ is $\geq3$. Therefore $1=Q_{23}Q_{32}=q^3$,  which is a contradiction, so the diagram associated to $(Q_{rs})$ does not correspond to a
finite-dimensional Nichols algebra. In consequence $S$ is infinite-dimensional.
\smallskip

\noindent \vii Set $w:=\left[x_{iijk},x_{ij}\right]_c$, and denote as above $y_1=x_i$, $y_2=x_j$, $y_3=x_k$, $y_4=w$,
$W$ the subspace generated by these elements, and $h_i \in \Gamma$, $\eta_i \in \widehat{\Gamma}$, $i=1,2,3,4$ the corresponding elements: suppose again
that $\cB(W)$ is a finite-dimensional Nichols algebra. Set $\zeta=q_{ii}\in\G_3$. We analyze each possible case.
\begin{itemize}
  \item $q_{jj}=-1$, $\widetilde{q_{ij}}\widetilde{q_{jk}}=1$: the diagram of $(Q_{rs})$ becomes
  $$ \xymatrix{ \circ^{\zeta} \ar@{-}[r]^{\pm\zeta} \ar@{-}[rd]_{\zeta^2} & \circ^{-1} \ar@{-}[r]^{\mp\zeta^2}\ar@{-}[d]^{\zeta^2} & \circ^{q_{kk}}
\ar@{-}[ld]^{q_{kk}^2\zeta}  \\ & \circ^{q_{kk}\zeta}  &  }. $$
As $Q_{12}Q_{21}$, $Q_{14}Q_{41}$, $Q_{42}Q_{24}\neq 1$, and the product of these three scalars is not $1$, such diagram is not in Heckenberger's list,
by \cite[Lemma 9]{H-classif RS}.
\item $q_{jj}^{-1}=\widetilde{q_{ij}}=\widetilde{q_{jk}}\neq -1$: now we have the diagram
  $$ \xymatrix{ \circ^{\zeta} \ar@{-}[r]^{\pm\zeta} \ar@{-}[rd]_{\zeta^2} & \circ^{\pm\zeta^2} \ar@{-}[r]^{\pm\zeta} & \circ^{q_{kk}}
\ar@{-}[ld]^{q_{kk}^2\zeta^2}  \\   & \circ^{q_{kk}\zeta}  &  }. $$
The lack of 4-cycles in Heckenberger's list implies that $1=Q_{34}Q_{43}=q_{kk}^2\zeta^2$, so $q_{kk}\zeta=-1$, because $Q_{44}=q_{kk}\zeta\neq1$.
But this diagram does not appear in \cite[Table 3]{H-classif RS}.
\end{itemize}
We obtain a contradiction in all the cases, so $S$ is infinite-dimensional.
\epf
\medskip

\begin{lemma}\label{lema:gen-relacion triangulo}
Let $i,j,k\in\{1,\ldots,\theta\}$ be such that $\widetilde{q_{ik}}$, $\widetilde{q_{ij}}$, $\widetilde{q_{jk}} \neq 1$. Let
$$ w:=x_{ijk}-\frac{1-\widetilde{q_{jk}}}{q_{kj}(1-\widetilde{q_{ik}})}\left[x_{ik},x_j\right]_c-q_{ij}(1-\widetilde{q_{jk}}) \ x_jx_{ik}.$$
If $w \in\cP(S)\setminus\{0\}$, then $S$ is infinite-dimensional.
\end{lemma}
\pf Set $y_1=x_i$, $y_2=x_j$, $y_3=x_k$, $y_4=w$, $W$ the subspace
generated by these elements, $h_i \in \Gamma$, $\eta_i \in
\widehat{\Gamma}$, $i=1,2,3,4$ the corresponding elements, and
$(Q_{rs}=\eta_s(h_r))_{1\leq r,s\leq4}$: suppose as above that
$\cB(W)$ is finite dimensional. Note that
$$Q_{14}Q_{41}=q_{ii}^2\widetilde{q_{ij}}\widetilde{q_{ik}}=q_{ii}^2\widetilde{q_{jk}}^{-1}, $$
because $\widetilde{q_{ij}}\widetilde{q_{ik}}\widetilde{q_{jk}}=1$, by \cite[Lemma 9]{H-classif RS}. By the same Lemma at least one vertex is labeled with $-1$. Then, if
$q_{ii}=-1$, we have that $Q_{14}Q_{41}\neq-1$; the same holds for the other vertices, so exactly one vertex is labeled with $-1$ (we have no 4-cycles).
We look for possible braiding with these conditions in \cite[Table 3]{H-classif RS}, but no one coincides with this description. Therefore $\cB(W)$ is
infinite-dimensional, and $S$ too.
\epf

\begin{lemma}\label{lema:gen-relaciones super C3-G3}
\vi Let $i,j,k\in\{1,\ldots,\theta \}$ be such that one of the following conditions holds:
\begin{itemize}
 \item $q_{ii}=q_{jj}=-1$, $\widetilde{q_{ij}}^2= \widetilde{q_{jk}}^{-1}$, $\widetilde{q_{ik}}=1$, or
 \item $\widetilde{q_{ij}}=q_{jj}=-1$, $q_{ii}= -\widetilde{q_{jk}}^2\in\G_3$, $\widetilde{q_{ik}}=1$, or
 \item $q_{kk}=\widetilde{q_{jk}}=q_{jj}=-1$, $q_{ii}= -\widetilde{q_{ij}}\in\G_3$, $\widetilde{q_{ik}}=1$, or
 \item $q_{jj}=-1$, $\widetilde{q_{ij}}=q_{ii}^{-2}$, $q_{kk}=\widetilde{q_{jk}}^{-1}=-q_{ii}^{3}$, $\widetilde{q_{ik}}=1$, or
 \item $q_{ii}=q_{jj}=q_{kk}=-1$, $\pm\widetilde{q_{ij}}=\widetilde{q_{jk}}\in\G_3$,
 $\widetilde{q_{ik}}=1$,
\end{itemize}
If $\left[\left[x_{ij},x_{ijk}\right]_c,x_j\right]_c\in\cP(S)\setminus\{0\}$, then $S$ is infinite-dimensional.
\smallskip

\noindent\vii Let $i,j,k\in\{1,\ldots,\theta\}$ be such that $q_{ii}=q_{jj}=-1$, $(\widetilde{q_{ij}})^3=(\widetilde{q_{jk}})^{-1}$, $\widetilde{q_{ik}}=1$. If
$\left[\left[x_{ij},\left[x_{ij},x_{ijk}\right]_c\right]_c,x_j\right]_c\in\cP(S)\setminus\{0\}$, then $S$ is infinite-dimensional.
\end{lemma}
\pf\vi Set $y_1=x_i$, $y_2=x_j$, $y_3=x_k$, $y_4=\left[\left[x_{ij},x_{ijk}\right]_c,x_j\right]_c$, $W$ the subspace generated
by these elements, $h_i \in \Gamma$, $\eta_i \in \widehat{\Gamma}$, $i=1,2,3,4$ the associated elements and $(Q_{rs}=\eta_s(h_r))_{1\leq r,s\leq4}$ the
braiding matrix. We will consider the associated generalized Dynkin diagram for each case.

For the first case, we have the following diagram, where $q:=\widetilde{q_{ij}}$:
$$ \xymatrix{ \circ^{-1} \ar@{-}[r]^{q} \ar@{-}[rd]_{q^3} & \circ^{-1} \ar@{-}[r]^{q^{-2}} & \circ^{q_{kk}} \ar@{-}[ld]^{q_{kk}^2q^{-6}}
\\ & \circ^{-q_{kk}}  &  }. $$
Suppose that $\cB(W)$ is finite-dimensional. Then $Q_{33}Q_{32}Q_{23}=1$, so $q_{kk}=q^2$ and then $Q_{34}Q_{43}=q^{-2}\neq1$. In consequence such diagram
is of type super $F(4)$. Then $1=Q_{14}Q_{41}=q^3$, which is a contradiction.

For the second case, $Q_{12}Q_{21}=q_{ii}\in\G_3$, and $Q_{14}Q_{41}=-1$, so the diagram contains a 4-cycle and then $\cB(W)$ is infinite-dimensional. An analogous situation holds for the third case, because $Q_{12}Q_{21}=-q_{ii}\in\G_6$, and $Q_{14}Q_{41}=-1$.

For the fourth case, $\widetilde{Q_{23}}=\widetilde{Q_{24}}=\widetilde{Q_{34}}=\widetilde{q_{kj}}\neq 1$ and $\widetilde{Q_{14}}=\widetilde{Q_{12}}=\widetilde{q_{ij}}\neq 1$, so $\cB(W)$ is infinite-dimensional.

For the last case, $Q_{44}=1$, and then $\cB(W)$ is infinite-dimensional.

Therefore $S$ is infinite-dimensional in all the cases.
\medskip

\noindent\vii We use the same notation, but in this case $ y_4=\left[\left[x_{ij},\left[x_{ij},x_{ijk}\right]_c\right]_c,x_j\right]_c$. So we have the following diagram for $(Q_{rs})$:
$$\xymatrix{\circ^{-1} \ar@{-}[r]^{q}\ar@{-}[rd]_{q^4} &\circ^{-1}\ar@{-}[r]^{q^{-3}}&\circ^{q_{kk}}\ar@{-}[ld]^{q_{kk}^2q^{-6}}
\\ & \circ^{-q_{kk}} & }, $$
where $q=\widetilde{q_{ij}}$. Suppose that $\cB(W)$ is finite-dimensional. By \cite[Table 3]{H-classif RS}, this diagram cannot be connected. In consequence,
$1=Q_{41}Q_{14}=Q_{34}Q_{43}$, so $q_{kk}=\pm1$. But then $Q_{33}=1$, or $Q_{44}=1$, which is a contradiction to the fact that $\cB(W)$ is
finite-dimensional. So $S$ is infinite-dimensional.
\epf
\medskip

\begin{lemma}\label{lema:gen-relacion C4}
Let $i,j,k,l\in\{1,\ldots,\theta\}$ be such that $q_{jj}\widetilde{q_{ij}}= q_{jj}\widetilde{q_{jk}}=1$, $\widetilde{q_{jk}}^2= \widetilde{q_{kl}}^{-1}= q_{ll}$, $q_{kk}=-1$,
$\widetilde{q_{ik}}=\widetilde{q_{il}}=\widetilde{q_{jl}}$. If $\left[ \left[ \left[x_{ijkl},  x_k \right]_c, x_j \right]_c, x_k\right]_c\in\cP(S)\setminus\{0\}$, then $S$ is infinite-dimensional.
\end{lemma}
\pf
We use a similar notation and consider the corresponding subspace $W$ generated by the corresponding primitive elements. Suppose that $\cB(W)$ is
finite-dimensional. Its associated Dynkin diagram is
$$ \xymatrix{\circ^{q_{ii}} \ar@{-}[r]^{q^{-1}}\ar@{-}[rd]_{q_{ii}^2q^{-2}} &\circ^{q}\ar@{-}[r]^{q^{-1}} &\circ^{-1}\ar@{-}[r]^{q^2}
&\circ^{q^{-2}}\ar@{-}[lld]^{q^2}\\ & \circ^{-q_{ii}} & & }, \qquad q=q_{jj}.$$
Note that $1=Q_{15}Q_{51}$, because there are no 5-cycles and $q^2\neq 1$. Therefore $q_{ii}=q$, but this diagram is not in Heckenberger's list, so
$\cB(W)$ is infinite-dimensional, and $S$ too.
\epf
\medskip

\begin{lemma}\label{lema:gen-relaciones C3,G3 raices chicas}
Let $i,j,k \in \{1, \ldots, \theta \}$ be such that $q_{jj}=\widetilde{q_{ij}}^{-1}=\widetilde{q_{jk}}$.

\noindent\vi If $q_{jj}\in \G_3$ and $\left[\left[x_{ijk},x_j \right]_c x_j \right]_c\in\cP(S)\setminus\{0\}$, then $S$ is infinite-dimensional.

\noindent\vii If $q_{jj}\in \G_4$ and $\left[\left[\left[x_{ijk},x_j\right]_c,x_j \right]_c,x_j \right]_c\in\cP(S)\setminus\{0\}$, then $S$ is infinite-dimensional.
\end{lemma}
\pf
\vi Using the same notation as in previous Lemmata, we have the diagram
$$\xymatrix{\circ^{q_{ii}}\ar@{-}[r]^{\zeta^2}\ar@{-}[rd]_{q_{ii}^2}&\circ^{\zeta}\ar@{-}[r]^{\zeta}&\circ^{q_{kk}}\ar@{-}[ld]^{q_{kk}^2}
\\ & \circ^{q_{kk}q_{ii}} & },\qquad \zeta=q_{jj}\in\G_3, $$
for $y_1=x_i$, $y_2=x_j$, $y_3=x_k$, $y_4=\left[\left[x_{ijk},x_j \right]_c x_j \right]_c$, with corresponding matrix $(Q_{rs})$, and $W$
is the subspace generated by these elements. Note that
\begin{itemize}
  \item if $q_{ii}=q_{kk}=-1$, then $Q_{44}=1$;
  \item if $q_{ii},q_{kk}\neq-1$, then the diagram contains a 4-cycle;
  \item if $q_{ii}=-1$, $q_{kk}\neq-1$, or $q_{ii}\neq-1$, $q_{kk}=-1$, the diagram contains $\xymatrix{\circ^{q}\ar@{-}[r]^{q^2}&\circ^{-q}}$ as a
subdiagram (where $q=q_{ii}$ or $q=q_{kk}$), and this connected subdiagram of rank two is not in \cite[Table 1]{H-classif RS}.
\end{itemize}
In all the cases $\cB(W)$  is infinite-dimensional, so $S$ too.
\medskip

\noindent\vii The proof is analogous.
\epf
\medskip

\begin{lemma}\label{lemma:gen-casos F4 y C4 modificado}
\vi Let $i,j,k,l\in\{1,\ldots,\theta \}$ be such that
$q_{ll}=\widetilde{q_{lk}}^{-1}=
q_{kk}=\widetilde{q_{jk}}^{-1}=q^2$, $\widetilde{q_{ij}}=
q_{ii}^{-1}=q^3$ for some $q\in \ku^\times$, $q_{jj}=-1$,
$\widetilde{q_{ik}}=\widetilde{q_{il}}=\widetilde{q_{jl}}=1$.

\noindent If $\left[\left[\left[x_{ijk},x_j\right]_c,
\left[x_{ijkl},x_j\right]_c \right]_c, x_{jk} \right]_c
\in\cP(S)\setminus\{0\}$, then $S$ is infinite-dimensional.
\smallskip

\noindent \vii Let $i,j,k,l\in\{1,\ldots,\theta \}$ be such that
$\widetilde{q_{jk}}= \widetilde{q_{ij}}= q_{jj}^{-1}\in
\G_4'\cup\G_6'$, $q_{ii}=q_{kk}=-1$,
$\widetilde{q_{ik}}=\widetilde{q_{il}}=\widetilde{q_{jl}}=1$,
$\widetilde{q_{jk}}^3= \widetilde{q_{lk}}$. If $
\left[\left[x_{ijk},\left[x_{ijkl}, x_k \right]_c\right]_c, x_{jk}
\right]_c \in\cP(S)\setminus\{0\}$, then $S$ is
infinite-dimensional.
\end{lemma}
\pf We use the same notation and strategy as in the previous
Lemmata.

\noindent \vi We have that $\widetilde{Q_{15}}=q^3\neq 1$, and the
diagram corresponding to $(q_{st})$ does not admit extensions with
finite root systems. Therefore $S$ is infinite-dimensional.

\noindent \vii Now, $\widetilde{Q_{15}}=q^3\neq 1$, and we conclude
the same as in \vi. \epf

\medskip

\begin{lemma}\label{lemma:gen-varios casos rango 4}
Let $i,j,k,l\in\{1,\ldots,\theta \}$ be such that one of the
following situations hold:
\begin{itemize}
  \item[$\circ$] $q_{kk}=-1$, $q_{ii}=\widetilde{q_{ij}}^{-1}= q_{jj}^2$,
$\widetilde{q_{kl}}= q_{ll}^{-1}= q_{jj}^3$, $\widetilde{q_{jk}}=
q_{jj}^{-1}$,
$\widetilde{q_{ik}}=\widetilde{q_{il}}=\widetilde{q_{jl}}=1$,
  \item[$\circ$] $q_{ii}=\widetilde{q_{ij}}^{-1}= -q_{ll}^{-1}=-\widetilde{q_{kl}}$,
$q_{jj}=\widetilde{q_{jk}}=q_{kk}=-1$,
$\widetilde{q_{ik}}=\widetilde{q_{il}}=\widetilde{q_{jl}}=1$, or
  \item[$\circ$] $q_{jj}=\widetilde{q_{jk}}^{-1}\in\G_3$, $q_{ii}=\widetilde{q_{ij}}^{-1}=q_{ll}=\widetilde{q_{kl}}^{-1}=-q_{jj}$,
$q_{kk}=-1$,
$\widetilde{q_{ik}}=\widetilde{q_{il}}=\widetilde{q_{jl}}=1$.
\end{itemize}
If $ \left[\left[x_{ijkl}, x_j \right]_c, x_k \right]_c-
q_{jk}(\widetilde{q_{ij}}^{-1}-q_{jj})
\left[\left[x_{ijkl},x_k\right]_c, x_j \right]_c
\in\cP(S)\setminus\{0\}$, then $S$ is infinite-dimensional.
\end{lemma}
\pf Using the same notation and strategy as in the previous Lemmata,
we compute the diagrams for $(Q_{rs})$ in each case, and note that
the diagrams of $(q_{rs})$ does not admit extensions with finite
root systems.
\begin{itemize}
  \item In the first case, $\widetilde{Q_{35}}=q_{jj}\neq 1$.
  \item In the second case, $\widetilde{Q_{25}}=\widetilde{q_{ij}}\neq 1$.
  \item For the last one, $\widetilde{Q_{35}}=-q_{jj}^2\neq 1$.
\end{itemize}
Therefore $S$ is infinite-dimensional in any case. \epf

\medskip

\begin{lemma}\label{lemma:gen-relaciones fila 18, rango3}
\vi Let $i,j,k\in\{1,\ldots,\theta\}$ be such that
$q_{kk}=q_{jj}=\widetilde{q_{ij}}^{-1}=\widetilde{q_{jk}}^{-1}\in
\G_9$, $\widetilde{q_{ik}}=1$, $q_{ii}=q_{kk}^6$. If $ \left[ \left[
x_{iij} , x_{iijk} \right]_c, x_{ij} \right]_c
\in\cP(S)\setminus\{0\}$, then $S$ is infinite-dimensional.
\smallskip

\noindent \vii Let $i,j,k\in\{1,\ldots,\theta\}$ be such that
$q_{ii}=\widetilde{q_{ij}}^{-1}\in \G_9$,
$q_{jj}=\widetilde{q_{jk}}^{-1}=q_{ii}^5$, $\widetilde{q_{ik}}=1$,
$q_{kk}=q_{ii}^6$. If $ [\left[x_{ijk}, x_{j} \right]_c, x_k]_c -(1
+ \widetilde{q_{jk}})^{-1}q_{jk} \left[ \left[x_{ijk}, x_{k}
\right]_c , x_{j} \right]_c \in\cP(S)\setminus\{0\}$, then $S$ is
infinite-dimensional.
\end{lemma}
\pf We use the same notation and strategy as in the previous
Lemmata, and note also that these diagrams does not admit extensions
with finite root systems.

\noindent \vi In this case,
$\widetilde{Q_{34}}=\widetilde{q_{jk}}\neq -1$, so $S$ is
infinite-dimensional, because the diagram of $(Q_{rs})$ is connected
and contains the diagram of $(q_{ij})$.

\noindent \vii Now, $\widetilde{Q_{34}}=q_{jj}\neq -1$, and $S$ is
also infinite-dimensional. \epf

\medskip

\begin{lemma}\label{lema:gen-relaciones restantes rango 3}
\vi Let $i,j,k\in\{1,\ldots,\theta\}$ be such that $q_{ii}=q_{kk}=\widetilde{q_{ij}}=-1$, $q_{jj}=\widetilde{q_{jk}}^{-1}$, $\widetilde{q_{ik}}=1$. If
$\left[x_{ij},x_{ijk}\right]_c\in\cP(S)\setminus\{0\}$, then $S$  is infinite-dimensional.
\smallskip

\noindent\vii Let $i,j,k\in\{1,\ldots,\theta\}$ be such that $q_{ii}=q_{kk}=-1$, $\widetilde{q_{ij}}\in\G_3$, $q_{jj}=\pm\widetilde{q_{ij}}=-\widetilde{q_{jk}}$, $\widetilde{q_{ik}}=1$. If $\left[x_i,x_{jjk}\right]_c-(1+q_{jj}^2)q_{kj}^{-1}\left[x_{ijk},x_j\right]_c-(1+q_{jj})(1+q_{jj}^2)q_{ij}\, x_jx_{ijk}\in\cP(S)\setminus\{0\}$, then $S$  is infinite-dimensional.
\smallskip

\noindent\viii Let $i,j,k\in\{1,\ldots,\theta\}$ be such that $\widetilde{q_{jk}}=1$, $q_{ii}=\widetilde{q_{ij}}=-\widetilde{q_{ik}}\in\G_3$. If
$\left[x_i, \left[ x_{ij},x_{ik} \right]_c \right]_c+q_{jk}q_{ik}q_{ji} \left[ x_{iik} ,x_{ij} \right]_c+q_{ij}\, x_{ij} x_{iik}\in\cP(S)\setminus\{0\}$,
then $S$ is infinite-dimensional.
\smallskip

\noindent\viv Let $i,j,k\in\{1,\ldots,\theta\}$ be such that $q_{jj}=q_{kk}=\widetilde{q_{jk}}=-1$, $q_{ii}=-\widetilde{q_{ij}}\in\G_3$, $\widetilde{q_{ik}}=1$. If $\left[x_{iijk},x_{ijk}\right]_c\in\cP(S)\setminus\{0\}$, then $S$  is infinite-dimensional.
\end{lemma}

\pf
We consider the same notation as before. We consider the subspace $W$ generated by $y_1=x_i$, $y_2=x_j$, $y_3=x_k$ and $y_4$
(the primitive element corresponding to the relation), where $y_i\in W_{h_i}^{\eta_i}$, for some $h_i\in\Gamma$, $\eta_i\in\widehat\Gamma$, and set
$(Q_{rs}=\eta_s(h_r))$. We will prove again that $\cB(W)$ is infinite-dimensional.
\smallskip

\noindent\vi The corresponding diagram of $(Q_{rs})$ is
$$ \xymatrix{\circ^{-1}\ar@{-}[r]^{-1}&\circ^{q}\ar@{-}[d]^{q^3}\ar@{-}[r]^{q^{-1}}&\circ^{-1}\ar@{-}[ld]^{q^{-2}}\\ & \circ^{-q^2}& }, \qquad
q=q_{jj}\in\G_3\cup\G_4\cup\G_6.$$
If $q\in\G_4$ then $Q_{44}=1$, so $\cB(W)$ is infinite-dimensional. If $q\in\G_3\cup\G_6$, the diagram is not in \cite[Table 3]{H-classif RS}, so $\cB(W)$ is also again infinite-dimensional.
\smallskip

\noindent\vii We note that $Q_{14}Q_{41}=Q_{34}Q_{43}=q_{jj}^2\in\G_3$. Therefore, the diagram corresponding to $(Q_{rs})$ contains a 4-cycle, and then $\cB(W)$ is infinite-dimensional.
\smallskip

\noindent\viii In this case, $Q_{24}Q_{42}=-q_{ii}^2$ and $Q_{34}Q_{43}=-q_{kk}^2\neq 1$, because $q_{kk}\in\{-1,\widetilde{q_{ik}}^{-1}\}$. The diagram corresponding to $(Q_{rs})_{r,s=2,3,4}$ is a 3-cycle such that $\widetilde{Q_{24}}\widetilde{Q_{34}}\widetilde{Q_{23}}\neq 1$. \cite[Lemma 9]{H-classif RS} implies that $\cB(W)$ is infinite-dimensional.
\smallskip

\noindent\viv $\cB(W)$ is infinite-dimensional because $Q_{44}=1$.
\epf
\medskip

\begin{lemma}\label{lema:gen-relaciones restantes rango 2}
Let $i,j\in\{1,\ldots,\theta\}$ be such that the satisfy one of the following conditions:

\noindent\vi $-q_{ii}, -q_{jj}, q_{ii}\widetilde{q_{ij}}, q_{jj}\widetilde{q_{ij}} \neq 1$,
$$ \left[x_i, \left[x_{ij}, x_j \right]_c \right]_c - \frac{(1+q_{jj})(1-q_{jj}\widetilde{q_{ij}})}{(1-\widetilde{q_{ij}})q_{jj}q_{ji}}x_{ij}^2\in\cP(S)\setminus\{0\};$$
\noindent\vii $q_{jj}=-1$, $q_{ii}\widetilde{q_{ij}} \notin \G_6$, and also $m_{ij}\in \{4,5\}$, or $m_{ij}=3$, $q_{ii} \in \G_4$,
$$\left[x_i, x_{3\alpha_i+2\alpha_j} \right]_c -\frac{1-q_{ii}\widetilde{q_{ij}}-q_{ii}^2\widetilde{q_{ij}}^2q_{jj}} {(1-q_{ii}\widetilde{q_{ij}})q_{ji}} x_{iij}^2\in\cP(S)\setminus\{0\};$$
\noindent\viii $4\alpha_i+3\alpha_j\notin \Delta_+^\chi$, $q_{jj}=-1$ or $m_{ji}=2$, and also $m_{ij}\geq 3$ or $m_{ij}=2$, $q_{ii}\in\G_3$,
$$[x_{3\alpha_i+2\alpha_j}, x_{ij} ]_c\in\cP(S)\setminus\{0\};$$
\noindent\viv $3\alpha_i+2\alpha_j\in\Delta_+^\chi$, $5\alpha_i+3\alpha_j\notin\Delta_+^\chi$, and $q_{ii}^3\widetilde{q_{ij}}, q_{ii}^4\widetilde{q_{ij}}\neq 1$, $ [x_{iij}, x_{3\alpha_i+2\alpha_j}]_c\in\cP(S)\setminus\{0\}$;
\smallskip

\noindent\vv $4\alpha_i+3\alpha_j\in\Delta_+^\chi$, $5\alpha_i+4\alpha_j\notin\Delta_+^\chi$, $[x_{4\alpha_i+3\alpha_j},x_{ij} ]_c\in\cP(S)\setminus\{0\}$;
\smallskip

\noindent\vvi $5\alpha_i+2\alpha_j\in\Delta_+^\chi$, $7\alpha_i+3\alpha_j \notin \Delta_+^\chi$, $ [[x_{iiij}, x_{iij}], x_{iij} ]_c\in\cP(S)\setminus\{0\}$;
\smallskip

\noindent\vvii $q_{jj}=-1$, $5\alpha_i+4\alpha_j\in\Delta_+^\chi$, $[x_{iij},x_{4\alpha_i+3\alpha_j}]_c- a x_{3\alpha_i+2\alpha_j}^2\in\cP(S)\setminus\{0\}$,
for some $a\in\ku^\times$.
\smallskip

Then $S$ is infinite-dimensional.
\end{lemma}
\pf Firstly we note that there exists just one connected generalized Dynkin diagram of rank three such that $3\alpha_i+2\alpha_j\in\Delta_+^\chi$, for some
pair $i,j$, which is exactly the unique one such that $m_{kl}\geq 3$ for some pair $k,l$. Moreover, $4\alpha_i+3\alpha_j$, $5\alpha_i+4\alpha_j$, are not
positive roots for any pair $i,j$ and any connected Dynkin diagram of rank 3.

We consider as above the subspace $W$ generated by $y_1=x_i$,
$y_2=x_j$ and $y_3$, the relation which is a primitive element by
hypothesis, and analyze its generalized Dynkin diagram.

\noindent\vi If $Q_{13}Q_{31}\neq1$ or $Q_{23}Q_{32}\neq1$, then $\cB(W)$  is infinite-dimensional. In other case,
$$ Q_{13}Q_{31}=q_{ii}^4q_{ij}^2q_{ji}^2=1, \qquad Q_{23}Q_{32}=q_{jj}^4q_{ij}^2q_{ji}^2=1, $$
so $Q_{33}=q_{ii}^4q_{ij}^4q_{ji}^4q_{jj}^4=1$, and $\cB(W)$ is also infinite-dimensional.
\smallskip

\noindent\vii If $q_{ii}\in\G_4$, $\widetilde{q_{ij}}=q_{ii}=q_{jj}^{-1}$ (and then $(q_{rs})$ is Cartan of type $G_2$), then
$$Q_{33}=q_{ii}^{16}q_{ij}^8q_{ji}^8q_{jj}^4=1,$$
so $\cB(W)$  is infinite-dimensional. In other case, $Q_{13}Q_{31}\neq1$, or $Q_{23}Q_{32}\neq1$, or
$$ Q_{13}Q_{31}=q_{ii}^8q_{ij}^2q_{ji}^2=1, \ Q_{23}Q_{32}=q_{ij}^4q_{ji}^4q_{jj}^4, \qquad \mbox{so }Q_{33}=1, $$
and therefore $\cB(W)$  is infinite-dimensional.
\smallskip

\noindent\viii Now we calculate
$$ Q_{33}=q_{ii}^{16}q_{ij}^{12}q_{ji}^{12}q_{jj}^9, \qquad Q_{13}Q_{31}=q_{ii}^8q_{ij}^3q_{ji}^3, \qquad Q_{23}Q_{32}=q_{ij}^4q_{ji}^4q_{jj}^6. $$
If $(q_{rs})$ is Cartan of type $G_2$ and $q_{ii}\in\G_6$, $\widetilde{q_{ij}}=q_{jj}=-1$, then $\cB(W)$  is infinite-dimensional, because we have a connected
diagram of rank three such that $M_{12}=3$, and it is not of type super $G(3)$. In other case, we will prove that $Q_{13}Q_{31}\neq1$ or $Q_{23}Q_{32}\neq1$
to conclude that $\cB(W)$  is infinite-dimensional. If $m_{ji}\geq2$, we have the following possible cases:
\begin{itemize}
  \item $q_{ii}=-\zeta$, $\widetilde{q_{ij}}=\zeta^7$, $q_{jj}=\zeta^3$, $\zeta\in\G_9$; in such case, $Q_{23}Q_{32}=\zeta$.
  \item $q_{ii}=-\zeta$, $\widetilde{q_{ij}}=-\zeta^{12}$, $q_{jj}=\zeta^5$, $\zeta\in\G_{15}$; therefore, $Q_{23}Q_{32}=\zeta^3$.
\end{itemize}
Also, if $q_{ii}=\zeta^8$, $\widetilde{q_{ij}}=\zeta^3$, $q_{jj}=-1$, $\zeta\in\G_{12}$, then $Q_{13}Q_{31}=\zeta$. In all the remaining cases, $q_{jj}=-1$ and
$\widetilde{q_{ij}}\notin\G_4$, so $Q_{23}Q_{32}\neq1$.
\smallskip

\noindent\viv This relation is not redundant just in the following two cases:
$$ \xymatrix{\circ^{\zeta^3}\ar@{-}[r]^{\zeta^8}&\circ^{-1}},\ \zeta\in\G_9,\qquad\xymatrix{\circ^{\eta^3}\ar@{-}[r]^{-\eta^4}& \circ^{-\eta^{-4}}}, \ \eta\in\G_{15}.$$
Note that they are not contained in any connected diagram of rank three in \cite[Table 2]{H-classif RS}, so it is enough to verify that $Q_{13}Q_{31}\neq1$
or $Q_{23}Q_{32}\neq1$ to conclude that $\cB(W)$ is infinite-dimensional. For the first diagram, $Q_{23}Q_{32}=\zeta^4\neq 1$; and for the second one,
$Q_{23}Q_{32}=-\eta^{-4}\neq1$.
\smallskip

\noindent\vv The proof is analogous to \viii. Note that
$$ Q_{33}=q_{ii}^{25}q_{ij}^{12}q_{ji}^{20}q_{jj}^{16},\qquad Q_{13}Q_{31}=q_{ii}^{10}q_{ij}^4q_{ji}^4, \qquad Q_{23}Q_{32}=q_{ij}^5q_{ji}^5q_{jj}^8. $$
We have that $q_{jj}=-1$ for every diagram satisfying the conditions for this item. Also, if $q_{ii}=\zeta\in\G_5$, $\widetilde{q_{ij}}=\zeta^2$, it follows that
$Q_{33}=1$. In the remaining cases, $\widetilde{q_{ij}}\notin\G_5$, so $Q_{32}Q_{23}\neq-1$, and then $\cB(W)$ is infinite-dimensional.
\smallskip

\noindent\vvi In this case,
$$ Q_{33}=q_{ii}^{49}q_{ij}^{21}q_{ji}^{21}q_{jj}^9,\qquad Q_{13}Q_{31}=q_{ii}^{14}q_{ij}^3q_{ji}^3, \qquad Q_{23}Q_{32}=q_{ij}^7q_{ji}^7q_{jj}^6. $$
We have that $q_{jj}=-1$ for every diagram satisfying the conditions for this item, and also $\widetilde{q_{ij}}\notin\G_7$, so $Q_{23}Q_{32}\neq1$. Therefore
$\cB(W)$ is infinite-dimensional.
\smallskip

\noindent\vvii The proof is analogous to the one for \vi.
\smallskip

We conclude that $S$  is infinite-dimensional in all the cases.
\epf
\bigskip

Now we can prove the main Theorem of this Section.

\begin{theorem}\label{thm:generationdegree}
Let $S=\oplus_{n\geq0}S(n)$ be a finite-dimensional graded Hopf algebra in $^{\ku \Gamma}_{\ku \Gamma} \mathcal{YD}$, where $\Gamma$ is a finite abelian group, such that $S(0)=\ku 1$. If $S$ is generated as an algebra by $S(0)\oplus S(1)$, then $S \cong \cB(V)$.
\end{theorem}
\pf
Fix a basis $x_1,\ldots,x_{\theta}$ of $V:=S(1)$, such that $x_i \in S(1)^{\chi_i}_{g_i}$ for some $g_i \in \Gamma$ and $\chi_i \in \widehat{\Gamma}$, and set $q_{ij}:=\chi_j(g_i)$.

As $S$ is generated as an algebra by $S(0)\oplus S(1)$, the canonical projection $T(V)\twoheadrightarrow\cB(V)=T(V)/I(V)$ induces a surjective morphism $\pi: S \twoheadrightarrow \cB(V)$ of graded braided Hopf algebras; we can consider $S=T(V)/I$, for some graded braided Hopf ideal $I$ of $T(V)$, generated by homogeneous elements of degree $\geq 2$, $I \subseteq I(V)$.

Suppose that $I(V)\supsetneqq I$. Then at least one of the generators of $I(V)$ from Theorem \ref{thm:presentacion minima} does not belong to $I$. We can assume that $\xb\in I(V)\setminus I$ is one of these generators, of minimal degree $k$. Then $\xb$ is primitive in $S$ by Lemma \ref{Lema:minimal es primitivo}.

By Proposition \ref{prop:gen-relacion qSerre} and Lemmata \ref{Lema:gen-relacion vertice -1}-\ref{lema:gen-relaciones restantes rango 2}, we deduce that $\xb=x_\alpha^{N_{\alpha}}$ for some $\alpha\in\cO$, or a simple root $\alpha=\alpha_i$ such that $i$ is not a Cartan vertex, or $\alpha=\alpha_i+\alpha_j$, such that $N_\alpha=2$, $q_{ii}=q_{jj}=\widetilde{q_{ij}}=-1$. If $g_\alpha\in\Gamma$, $\chi_\alpha\in\widehat\Gamma$ are the associated elements, we have that $q_\alpha=\chi_\alpha(g_\alpha)$, which is a root of unity of order $N_\alpha$. Therefore $g_\alpha^{N_\alpha}\in\Gamma$ and $\chi_\alpha^{N_\alpha}\in\widehat\Gamma$ are the associated elements to $\xb$, and
$$ c(\xb\ot\xb)= g_\alpha^{N_\alpha}\cdot\xb \ot\xb= \chi_\alpha^{N_\alpha}\left(g_\alpha^{N_\alpha}\right)\xb \ot\xb=\xb\ot\xb, $$
so $\xb$ generates in $S$ an infinite-dimensional braided Hopf subalgebra, and we obtain a contradiction. In consequence, $S\cong \cB(V)$.
\epf

\begin{remark}
Note that we just use the fact that the braiding is diagonal, so we can generalize this Theorem to a general braided Hopf algebra $R$ in $\yd$, where $H$ is a finite-dimensional Hopf algebra which acts diagonally over $R(1)$.
\end{remark}

The following Theorem answers positively Conjecture \ref{conj:AS, generacion en grado 1} in the case that the group of group-like elements is abelian. It extends
\cite[Thm. 5.5]{AS Class}.

\begin{theorem}\label{thm:generationdegree1}
Let $H$ be a finite-dimensional pointed Hopf algebra over an abelian group $\Gamma$. Then $H$ is generated by its group-like and skew-primitive elements.
\end{theorem}
\pf
Let $\gr H=R\#\ku \Gamma$, $V=R(1)$. Then $H$ is generated by its group-like and skew-primitive elements if and only if $\gr H$ satisfies this condition, which is equivalent to the fact that $R$ is the Nichols algebra $\cB(V)$. Let $S$ be the graded dual $R^*$ in the category $\gyd$, which is generated as an algebra by $S(1)=V^*$. By \cite[Lemma 2.3]{AS Pointed HA} it is enough to prove that $S$ is the Nichols algebra $\cB(V^*)$, which follows by Theorem \ref{thm:generationdegree}.
\epf

\end{document}